\documentclass{article}
\usepackage{amsmath,amsthm,amsfonts,amscd,amssymb,eucal,latexsym,hyperref} 
\usepackage[all]{xy}
\def\cA{{\cal A}}

\def\uu{\underline u}
\def\uv{\underline v}

\def\cH{{\cal H}}

\def\cC{{\cal C}}

\begin{document}

\title{A theory of induction and classification of tensor $C^*$--categories}
\author{Claudia Pinzari$^1$, John E. Roberts$^2$\\ \\
$^1\,$Dipartimento di Matematica, Universit\`a di Roma ``La Sapienza''\\
00185--Roma, Italy\\
$^2\,$Dipartimento di Matematica, Universit\`a di Roma ``Tor Vergata'',\\
00133 Roma, Italy}
\date{}
\maketitle
                    
\begin{abstract}

 This paper addresses the problem of describing the structure  of tensor 
$C^*$--categories ${\cal M}$ with conjugates  and irreducible tensor unit. 
No assumption on 
the existence of a   braided symmetry or on amenability is made.  Our 
assumptions are motivated by the remark that these  
categories  often contain non-full tensor $C^*$--subcategories  with 
conjugates and the same objects admitting an embedding 
into the Hilbert spaces. Such an embedding defines a compact 
quantum group by 
Woronowicz duality. An important example 
 is the Temperley--Lieb category canonically contained in a tensor 
$C^*$--category generated by a single real or pseudoreal object of 
dimension $\geq2$. The associated quantum groups are the universal 
orthogonal quantum groups of Wang and Van Daele.
  
Our main result asserts that
there is a full and faithful tensor functor 
from ${\cal M}$ to a category of Hilbert bimodule 
representations of the   compact quantum group.
In the classical case,  these bimodule representations reduce to the
  $G$--equivariant 
Hermitian bundles over compact homogeneous $G$--spaces, with $G$ a 
compact group.  
Our structural results shed light on   the problem of whether there is
an embedding functor of ${\cal M}$ into the  Hilbert 
spaces. We show  that this is related to
the problem of whether a classical compact Lie group
can act ergodically on a non-type $I$ von Neumann algebra. 
In particular, combining this with a result of Wassermann shows that an 
embedding exists if ${\cal M}$ is generated by  a pseudoreal  
object of dimension $2$.

\end{abstract}

\tableofcontents   

\section{Introduction}  

By Tannaka--Krein duality theory, a semisimple   rigid tensor category
admits an   embedding functor  into the 
category of  finite dimensional vector spaces if and only if it is the  representation category of a quantum 
group.   
There are variants of this result depending on the 
framework under consideration.
This characterization does  however not allow one to tell whether a given tensor category 
admits such an embedding.

A positive result, motivated by algebraic QFT,  asserts that
a symmetric tensor $C^*$--category with conjugates $\cal M$, 
is equivalent, after completion under direct sums
and subobjects, to the symmetric tensor $C^*$--category of finite-dimensional
representations of a unique compact group \cite{DR}. There is a similar 
well known result in the context of
algebraic geometry \cite{Deligne}.

However,  different finite groups may have equivalent
representation categories   \cite{EG, IK}, so the symmetry is crucial for uniqueness.

Another positive result  is  a 
theorem of  \cite{Hayashi}, asserting that a 
semisimple rigid tensor category with
finitely many irreducibles is  equivalent to a representation category of  
a finite dimensional, but not unique, weak Hopf algebra in the sense of 
\cite{Szlachanyi2},  see also \cite{Ostrik}.
However, this theorem does not say anything about whether the weak Hopf algebras
can be chosen to be a quantum group or even a group. Moreover, this approach does not generalize 
easily to categories with infinitely many 
irreducibles, as difficulties of an analytic nature arise.

In this paper we consider  the analytic framework of tensor $C^*$--categories with conjugates. Now, an embedding functor is naturally required to take values in the category $\text{Hilb}$ of 
Hilbert spaces. Here, the quantum groups of the Tannaka--Krein duality theory
are compact quantum groups \cite{WoronowiczTK}. However, even in this case, there 
are situations where there can be no such embedding.

 For example, if a tensor $C^*$-category 
with conjugates $\cal M$ has an object with an intrinsic dimension strictly 
between 1 and 2, it obviously cannot be embedded into Hilb. 
These cases are often related to quantum  groups at roots of unity 
by Jones fundamental result on the 
restriction of the index values
\cite{Jones1, Wenzl}.  Ocneanu  
\cite{O} indicated that they should be understood as `quantum subgroups' of the deformed $SU(2)$, a program developed in  \cite{KO}.

But other classes are known even when there are objects with intrinsic dimension $>2$. 
For example,   consider    
an irreducible  inclusion of $II_1$ factors $N\subset M$ with finite 
Jones index $[M:N]$ and consider
the tensor $C^*$--category    of $N$--bimodules 
generated by ${}_NM_N$. This category has conjugates and the intrinsic 
dimension of the generator ${}_NM_N$ is  $[M:N]$, hence 
$\geq2$ if the inclusion is proper, again by the restriction of the index 
values.   By an easy consequence of Popa's work \cite{Popa}, this category 
cannot 
be embedded in Hilb  whenever 
the index is not an integer and the inclusion is amenable in the sense of Popa,
 see \cite{PR2}. There is a similar result  for an amenable object in a tensor
$C^*$--category with non-integral dimension \cite{LR}. 
(The precise relation between the two notions of amenability has not been clarified.)

On the other hand, compact quantum groups provide examples of 
tensor $C^*$--categories with conjugates that are embedded but not 
amenable, as
 the intrinsic dimensions of unitary representations are often not integral.

 This paper addresses the problem of describing the structure  of tensor $C^*$--categories with conjugates  and irreducible tensor unit. We are interested  in the case where  there is a a generator  with intrinsic dimension $\geq 2$. 
No assumption on 
the existence of a   braided symmetry or on amenability is made.

  The tensor $C^*$--categories ${\cal M}$ arising from subfactors have as 
  objects the tensor powers of an irreducible selfconjugate object $x$, a property expressed in 
terms of an intertwiner $R\in(\iota,x\otimes x)$. 
The tensor $^*$--subcategory generated by $R$ is a Temperley-Lieb category, 
admitting an embedding if $d(x)\geq2$. All such embeddings may be  classified; they correspond to
suitable compact quantum groups $A_o(F)$ of Wang and Van Daele. In particular,
if  $F\in M_2({\mathbb C})$ these quantum groups are the quantum $SU(2)$ groups of Woronowicz 
for  deformation parameters uniquely determined by the 
dimension and the reality character  of $x$.
A similar result holds 
if the objects of $\cal M$ are the semigroup with unit generated by an object 
$x$   and its conjugate $\overline x$,
the quantum groups involved are certain $A_u(F)$.

These remarks show that although tensor $C^*$--categories with conjugates cannot be embedded
 generically, they may contain an embeddable subcategory with conjugates.  

Abstracting from the above, 
we start from two tensor $C^*$--categories with conjugates 
$\cal A$ and $\cal M$, a tensor $^*$--functor $\tau:\cA\to$ Hilb and a quasitensor functor 
$\mu:\cal A\to\cal M$. We may suppose that the objects of $\cal M$ are tensor products 
of objects in the image of $\mu$. By \cite{WoronowiczTK}, 
$\tau$ determines a compact quantum group $G_\tau$.  We showed in \cite{PR} that the pair $\mu,\tau$ 
canonically defines an ergodic action of $G_\tau$ on a $C^*$--algebra $\cal C$.
When $\mu$ is the functor  restricting a representation to a subgroup $K$ of $G_\tau$, this action reduces to the translation action on the quantum quotient space $K\backslash G_\tau$.

If $\cal M$ is not embeddable, the associated ergodic action cannot correspond to a 
true quantum subgroup \cite{Pinzari3}. 
Borrowing a notion due to Mackey \cite{MackeyVirtual}, 
we  may talk of a virtual subgroup.

The notions quasitensor functor and relaxed tensor functor are recalled in Sect.\ 2, 
this extra generality is motivated by their role in the duality theorem for 
ergodic actions of compact quantum groups, where they arise as spectral functors. (The spectral functor of an ergodic action  maps a representation to the corresponding spectral space, thus, in particular, it takes values in the category of Hilbert spaces \cite{PR}.)  Note that, unlike relaxed tensor functors, quasitensor functors may take a non-zero object to the zero object.

One of the aims of this paper is to describe ${\cal M}$ as  a 
category of representations of the virtual subgroup.
To handle the non-embeddable case, 
we introduce the notion of a representation of a  compact quantum group on a Hilbert bimodule over an ergodic $C^*$--algebra. This is the
noncommutative analogue of the bimodule of continuous sections of an equivariant 
Hermitian bundle over a compact homogeneous space. 
  We show that   these bimodule representations form a tensor $C^*$--category with irreducible tensor unit, Theorem 3.1. 

 Given an object of ${\cal M}$, we construct a Hilbert bimodule representation
 of $G_\tau$, that we regard as a representation induced from the virtual subgroup, as
 the associated bimodule generalizes the bimodule of sections of the equivariant vector 
bundle induced from a representation of a subgroup.
As a right module, the induced module  turns out to be finite 
projective and, if $\mu$ is tensorial, even free and 
finitely generated, a result generalizing Swan's theorem to a noncommutative framework.

We show that the bimodule  construction yields 
 a full and faithful  tensor $^*$--functor, the induction functor, from $\cal M$ 
to the category $\text{Bimod}_\alpha(G_\tau)$ of Hilbert bimodule representations
of $G_\tau$. Note that this tensoriality property holds despite the fact that
$\mu$ was only assumed to be quasitensor, Theorems 6.2 and  6.4. 

In particular, if ${\cal M}$ is generated as a tensor $C^*$--category by a real or pseudoreal object $x$ of $d(x)\geq2$, then it may be identified with a category
of bimodule representations for any one of  a class of compact quantum groups $A_o(F)$. 
A similar result holds if $x$ is not selfconjugate, Theorems 6.5 and 6.6.

We then use these abstract results to investigate the case where ${\cal M}$ is an extension
of the representation category of a compact Lie group $G$. We show that if the associated
ergodic $C^*$--algebra ${\cal C}$ yields a finite type $I$ von Neumann algebra after completion 
in the GNS-representation  of the unique invariant trace, then 
  $\cal M$ admits an embedding functor. In fact, we construct an explicit full tensor functor  to the representation category of a closed subgroup of $G$, Theorem 6.7. 

This last result shows that the question of the existence of an embedding into the Hilbert spaces is related to the open problem posed in  \cite{HLS} and mentioned in the abstract on the existence of ergodic actions of classical compact Lie groups on non-type-$I$ von Neumann algebras.  A negative answer for   $G$ would imply  that the associated ergodic von Neumann algebra is of finite type $I$ and hence that  ${\cal M}$ is embeddable.

This  is known  for $SU(2)$  \cite{Wassermann3}. This negative result shows that 
if the  objects of ${\cal M}$  are generated by a single pseudoreal object of intrinsic dimension $2$,
${\cal M}$ can be embedded, and admits a full and faithful tensor $^*$--functor to the category of 
representations of a closed subgroup  of $SU(2)$, Theorem 6.9.
We would like to point out 
  the analogy of this result with the well known classification of subfactors of index $4$ in terms of closed subgroups of $SU(2)$. 
  
  The notion of full bimodule representation plays a role in  this paper.  An object 
of Bimod$_\alpha(G)$ is full if every fixed vector for the action is central
(see Sect.\ 5). This guarantees that the left module structure is 
naturally compatible with that of right module representation. We show that the induced 
bimodule representations are full (Prop.\ 8.9, 9.4) and
use this to show the embedding result. Furhermore we use this property 
to show that certain ergodic actions cannot arise from  
a pair of  {\it tensor} functors $\mu$, $\tau$.
We show for example  that  
neither the adjoint action of a non-trivial irreducible representation of $SU(2)$
nor those with full spectrum and low multiplicity can arise, Sect.\ 11. 

Here we have interpreted induction in terms of bimodule representations. 
However we may also induce from representations of the virtual subgroup 
to Hilbert space representations of $G_\tau$, except that such representations need no 
longer be finite dimensional. However after completing our categories 
under infinite direct sums, we show that induction and restriction are a 
pair of adjoint functors, Theorem\ 10.1. 

In conclusion, we draw the reader's attention to an incomplete list of papers where related results may be found,  although some from different perspectives, 
\cite{BS,   Szlachanyi2,   Cuntz,    DPR, ENO, Izumi2,  LongoHopf, Vainerman1, Vainerman2, Pal,  Pinzari3,  Szymanski,   Vaes}.

The paper is organized as follows. Sect.\ 2 establishes notation and recalls results that we shall need. In Sect.\ 3, we explain the notion of a representation of a compact quantum group on a Hilbert $C^*$--bimodule and we introduce the tensor $C^*$--category $\text{Bimod}_\alpha(G)$.  In Sect.\ 4 
we review Mackey's induced representation and Frobenius reciprocity 
from the standpoint of bimodule representations.
In the next section, we introduce the notion of   full
bimodule representation and discuss the example of quantum quotients. 
In Sect.\ 6 we illustrate the  main ideas and results of 
this paper. Sections 7-9 are dedicated to the induction functor into 
the tensor category of Hilbert $C^*$--bimodules. In Sect.\ 7 we 
give the algebraic construction of bimodule and introduce an inner product  starting from a pair 
$(\tau, \mu)$ and show positivity of the inner product when $\tau$ is tensorial in sect.\ 8,
leading to the Hilbert $C^*$--bimodule representation of the compact quantum group $G_\tau$ associated with $\tau$. 
In Sect.\ 9 we show that there is a unique extension of the the induction functor to a  
 tensor functor.
 In Sect.\ 10 we show that if we instead define an induction functor taking 
values in the category of unitary representations of $G_\tau$ on Hilbert spaces, then 
$(\mu,\text{Ind})$ is an adjoint pair of functors.
Sect.\ 11 and 12 are dedicated to the analysis of ergodic actions of compact groups.
In Sect.\ 11 we classify full bimodule representations of compact groups on finite type $I$ von Neumann 
algebras and use the classification in the following section to derive results on the 
problem of embedding into Hilbert spaces. A few computations 
in an appendix conclude the paper.

\section{Notation and Preliminaries}

\noindent {\it 2.1. Tensor $C^*$--categories with conjugates}
We shall work with tensor $C^*$--categories defined as in \cite{DR}. By MacLane's theorem  \cite{MacLane}, we may and shall assume that the tensor product is strictly associative. 
The tensor product between objects
$u$ and
$v$ will be  denoted by $u\otimes v$
and between
arrows $S$ and $T$ by $S\otimes T$.
The $n$--th tensor power of an object $u$ will be denoted $u^n$. 
The tensor unit, denoted by $\iota$, will always be
assumed irreducible:
$$(\iota,\iota)={\mathbb C}.$$
An object $\overline{u}$ of a tensor $C^*$--category ${\cal A}$ is a conjugate of $u$ if
there are arrows
$$R\in(\iota,\overline{u}\otimes u),\quad \overline{R}\in(\iota, u\otimes\overline{u})$$
satisfying the conjugate equations
$$\overline{R}^*\otimes 1_u\circ 1_u\otimes R=1_u,\quad {R}^*\otimes 1_{\overline{u}}\circ 1_{\overline{u}}\otimes \overline{R}=1_{\overline{u}}.$$
If $R$, $\overline{R}$ is a solution of the conjugate equations for 
$u$, any other solution is of the form ${X^*}^{-1}\otimes 1_u\circ R$ and $1_u\otimes X\circ 
\overline{R}$, where $X\in(\overline{u},\tilde{u})$ is an invertible intertwiner. $\overline R$ is 
uniquely determined by $R$. 

We will always take $1_\iota$ as the solution of the conjugate equations 
for $\iota$.

An object $u$ is called {\it real} or {\it pseudoreal} if we may choose
$\overline{u}=u$ (i.e. $u$ selfconjugate), and   a solution of the form $\overline{R}=R$ or $\overline{R}=-R$ respectively.

We shall say that ${\cal A}$ has conjugates if every object has a conjugate.
 In this case, every object is a direct sum of minimal projections.
 A  solution of the conjugate equations $(R,\overline{R})$  is said to be {\it standard} if 
$$R^*\circ 1_{\overline u}\otimes Y\circ R=\overline R^*\circ Y\otimes 1_{\overline u}\circ\overline R,\quad Y\in(u,u).$$
The $X\in(\overline u,\tilde u)$ taking one standard solutions to another is unitary. The intrinsic dimension of an object is defined as $d(u)=\|R\|^2$, where $R$ is part of a standard solution. Equivalently, $d(u)$ is the minimal value  of $\|R\|\|\overline{R}\|$ for all solutions.
We refer to \cite{LR} for details.

Fix objects $u$, $v$ of ${\cal A}$ and 
pick a solution $R_u$, $\overline{R}_u$ and  $R_v$, $\overline{R}_v$ of the conjugate equations for $u$ and $v$ respectively, and define the associated antilinear map,
$$A\in(v,u)\to
A^\bullet:=R_v^*\otimes 1_{\overline{u}}\circ 1_{\overline{v}}\otimes 
A^*\otimes 1_{\overline{u}}\circ1_{\overline{v}}\otimes\overline{R}_u,
\quad \in(\overline{v},\overline{u}).$$
This map depends on the choice of conjugates: changing the 
solution of the conjugate equations using invertibles 
$X\in(\overline{u},\tilde{u})$ and  $Y\in(\overline{v},\tilde{v})$, $A^\bullet$ becomes  $X\circ A^\bullet \circ Y^{-1}$.

We stress that the notation $R_u$ refers to a particular solution
of the conjugate equations for $u$ but does not necessarily
imply a choice for each object $u$.

\medskip

\noindent{\it Example} In the category $\text{Hilb}$ any solution of the conjugate equations
for a f.d. Hilbert space $H$ is of the form
$${R}=\sum_h\phi_h\otimes j^{-1}\phi_h,\quad
{\overline{R}}=\sum_k\psi_k\otimes j\psi_k,$$
where $j$ is an antilinear invertible map to another Hilbert space $\overline{H}$
and $(\psi_k)$ and $(\phi_h)$ are orthonormal bases 
of $H$ and $\overline{H}$ respectively. We shall use the notation $j_H$ to emphasize
that $j$ refers to the object $H$.
For $A\in(H,K)$, the associated $\bullet$  is given by
$A^\bullet=j_KAj_H^{-1}$. In particular, $\psi^\bullet=j_H\psi$ for $\psi\in H$.
\medskip

If $(R_u,\overline{R}_u)$, $(R_v,\overline{R}_v)$ are solutions for $u$ and $v$ respectively, $R_{u\otimes v}:=1_{\overline{v}}\otimes R_u\otimes 1_v\circ R_v$ and 
$\overline{R}_{u\otimes v}:=1_u\otimes\overline{R}_v\otimes 
1_{\overline{u}}\circ \overline{R}_u$ is a solution for $u\otimes v$, 
called the {\it tensor product solution}.
Similarly, $R_{\overline{u}}:=\overline{R}_u$ and 
$\overline{R}_{\overline{u}}:=R_u$ the solution for $\overline{u}$, called the
{\it conjugate solution}. 

The  main properties of $\bullet$ 
are the following:
$$(A\circ B)^\bullet=A^\bullet\circ B^\bullet,$$ 
$$(A\otimes B)^\bullet=B^\bullet\otimes A^\bullet,$$
for the tensor product solution.

\medskip

\noindent{\it 2.2. Quasitensor functors}
Although all tensor categories may be assumed to be strict, it is well known that one may meet functors preserving the tensor structure only up to a natural equivalence. Here we need the notion of {\it quasitensor} functor whose definition we recall. 
A $^*$--functor $\mu:{\cal A}\to{\cal M}$ is called quasitensor if  there are   isometries
$\tilde\mu_{u,v}\in(\mu_u\otimes\mu_v, \mu_{u\otimes v})$,
 such that
$${\mu}_\iota=\iota,\eqno(2.1)$$
$$\tilde\mu_{u,\iota}=\tilde\mu_{\iota, u}=1_{\mu_u},\eqno(2.2)$$
$$\tilde\mu_{u,v\otimes w}^*\circ\tilde\mu_{u\otimes
v,w}=1_{\mu_u}\otimes\tilde\mu_{v,w}
\circ{\tilde\mu_{u,v}}^*\otimes 1_{\mu_w}\eqno(2.3)$$
and natural in $u$, $v$, 
$$\mu({S\otimes T})\circ
\tilde\mu_{u,v}=\tilde\mu_{u',v'}\circ\mu(S)\otimes\mu(T),\eqno(2.4)$$
 for  objects $u$, $v$, $w$, $u'$, $v'$ of ${\cal A}$ and arrows
$S\in(u, u')$,
$T\in(v,v')$.
The above definition was given in \cite{PR}
in a different form, in connection with the study of ergodic actions of compact quantum groups on unital $C^*$--algebras. The equivalence with the above definition was shown in \cite{Ergodic}. If all the isometries $\tilde\mu_{u,v}$ are unitary,
$(\mu,\tilde{\mu})$ will be called a {\it relaxed tensor} functor.  In particular, a
{\it strict tensor} functor, or simply a tensor functor, is a quasitensor functor with $\tilde\mu_{u,v}:=1_{\mu_u\otimes\mu_v}$ for all objects $u$, $v$.
Note that a 
quasitensor functor may take a non-zero object to the zero object. Examples arise from ergodic actions  (cf. also subsect.\ 2.4).

Notice that, since we are dealing with isometries, $(2.3)$ implies the associativity property,
$$\tilde{\mu}_{u\otimes v, w}\circ\tilde{\mu}_{u,v}\otimes 1_{\mu_w}=\tilde{\mu}_{u, v\otimes w}\circ 1_{\mu_u}\otimes\tilde{\mu}_{v,w}.$$
Hence both sides of this equation  define the same intertwiner $\tilde{\mu}_{u,v,w}\in(\mu_{u}\otimes\mu_v\otimes\mu_w,\mu_{u\otimes v\otimes w})$. Iterating, we get, for any finite sequence $\uu=(u_1,\dots,u_n)$ of objects of ${\cal A}$, with $n\geq2$, an unambiguous arrow 
$$\tilde{\mu}_{u_1,\dots, u_n}\in(\mu_{u_1}\otimes\dots\otimes\mu_{u_n},\mu_{u_1\otimes\dots\otimes u_n}).$$
We set $\tilde{\mu}_{u}=1_{\mu_u}$ for a sequence of length $1$.
$\tilde{\mu}$
is a natural transformation, i.e. for $S_i\in(u_i,v_i)$,
$$\mu(S_1\otimes S_2\otimes\cdots S_n)\circ\tilde\mu_{u_1,u_2,\dots,u_n}=\tilde\mu_{v_1,v_2,\dots,v_n}\circ\mu(S_1)\otimes\mu(S_2)\otimes\cdots\otimes\mu(S_n).$$

We remark however that $(2.3)$ is stronger than associativity. It also implies a categorical analogue of Popa's
commuting square condition in the theory of subfactors, see  \cite{PR} and references there. If $\cal A$ has conjugates, so do objects in the image of $\mu$. 
In detail, if $(\mu,\tilde{\mu}):{\cal A}\to {\cal M}$ is a quasitensor functor
and if $(R,\overline{R})$ is a solution of the conjugate equations for $u$ and $\overline{u}$  in ${\cal A}$ then
$\hat{R}:=\tilde{\mu}_{\overline{u}, u}^*\circ\mu(R)$, $\hat{\overline{R}}:=\tilde{\mu}_{\overline{u}, u}^*\circ\mu(\overline{R})$
 is a solution of the conjugate equations for $\mu_u$ in ${\cal M}$, called the {\it image solution}. In particular,
$d(\mu_u)\leq d(u)$.

 Image solutions $\hat{R}_u$, $\hat{\overline{R}_u}$ and 
  $\hat{R}_v$, $\hat{\overline{R}_v}$ associated to $u$ and $v$ define a map $\bullet$ on the the arrow space $(\mu_u,\mu_v)$ in ${\cal M}$ and we have
$$\mu(A)^\bullet=\mu(A^\bullet).$$

If $(\mu,\tilde{\mu}):{\cal A}\to{\cal M}$ and $(\nu, \tilde{\nu}):{\cal M}\to{\cal N}$ are 
quasitensor functors, the composition $\nu\mu$ becomes   quasitensor with natural transformation  $\nu(\tilde{\mu}_{u,v})\circ\tilde{\nu}_{\mu_u,\mu_v}$ \cite{Ergodic}. 
A composition  $\nu\mu$ of two quasitensor functors will always be implicitly understood
as a quasitensor functor  with this natural transformation.

\medskip

\noindent{\it 2.3.  Ergodic $C^*$--actions of compact groups.}
Let   $\alpha:G\to\text{Aut}({\cal C})$ 
be a continuous ergodic action of  a compact group $G$ on a unital $C^*$--algebra ${\cal C}$. 
The   {\it finiteness theorem}  for the noncommutative ergodic space 
${\cal C}$ and the {\it 
multiplicity bound theorem}  assert respectively  that the unique $G$--invariant state of ${\cal C}$ is a trace,  and that the   multiplicity    of an irreducible representation of $G$  in
 $\alpha$ is bounded above by its dimension. Furthermore, any von Neumann algebra with an ergodic action of a compact group is necessarily hyperfinite \cite{HLS}.
 
 Recall that if $\beta$ is an automorphic action of a closed subgroup $K$ of $G$ on a von Neumann algebra 
${\cal F}$, the induced von Neumann algebra is defined by:
$$\text{Ind}({\cal F}):=\{f\in L^\infty(G, {\cal F}): f(kg)=\beta_k(f(g)), k\in K, g\in G\}=$$
$$(L^\infty(G)\overline{\otimes} {\cal F})^{\lambda\otimes\beta},$$
where $\lambda$ is left translation of $K$ on $L^\infty(G)$. 
If ${\cal F}$ is a $C^*$--algebra,  
the von Neumann tensor product
$\overline{\otimes}$ is replaced by the minimal one, and $L^\infty$-functions by 
continuous ones. 
The induced algebra carries
the induced action $\rho$ of $G$ given by right translation.

As recalled in \cite{Wassermann1}, combining the above results with an imprimitivity theorem of  Takesaki  \cite{Takesaki} for locally compact group actions on von Neumann algebras allows one to reduce the study of ergodic actions on von Neumann algebras to those on finite factors.
Indeed, any ergodic   action of a compact group 
$G$ on a von Neumann
algebra ${\cal C}$ is induced by an action of a closed subgroup $K$ 
 on a matrix 
algebra or on the hyperfinite $II_1$ factor $R$.

Wassermann has shown the important result that  $G=SU(2)$ 
acts ergodically   only on (finite) type $I$ von Neumann algebras \cite{Wassermann3}.
 For more results in this direction see also \cite{Ocairbre}.  
It is not yet known whether there are any ergodic actions of  
compact classical Lie groups on the hyperfinite $II_1$ factor $R$, a problem  raised in \cite{HLS}.
\medskip

\noindent{\it 2.4.  Ergodic $C^*$--actions of compact quantum groups.} We refer to \cite{Wlh} for the general definition of a compact quantum group.
 If $G=({\cal Q},\Delta)$ is a compact quantum group, $\text{Rep}(G)$ will denote the category of unitary finite dimensional representations of $G$.

The theory of ergodic actions of compact quantum groups  on unital $C^*$--algebras has been initiated in \cite{Boca, Podles}. Recall from \cite{Podles} that an action 
of $G$ 
on a unital $C^*$--algebra ${\cal C}$ is a unital $^*$--homomorphism
$$\alpha:{\cal C}\to{\cal C}\otimes {\cal Q},$$
where $\otimes$ denotes the minimal tensor product of $C^*$--algebras,
such that $\alpha\otimes \iota\circ\alpha=\iota\otimes\Delta\circ\alpha$ and with the property that 
$\alpha({\cal C}){\mathbb  C}\otimes Q$ is dense in ${\cal C}\otimes Q$.
The action is called {\it ergodic} if ${\cal C}^\alpha:=\{c\in{\cal C}:\alpha(c)=c\otimes I\}={\mathbb C}$. Recall that ${\cal C}$ has a unique faithful state invariant under the action of $G$,  but, unlike in the group case, is not a trace in general \cite{Woronowicz}, see also
\cite{Wang}.

\medskip

\noindent{\it 2.5. Spectrum, multiplicity maps, spectral functor and duality theorem.}
The {\it spectrum} of an action $\alpha$ of a compact 
quantum group $G$ on ${\cal C}$, denoted $\text{sp}(\alpha)$,
 is   the set of all unitary
representations $u$ of $G$ for which there is a faithful linear map
$T: H_u\to{\cal C}$ intertwining the representation $u$ with the action 
$\alpha$. This means that
  if $u_{ij}$ are the 
coefficients of $u$ in some orthonormal basis of $H$,
there are linearly independent elements
$c_1,\dots, c_d\in{\cal C}$, with $d$ the dimension of $u$, such that
$\alpha(c_i):=\sum_j c_j\otimes u_{ji}$. The linear span of all the $c_i$'s, 
 as $u$ varies in the spectrum, is a dense invariant $^*$--subalgebra of ${\cal C}$, denoted ${\cal C}_{\text{sp}}$
\cite{Podles}.

Examples of ergodic actions are the quantum quotient spaces $C(K\backslash G)$ by a compact quantum subgroup $K$.  As in the classical case, $C(K\backslash G)$ is the fixed point algebra under a suitable action of $K$ on   the Hopf $C^*$--algebra of $G$,
with action of $G$   given by restricting the coproduct. This action is usually called the {\it translation action}  \cite{Wang}, see also \cite{Pinzari3}.
 $C(K\backslash G)_{\text{sp}}$ 
is linearly spanned by the matrix coefficients $\{u_{k,\psi_i}\}$, where  $u$ varies in the set of unitary f.d. representations of $G$,  $k$ in the set of fixed vectors  for the restriction $u\upharpoonright_K$ and $(\psi_i)$ is an orthonormal basis. 

For any representation $u$, consider the space $L_u$ of all linear intertwiners $T$, not necessarily faithful,  between $u$ and $\alpha$.
$L_u$ becomes a Hilbert space with inner product
$<S,T>:=\sum_iT(\psi_i)S(\psi_i)^*$, with $(\psi_i)$ an orthonormal basis of $H_u$.
For an 
irreducible $u$, $L_u\neq0$ precisely
when $v\in\text{sp}(\alpha)$.
 The dimension of $L_u$ is the {\it multiplicity}
of $u$ in ${\alpha}$.
$L_u$ is known to be finite dimensional if $\alpha$ is ergodic \cite{Boca}. 
The complex conjugate vector space $\overline{L_u}$,
 is called the {\it spectral space} associated with $u$.
 For any $u\in\text{Rep}(G)$, associate  the map
$$c_u: H_u\to\overline{L_u}\otimes{\cal C}\quad c_u(\psi):=\sum_k\overline{T_k}\otimes T_k(\psi),$$
where $T_k$ is any orthonormal basis of $L_u$. Note that $c_u$ does not depend on the choice of  orthonormal basis.
The  $c_u$'s are called {\it multiplicity 
maps} in \cite{PR}.

We can represent $c_u$ as a rectangular matrix whose $i$-th
row is given by the multiplet $T_i=(T_i(\psi_1)\dots T_i(\psi_d))$
transforming like $u$ under $\alpha$.

The set of all coefficients 
 $\{c^u_{i,j}=T_i(\psi_j), i,j\}$
forms a linear basis for the dense $^*$--subalgebra 
${\cal C}_{\text{sp}}$  
as $u$ varies over a complete set of irreducible representations of 
$\text{sp}(\alpha)$, \cite{Podles, BDV, PR}.

 The map $u\mapsto{L_u}$ 
can be extended to a {\it functor} 
${L}:\text{Rep}(G)\to\text{Hilb}$
from the category $\text{Rep}(G)$ of  unitary f.d. representations of $G$ 
to the category   of Hilbert 
spaces. 
If $A\in(u,v)$ and $T\in L_v$ then $T\circ A: H_u\to{\cal C}$
lies in $L_u$.  $L$ is a contravariant $^*$--functor, hence it is 
convenient to pass to the dual Hilbert spaces, that we naturally identify with the spectral spaces
$\overline{L}_u$. We thus get    a covariant 
$^*$--functor,
$\overline{L}$, 
the {\it spectral functor} of the ergodic action. The spectral functor  and the multiplicity maps   
are related by
$\overline{L}_A\otimes I\circ c_u=c_v\circ A,\quad A\in(u,v),$
for any $u,v\in\text{Rep}(G)$.

 For example, the spectral functor of a quantum quotient space $C(K\backslash G)$ maps the representation $u$ of $G$ to the Hilbert space of fixed vectors of the restriction $u\upharpoonright_K$.

There is a 
natural isometric inclusion
$$\tilde{\overline L}_{u,v}:\overline{L}_u\otimes\overline{L}_v\to\overline{L}_{u\otimes 
v}$$
identifying a simple tensor $\overline{S}\otimes\overline{T}$ with
the complex conjugate of the element of $L_{u\otimes v}$ defined by
$\psi\otimes\phi\in H_u\otimes H_v\to S(\psi)T(\phi).$
It has been shown in \cite{PR} that $(\overline{L},\tilde{\overline{L}})$ is a quasitensor functor and that $({\cal C},\alpha)$ may be reconstructed almost entirely from $\overline{L}$. In detail,  $\overline{L}$ keeps complete information on the dense $^*$--subalgebra 
${\cal C}_{\text{sp}}$, its maximal $C^*$--norm and the restricted
action $\alpha\upharpoonright_{{\cal C}_{\text{sp}}}$.
(However, $\overline{L}$ does not keep track of the original $C^*$--norm of ${\cal C}$, a feature already present in Woronowicz' version of Tannaka--Krein duality \cite{WoronowiczTK}.)
Moreover, any quasitensor functor $(\mu,\tilde{\mu}):\text{Rep}(G)\to\text{Hilb}$ is the spectral functor of  an ergodic action of $G$ on a unital $C^*$-algebra.

 For completeness, we recall that the spectral functor  is a relaxed tensor functor if and only if the  quantum multiplicity
of every irreducible equals its quantum dimension. We refer 
to \cite{BDV} for  the notion of quantum multiplicity and to \cite{PR} for the proof of this fact. In the group case, this means that the multiplicity of every irreducible equals its dimension.

Note that the finiteness theorem fails for 
compact quantum groups, as the Haar measure is not a trace in general. On the other hand, the multiplicity bound theorem holds, provided multiplicity and dimension are replaced by their noncommutative analogues \cite{Boca}, \cite{BDV}.
\medskip

\noindent{\it 2.6. Hilbert modules and Hilbert bimodules over $C^*$--algebras.}
We refer to 
\cite{Blackadar, Lance} for the definition of a (right) Hilbert module $X$ over a $C^*$--algebra
${\cal C}$. The ${\cal C}$--valued inner product will be denoted by $<\cdot{}\ ,\ \cdot{} >$ and we shall assume ${\cal C}$--linearity on the right.
We recall in particular that if ${\cal C}={\mathbb C}$, $X$ is just a Hilbert space. 
Any $C^*$--algebra ${\cal C}$ gives rise to the Hilbert module, $X={\cal C}$, with inner product $<c,c'>:=c^*c'$. This is usually called the {\it trivial} Hilbert module. 
More generally, if $H$ is a Hilbert space, we may form the
algebraic tensor product $H\odot {\cal C}$, which is a pre-Hilbert module with the obvious module structure and inner product
$$<\psi\otimes c,\psi'\otimes c'>=<\psi,\psi'>c^*c'.$$ 
The completion will be denoted by
$H\otimes {\cal C}$. We shall only consider Hilbert modules over unital $C^*$--algebras.
A Hilbert ${\cal C}$--module $X$ is called {\it full} if the inner products $<x,x'>$ span a dense subspace of ${\cal C}$.

If $X$ and $X'$ are Hilbert modules over ${\cal C}$, ${\cal L}_{\cal C}(X, X')$ denotes
the Banach space of bounded adjointable maps from $X$ to $X'$.
 
A Hilbert module $X$ over ${\cal C}$ will be called a  {\it  Hilbert bimodule} if there is
a left  action of ${\cal C}$ on $X$ given by a unital $^*$--homomorphism
${\cal C}\to{\cal L}_{\cal C}(X,X)$. For example, the trivial module becomes a Hilbert bimodule in the obvious way. If $X'$ is another Hilbert bimodule
over ${\cal C}$, ${}_{\cal C}{\cal L}_{\cal C}(X, X')$ denotes the set of elements 
$T\in{\cal L}_{\cal C}(X, X')$ commuting with the left actions. 

The great advantage of Hilbert bimodules versus Hilbert modules is that we may form tensor products, $X\otimes_{\cal C} Y$ see \cite{Blackadar, Connes}.  The category with objects Hilbert bimodules over ${\cal C}$ and arrows $(X, X'):={}_{\cal C}{\cal L}_{\cal C}(X, X')$ is a tensor $C^*$--category with tensor unit given by the trivial Hilbert bimodule.

\section{Bimodule representations of compact quantum groups}

In this section we define unitary  representations of compact quantum groups 
on Hilbert modules or Hilbert bimodules over  unital $C^*$--algebras.  These representations may be regarded as the noncommutative analogues of 
 the  $G$--equivariant Hermitian bundles over compact spaces introduced by Segal 
\cite{Segal}, where 
$G$ is a compact group.

 In the following general definition we shall not assume  our modules to be 
finite projective 
 (this would correspond to local triviality in the commutative case, by 
Swan's theorem \cite{Swan}), even though 
we shall  eventually be 
interested in finite projective Hilbert modules.

Let us fix an action $({\cal C},\alpha)$ of a compact quantum group $G=({\cal Q}, \Delta)$ 
on   a unital $C^*$--algebra ${\cal C}$.   By a  {\it Hilbert
module 
representation} of $G$, or simply a {\it module representation}, 
we mean   a ${\Bbb C}$--linear map
$$v: X_v\to X_v\otimes{\cal Q},$$
where $X_v$ is a   Hilbert ${\cal C}$--module, ${\cal Q}$ is regarded as the trivial Hilbert ${\cal Q}$--module and 
$X_v\otimes {\cal Q}$ denotes the exterior tensor product of Hilbert modules, a Hilbert module over ${\cal C}\otimes{\cal Q}$, see \cite{Blackadar} for details,
such that
$$<v(x), v(x')>_{{\cal C}\otimes{\cal Q}}=\alpha(<x,x'>_{\cal C}),\quad
x,x'\in
X_v,\eqno(3.1)$$
$$v(xc)=v(x)\alpha(c),\quad x\in X_v, c\in {\cal C},\eqno(3.2)$$
$$v\otimes 1_{{\cal Q}}\circ v=1_{X_v}\otimes\Delta\circ v,\eqno(3.3)$$
$$v(X_v)1_{X_v}\otimes{\cal Q}\text{ is dense in } X_v\otimes{\cal Q}.\eqno(3.4)$$
The simplest example of a module representation is the {\it trivial representation},
$v=\alpha$ on the trivial Hilbert ${\cal C}$--module. It will be denoted by $\iota$.

Note that if  ${\cal C}={\mathbb C}$ this definition reduces to the notion of a  strongly continuous unitary
representation  of a compact quantum group  on a Hilbert space, see \cite{BS, Wlh}. 

On the other hand, as mentioned at the beginning of the section, if ${\cal C}$ is commutative and $G$ is a compact group, this notion
reduces to that of a $G$--equivariant Hermitian bundle, the equivariance property being expressed by $(3.2)$.

One can form
  the $C^*$--category $\text{Mod}_\alpha(G)$ with
objects the module representations of $G$
and arrows 
$$(v,v'):=\{T\in{\cal L}_{\cal C}(X_v, X_{v'}): v'\circ T=
T\otimes
1_{\cal Q}\circ v\}.$$
Note that  $(\iota,\iota)$ can be identified with the fixed
point algebra ${\cal C}^\alpha$.

\medskip

We are interested in module representations of a compact quantum group where $X_v$ is a 
Hilbert  bimodule.
 $v$ will be called a 
 {\it (Hilbert) bimodule representation}  of $G$ if in addition to $(3.1)$--$(3.4)$,
$$v(cx)=\alpha(c)v(x),\quad c\in{\cal C}, x\in X_v.\eqno(3.5)$$

As an example,
the trivial representation $\iota$ is a bimodule representation, that we shall denote by the same symbol.
\medskip

We denote by  $\text{Bimod}_\alpha(G)$ the category with objects Hilbert bimodule representations and  arrows  $(v,v')$ between two of them  the space of  intertwining operators in $\text{Mod}_\alpha(G)$ which in addition intertwine  the left actions of ${\cal C}$.
If $u,v$ are two objects of 
$\text{Bimod}_\alpha(G)$ 
we define,
for $x\in X_u$, $y\in X_v$,
the {\it tensor product bimodule representation}
$u\otimes v$  by
$$u\otimes v (x\otimes y)=u(x)_{13}v(y)_{23},$$
an element of $X_u\otimes_{\cal C} X_v\otimes{\cal Q}$.
$(3.2)$ and $(3.5)$ show that
$u\otimes v$ is well defined on the algebraic bimodule tensor product $X_u\odot_{\cal C} X_v$ and that $(3.1)$, $(3.2)$, $(3.3)$ and $(3.5)$ 
hold. The validity of $(3.1)$ implies that $u\otimes v$
extends uniquely to a bounded ${\mathbb C}$--linear map
$$u\otimes v: X_u\otimes_{\cal C} X_v\to X_u\otimes_{\cal C}
X_v\otimes{\cal Q},$$ and the above equations still hold whereas $(3.4)$ holds by construction. The tensor product of two 
intertwiners is now well defined and intertwines the tensor product 
representations. Note that if $v$ is an object of $\text{Bimod}_\alpha(G)$, $v\otimes\iota$ and $\iota\otimes v$ are equivalent to $v$ 
in $\text{Bimod}_\alpha(G)$.
This leads to the following result.
\medskip

\noindent{\bf 3.1. Theorem} {\sl The category
$\text{Bimod}_\alpha(G)$ with objects Hilbert bimodule 
representations  of $G$  and arrows the
bimodule intertwining operators is
a tensor $C^*$--category. The tensor unit  
is the trivial representation $\iota$ and 
$(\iota,\iota)=\{c\in Z({\cal C}): \alpha(c)=c\otimes I\},$ the set of 
central  fixed points.
There is an obvious  faithful $^*$--functor
$\text{Bimod}_\alpha(G)\to\text{Mod}_\alpha(G).$}\medskip

We shall only consider 
ergodic  actions, hence  
 $(\iota,\iota)={\mathbb C}$. \medskip

\section {The induced $C^*$--bimodules for compact groups}

We  recall  Mackey's definition \cite{Mackey} of a
 representation induced from a closed subgroup of a compact 
group and the 
Frobenius reciprocity theorem in the form later generalized to tensor $C^*$--categories.
 
Our discussion has 
points in common with
 \cite{Rieffel}. The main point is that we shall pass from Hilbert space representations
 to Hilbert bimodule representations.
The module approach to induction is particularly convenient 
in the compact case as it provides finite dimensional objects
by Swan's theorem \cite{Swan} and moreover the induction functor has good tensorial properties
(cf. Theorem 4.1).
\medskip

\noindent {\it 4.1. Mackey's  induced representation and Frobenius reciprocity.}
Let   $K$ be  a closed subgroup of a compact group $G$ and 
$v$ a (unitary,  finite dimensional) representation of $K$ on the Hilbert space $H_v$. Mackey's induced 
representation $\text{Ind}(v)$ is defined as right translation by 
elements of $G$ on the Hilbert space  of $L^2$ 
functions $\zeta$ on $G$ with values in $H_v$ satisfying 
$$\zeta(kg)=v(k)\zeta(g),\quad k\in K, g\in G,$$
where the inner product 
$<\zeta,\zeta'>=\int_{K\backslash 
G}<\zeta(g),\zeta'(g)>d\mu$
involves the unique normalized
$G$--invariant measure
$\mu$ on $K\backslash G$. 

The main result is the Frobenius reciprocity 
 theorem, asserting that there is an explicit linear isomorphism from  the intertwining space
$(u\upharpoonright_K,v)$
to
$(u,\text{Ind}(v))$, taking an 
intertwiner $S$ to the intertwiner $T$, where $T(\psi)(g)=S(u(g)\psi)$. 
The Frobenius  isomorphism  is natural  in $u$ and $v$, and hence makes restriction and induction into a pair of adjoint functors. (In Sect.\ 10, we will briefly recall the notion of an adjoint pair of functors. For details we refer the reader
to  MacLane's book \cite{MacLane}.)
Consequently, the spectrum of the induced representation $\text{Ind}(v)$ is the set of all irreducible 
$G$--representations $u$ for which $(u\upharpoonright_K, 
v)$ is nonzero. Another consequence  is that any 
irreducible representation $v$ of $K$, and hence any $v$, is a 
subrepresentation of some 
restriction to $K$ of a representation $u$ of $G$. 
 Finally, the explicit form of the isomorphism shows that 
if $T\in (u, \text{Ind}(v))$ all the functions $T(\psi)$ are   {\it continuous}.   This last remark leads to the next step.
\medskip

\noindent{\it 4.2. Replacing Hilbert spaces with Hilbert bimodules.} Since  we do not loose any information on the arrows, we may replace Hilbert spaces with Hilbert bimodules.  More precisely, we
 pass    from  the  Hilbert space of the induced representation to
the space $C_v$
of continuous $H_v$--valued  functions $\zeta$ as above, which is  a 
 bimodule
over the commutative $C^*$--algebra 
$C(K\backslash 
G)$ of continuous functions on the quotient space in the obvious way.
$C_v$ has an inner 
 product given by pointwise evaluation of the inner product of $H_v$,
 $$<\zeta,\zeta'>(g):=<\zeta(g),\zeta'(g)>.$$
This inner product is constant on each left coset $Kg$ as $v$ is unitary, 
and $C_v$ becomes a {\it Hilbert bimodule}  over 
$C(K\backslash G)$. 
Hence $\text{Ind}(v)$ becomes  a Hilbert $C(K\backslash 
G)$--bimodule representation of $G$
in the sense of the previous section, where $\alpha$ is
given by right translation by elements of $G$ on the quotient space.
Note that $C_v$ is the bimodule of continuous sections of   the classical equivariant 
vector bundle induced from 
$v$.
\medskip

\noindent {\it 4.3. The induction functor and Swan's theorem.}
We thus have a  $^*$--functor  
$$\text{Ind}: \text{Rep}(K)\to\text{Bimod}_\alpha(G),$$  
taking an object $v$ of $\text{Rep}(K)$ to $\text{Ind}(v)$ and an arrow $T\in(v,v')$ 
 to the  arrow $\text{Ind}(T)\in(\text{Ind}(v), 
\text{Ind}(v'))$ defined by $\text{Ind}(T)\xi(g)=T\xi(g)$.
We shall refer to $\text{Ind}$ as the {\it induction functor}.

If $u\upharpoonright_K$ is the restriction of a representation $u$ of $G$ to $K$, there is  a natural faithful bimodule map,
$$U: C_{u\upharpoonright_K}\to H_u\otimes C(K\backslash G),$$
$$U\zeta(g)=u(g^{-1})\zeta(g).$$
$U$ is  invertible, and hence surjective, with inverse
given by $U^{-1}\xi(g)=u(g)\xi(g)$.  Hence $C_{u\upharpoonright_K}$ is free as a 
right  (and left) module. Moreover,
 $U$ becomes unitary when $H_u\otimes C(K\backslash G)$ is regarded as 
a Hilbert bimodule. 
Note that the tensor product action $u\otimes\alpha$ of $G$ on
$H_u\otimes C(K\backslash G)$
is a 
Hilbert bimodule representation 
of $G$ and $U$ becomes a unitary intertwiner from
$\text{Ind}(u\upharpoonright_K)$ to $u\otimes\alpha$ in 
$\text{Bimod}_\alpha(G)$.
 \medskip

 Let now $v$ be a generic unitary  finite dimensional representation of $K$, and consider a restricted representation $u\upharpoonright_K$ containing   $v$ as a subrepresentation. An isometric intertwiner in $(v, u\upharpoonright_K)$ defines an isometric 
 intertwiner in $(\text{Ind}(v),\text{Ind}(u\upharpoonright_K))$ between the associated  Hilbert bimodules representations,
 via the induction functor. Therefore   $\text{Ind}(v)$ 
 is a subobject of $\text{Ind}(u\upharpoonright_K)$. 
 Moreover, composition with 
the map $U$ as above, gives rise to an isometry of Hilbert modules from $C_v$ to  
$H_u\otimes C(K\backslash G)$. Hence
$C_v$ is a {\it finite projective} module.
This is essentially Rieffel's proof \cite{Rieffel} of  Swan's 
theorem \cite{Swan}.

 Taking into account the previous remark, we may conclude that bimodule representations of the form 
$u\otimes\alpha$ suffice to generate the category of all induced bimodule representations via subobjects. 
This viewpoint will play a  role   in the next sections.

On the other hand, the naturality of the Frobenius isomorphism shows that any element of $(\text{Ind}(v), \text{Ind}(v'))$ is of the form $\text{Ind}(T)$ for a unique $T\in(v,v')$, and so   $\text{Ind}$ is a {\it full} functor.

\medskip

\noindent{\it 4.4. Tensorial properties of the induction functor.}
We next analyse the behaviour of $\text{Ind}$ under tensor products. 
We may consider the  
tensor product of Hilbert bimodules,
$C_u\otimes C_v:=C_u\otimes_{C(K\backslash G)}C_{v}$.
 There are obvious  isometric inclusions of $G$--bimodule representations
$C_u\otimes  C_v\to C_{u\otimes v}$.
These maps are in fact surjective, and hence unitary, as, given 
isometries
$S\in(u, u'\upharpoonright_K)$ and 
$T\in(v,v'\upharpoonright_K)$, module bases for $C_u$ 
and $C_v$ are given by the functions 
$(x_i:=g\to S^*u'(g)\psi_i)$ 
and
$(y_j:=g\to T^*v'(g)\phi_j)$  respectively, where $\psi_i$ and $\phi_j$ are 
orthonormal bases of the Hilbert spaces of $u'$ and $v'$. 
Hence, as  $u\otimes v$ is a 
subrepresentation 
of $(u'\otimes v')\upharpoonright_K$,  a module basis for 
$C_{u\otimes v}$ is given by the functions $g\to S^*\otimes T^*u'\otimes 
v'(g)\psi_i\otimes\phi_j=x_i(g)\otimes y_j(g)$. 
Thus $C_{u\otimes v}$ can be  naturally identified  with $C_u\otimes C_v$.
 It follows that $\text{Ind}$ is a  relaxed tensor functor  in the sense recalled in Sect. 2. 
We summarize this discussion in the following 
theorem, essentially  a geometric form of the Frobenius 
reciprocity theorem.
\medskip

\noindent{\bf 4.1. Theorem} {\sl The induction functor $\text{Ind}:\text{Rep}(K)\to\text{Bimod}_\alpha(G)$ is a full and faithful $^*$--functor
 into the bimodule 
representation category of $G$, where  $\alpha$ is 
given by right translation of $G$ on $C(K\backslash G)$. For any 
$v\in\text{Rep}(K)$, the $C(K\backslash G)$--bimodule $C_v$ of 
$\text{Ind}(v)$ is finite 
projective. In particular, if $v$ is the restriction of a 
representation of $G$, $C_v$ is free.
 The natural unitaries $C_u\otimes C_{v}\to 
C_{u\otimes v}$ make $\text{Ind}$ into a relaxed tensor functor.}\medskip

\section{Full bimodule representations }

In this section we introduce a  notion central to this paper, that of a
{\it  full bimodule representation } of a compact quantum group $G$. This is a compatibility condition between the left and right bimodule structure of a Hilbert $C^*$--bimodule carrying a representation of $G$.  As we shall see, in the classical case, triviality of the bimodule structure means that every induced bimodule representation is full. The importance of this 
notion is that all induced  bimodule representations of compact quantum groups constructed in this paper are full. 

In the noncommutative situation, $C(K\backslash G)$ with the translation action is replaced by a compact quantum group $G$ acting ergodically on a unital $C^*$--algebra ${\cal C}$. The action will be denoted by $\alpha$.
 It would be too restrictive to consider just quantum quotient spaces $K\backslash G$, as $G$ can act on far more noncommutative $C^*$--algebras. Following Mackey, we may regard 
the ergodic action $({\cal C},\alpha)$ as arising from a virtual subgroup.

Restricting or inducing a representation  now looses its strict meaning. 
 What is left   is the analogue of 
$\text{Ind}(u\upharpoonright_K)$, which may be defined as acting on the free Hilbert module $H_u\otimes{\cal C}$.
 More 
precisely, it is easy to see that the map 
$u\otimes\alpha$ defined by
$$u\otimes\alpha(\psi\otimes c):=
u(\psi)_{13}\alpha(c)_{23}\in H_u\otimes{\cal C}\otimes{\cal Q},\quad \psi\in H_u, c\in {\cal C},$$
is  a Hilbert module representation of $G$ on $H_u\otimes{\cal C}$.  

Given $T\in(u,v)$, with $u$, $v\in\text{Rep}(G)$, define   $ T\otimes 1_{\cal C}: H_u\otimes{\cal C}\to H_v\otimes{\cal C}$ by ${T\otimes 1_{\cal C}}(\psi\otimes c)=T\psi\otimes c$ then 
$ T\otimes 1_{\cal C}\in(u\otimes\alpha, 
v\otimes\alpha)$ in $\text{Mod}_\alpha(G)$.
\medskip

\noindent{\bf 5.1. Proposition} {\sl The map $\text{Rep}(G)\to \text{Mod}_\alpha(G)$,
taking $u\to u\otimes\alpha$ and $T\in(u, v)\to T\otimes 1_{\cal C}$, 
 is a faithful $^*$--functor between $C^*$--categories.}\medskip

The previous proposition is a very weak form of Theorem 4.1 and
our aim is to generalize it  to the noncommutative setting. This involves 
giving 
 $H_u\otimes{\cal C}$ a left module structure 
making
$u\otimes\alpha$ into a Hilbert  bimodule 
$G$--representation for all $u$. 
Although  classical induction corresponds to the simplest solution 
(trivial bimodule structure), 
in noncommutative case, it is not obvious how to select ergodic actions 
$({\cal C},\alpha)$ in such a way that  the module representation 
$u\otimes\alpha$ can be completed to a bimodule representation. 
Even if this is the case, many  left module structures will not be relevant, and we claim that a compatibility condition between left and right module structure is needed.
This may already be seen for a compact group $G$ acting on 
 a noncommutative $C^*$--algebra ${\cal C}$.

In this case the most natural left module structure on 
$H_u\otimes{\cal C}$  is the obvious left multiplication by elements of ${\cal C}$. With this left ${\cal C}$--action, the right Hilbert module representation $u\otimes\alpha$ becomes  a Hilbert bimodule representation.

However, an intertwiner between two such bimodule representations, 
being  a bimodule map, must  lie in 
${\cal B}(H_u, H_{u'})\otimes Z({\cal C}),$ where $Z({\cal C})$ is the centre of ${\cal C}$. Hence these intertwiners  
do not see the noncommutativity of ${\cal C}$, in contrast
to the module intertwining spaces, where  
$(\iota, u\otimes\alpha)$ is the space of 
multiplets 
$\xi=(c_1,\dots, 
c_d)$, with $c_i\in{\cal C}$ and $\alpha(c_i)=\sum_j c_j\otimes u_{ji}^*$, 
$d$ being the dimension of $u$. As $u$ varies 
over the irreducible spectral representations of $G$,  the  linear span of 
the corresponding  $c_i$'s is a dense invariant $^*$--subalgebra, cf. Sect. 2. 

Hence the natural left action 
on $H_u\otimes{\cal C}$  
gives rise to a tensor category which
does {\it not} allow one to reconstruct ${\cal C}$, but only its centre.    
It would be desirable to use instead a left ${\cal C}$--action on  
$H_u\otimes{\cal C}$  where all
the module $G$--intertwiners become bimodule $G$--intertwiners. 
This leads us to the    
notion of {\it full} bimodule representation.
\medskip

\noindent{\it Definition}
Let $G$ be a compact quantum group.
 A {\it fixed vector} $\xi$ for a module representation $v$ on $X_v$ is  an element
$\xi\in X_v$ such that $v(\xi)=\xi\otimes I$.
The set of fixed vectors for $v$ is precisely the intertwining space 
$(\iota, v)$ in $\text{Mod}_\alpha(G)$.
A bimodule 
representation $v$ will be called {\it full} if every fixed vector $\xi$ 
for the underlying module $G$--representation 
is central: $\xi c=c\xi$, for $c\in{\cal C}$.
\medskip

Note that the trivial representation is full since $\alpha$ is 
ergodic.

The next result shows that if $G$ is a group, classical induction
is characterized among functors $u\to u\otimes\alpha$ from $\text{Rep}(G)$ to $\text{Mod}_\alpha(G)$ by the property 
that under
the natural left action each 
$u\otimes\alpha$ becomes a
full bimodule representation. 
\medskip

\noindent{\bf 5.2. Proposition}   {\sl Let $G$  be a compact group and $\alpha$ an ergodic 
 action of 
$G$ on a unital $C^*$--algebra ${\cal C}$. Then the natural left ${\cal C}$--action
turns $u\otimes\alpha$ into
 a full representation for all 
$u\in\text{Rep}(G)$ if and only if ${\cal C}$ is commutative.
 In this case,  ${\cal C}=C(K\backslash G)$ for a closed 
subgroup
$K$, unique up to conjugation, where $\alpha$  acts by right translation,
$\alpha_gf(g')=f(g'g)$, $f\in C(K\backslash G)$. Hence $u\otimes\alpha$ 
corresponds to the classical induced representation $\text{Ind}(u\upharpoonright_K)$.
}\medskip

\noindent{\it Proof} If ${\cal C}$ is commutative, $c\eta=\eta c$ for 
$c\in{\cal C}$, $\eta\in H_u\otimes{\cal C}$ and $u\in\text{Rep}(G)$.
Hence 
any module intertwiner 
between $u\otimes\alpha$ and $u'\otimes\alpha$ is a bimodule intertwiner. In particular, 
$u\otimes\alpha$ is full for all $u$. Conversely, 
assume that 
all the $u\otimes\alpha$ are full. We 
have 
already seen that a fixed vector 
$\xi$ for $u\otimes\alpha$ 
has the form
$\xi=\sum_j\psi_j\otimes c_i$ for an orthonormal basis $(\psi_j)$ of 
$H_u$ where $c_j$ transforms like the complex conjugate representation $u_*=(u_{ij}^*)$   under 
$\alpha$. Since $\xi$ is supposed central, the elements $c_i$ are central 
in ${\cal C}$. If $u$ varies in the spectrum of $\alpha$, we 
get a dense commutative $^*$--subalgebra of ${\cal C}$, hence ${\cal C}$ 
is commutative.
As is well known, when $\alpha$ is an ergodic action on a unital commutative $C^*$--algebra,  the action is   right 
translation on
 $C(K\backslash G)$  by elements of $G$ for a closed subgroup $K$ of $G$, unique
up to conjugation.  Hence $u\otimes\alpha$ can be identified with $\text{Ind}(u\upharpoonright_K)$.
\medskip

Of course, we  expect that, requiring all $u\otimes\alpha$ to be full bimodule representations for some left module structure, selects a proper subclass of ergodic actions. This is the case, as we shall see (cf. Sect. 6 and references there to later sections).

\medskip

The situation becomes significantly worse if $G$ is a quantum group,
where the natural left action on $H_u\otimes{\cal 
C}$ may not even lead to a bimodule representation structure on 
$u\otimes\alpha$.
\medskip

 \noindent{\bf 5.3. Proposition} {\sl Let $G$ be a compact quantum group, 
$u\in\text{Rep}(G)$ and $\alpha$ an ergodic action of $G$ 
on a unital $C^*$--algebra ${\cal C}$. 
Then $u\otimes\alpha$ is a bimodule 
representation  
for the natural left module 
structure if and only if all coefficients of  the 
irreducibles in the spectrum of $\alpha$ commute with the 
coefficients of 
$u$.}\medskip

\noindent{\it Proof}   ${\cal C}$ is 
generated as a Banach space, by the entries of rectangular matrices 
$(c^v_{i,j})$
 transforming like irreducible $G$--representations 
$$\alpha(c^v_{i,j})=\sum_p c^v_{i, p}\otimes v_{p,j}.$$
These entries are linearly independent \cite{Podles}, 
\cite{BDV}, \cite{PR}
 and the conclusion follows from $(3.5)$.\medskip
 
 We next discuss  examples of full bimodule representations arising from quantum quotients.\medskip
 
 \noindent {\it   Examples from quantum quotients.}
Let $G=({\cal Q},\Delta)$ be a compact quantum group, and $K$ a quantum subgroup.
$G$ acts on the  quotient space 
$K\backslash G$ by right translation, given by restricting  
the coproduct $\Delta$ of $G$.
One can consider the 
left $C(K\backslash G)$--action on $H_u\otimes C(K\backslash G)$
 defined as follows. For $c\in C(K\backslash G)$,
consider 
the element $\lambda_u(c)\in{\cal L}(H_u)\otimes{\cal Q}$ defined
by
$$<\psi_i\otimes I,\lambda_u(c)\psi_j\otimes I>:=\sum_h u_{hi}^* 
c u_{hj},$$
where $(\psi_i)$ is  an orthonormal basis of $H_u$.
It is easy to check that this element is independent of the 
choice of orthonormal basis. One could  show directly 
that $\lambda_u(c)\in{\cal L}(H_u)\otimes C(K\backslash G)$ and that $\lambda_u$   
makes $u\otimes\alpha$ into a full  bimodule 
representation.
However, we refrain from giving complete details, as this will be 
proved
in more generality  
in Sect.\ 8. We just verify that it is full.
If $\xi:=\sum \psi_j\otimes c_j\in H_u\otimes{\cal C}$ is a fixed vector, 
i.e., $\Delta(c_i)=\sum_j c_j\otimes u^*_{ji}$,
then 
$$c\xi=\sum_k \psi_k\otimes<\psi_k\otimes I, 
\lambda_u(c)\xi>=\sum_{k,j,h}\psi_k\otimes u^*_{hk}cu_{hj}c_j$$
whereas
$$\xi c=\sum_k\psi_k\otimes c_k c,$$
hence we need to show that for every $c\in C(K\backslash G)$,
$$\sum_{h,j} u^*_{hk} c u_{hj}c_j=c_k c$$
for all $k$. On the other hand, for a quotient 
space, we 
can find a $K$--fixed vector $\eta$ of $H_u$ such that
$c_j=u^*_{\eta,j}$, and the desired equality follows from the unitarity of 
$u$, see, e.g., Sect. 2 in  \cite{PR}.
\medskip

\medskip

We are thus left  with the 
problem of finding full Hilbert bimodules that lead to an  induction functor with good tensorial properties. We shall show that tensor $C^*$--categories with conjugates provide  a natural solution to this problem, as well as 
large classes of ergodic actions of compact quantum groups, among them the compact quantum 
quotient spaces. However, many more will be discussed in Sect.\ 6.

\section{Main results}

In this section we illustrate  the main ideas and
results. Proofs will be given in later sections.

We start
with a pair of tensor
$C^*$--categories ${\cal A}$ and ${\cal M}$   related by a quasitensor functor $\mu:{\cal A}\to{\cal M}$.
The category ${\cal A}$  is assumed to be
embeddable into Hilbert spaces and we then pick a tensor functor $\tau:{\cal
A}\to\text{Hilb}$. We will   assume that ${\cal A}$ has conjugates.

The simplest  example   is provided by a closed subgroup $K$ of a compact 
group $G$.  We may choose  $\tau:\text{Rep}(G)\to\text{Hilb}$ the embedding functor and $\mu:\text{Rep}(G)\to\text{Rep}(K)$ the tensor  functor restricting a representation of $G$ to
$K$. 

Note however that this example has certain special features, like the 
fact that ${\cal M}=\text{Rep}(K)$ is embeddable, or that $\mu$ is 
tensorial.
In general, ${\cal M}$ is not assumed to be embeddable, and, as recalled 
in Sect.\ 
2, a quasitensor functor $\mu$, unlike a relaxed tensor functor, may take 
a nonzero object to a zero object.

By Woronowicz duality, the embedding $\tau$ defines  a compact quantum
group $G_\tau$  such that every object $u\in{\cal A}$ has an associated representation $\hat{u}\in\text{Rep}(G_\tau)$ on the Hilbert space $\tau_u$.  
The arrow spaces of $\text{Rep}(G_\tau)$ are the images under $\tau$ of the arrow spaces of ${\cal A}$,   $(\hat{u},\hat{v})=\tau((u,v))$.

The pair $(\tau,\mu)$   determines canonically a
unital $C^*$--algebra, ${\cal C}$ and an ergodic action $\alpha$ of
$G_\tau$ on ${\cal C}$.  This fact has been shown in \cite{PR} in the  
 special case where  $\mu:{\cal A}\to{\cal M}$ is a tensor functor. To see that this holds in our more  general setting 
consider the composed $^*$--functor
$$F:{\cal A}\stackrel\mu\rightarrow{\cal M}\to\text{Hilb}$$
where the second functor is the so called minimal functor $x\in{\cal 
M}\to(\iota, x)\in{\text{Hilb}}$, which is quasitensor, see \cite{PR}.
Since composition of quasitensor functors is quasitensor (cf. subsect.\ 2.2),  
so is $F$. Thus there is a unique quasitensor functor 
$(\nu,\tilde{\nu}):\text{Rep}(G_\tau)\to\text{Hilb}$ such that the 
following diagram commutes
\[
\xymatrix{
{\cal A} \ar[r]^F\ar[rd]_{\tau}&
\text{Hilb}\\
{}&\text{Rep}(G_\tau)\ar[u]_\nu
}
\]
 It follows   that $(\nu,\tilde{\nu})$ is the spectral functor of an 
ergodic $C^*$--action of $G_\tau$ (cf. subsect.\ 2.5), and this is 
 $({\cal C},\alpha)$. The following simple remark clarifies matters.
 \medskip
 
 \noindent{\bf 6.1. Proposition} {\sl  For a fixed $\tau:{\cal A}\to\text{Hilb}$, let $(\mu,\tilde{\mu}):{\cal A}\to{\cal M}$,
 $(\mu', \tilde{\mu'}):{\cal A}\to{\cal M}'$ be a pair of quasitensor functors. 
 If there is a full relaxed tensor functor $(\phi,\tilde{\phi}):{\cal M}\to{\cal M}'$
 such that $\phi\mu=\mu'$ 
 then the associated ergodic $C^*$--actions $({\cal C},\alpha)$
 and $({\cal C}',\alpha')$ are conjugate.}\medskip

 \noindent{\it Proof} By Prop. 8.4 in \cite{PR} two ergodic actions are conjugate if their  spectral functors
 are related by a quasitensor natural unitary transformation. Explicitly, if $(\overline{L},\tilde{\overline{L}})$ and $(\overline{L}',\tilde{\overline{L}'})$ are the spectral functors of  the actions, we  need a unitary
  $U_u:\overline{L}_u\to\overline{L}'_u$ for each $u\in\text{Rep}(G_\tau)$, 
natural in $u$, such that $U_{u\otimes v}\circ\tilde{\overline{L}}_{u,v}=
\tilde{\overline{L}'}_{u,v}\circ U_u\otimes U_v$ and $U_\iota=1_\iota$. 
Now the spectral space of the ergodic action
  constructed from $(\tau, \mu)$ is $\overline{L}_u:=(\iota,\mu_u)$, and similarly 
  for $(\tau,\mu')$. 
 Note that
  $\phi$, as a map  between the Hilbert spaces $(\iota, \mu_u)$
  and  $(\iota, \mu'_u)$ is isometric and in fact unitary, since $\phi$ is a  full functor. An easy computation shows that
  the collection of these unitaries satisfies the needed relations.
  \medskip

Our first result concerns the construction of induced bimodule representations.
 The   construction reduces to $({\cal C},\alpha)$  for $u=\iota$.

\medskip

\noindent{\bf 6.2. Theorem} {\sl Pick an object $u$ of ${\cal A}$ with $1_{\mu_u}\neq0$.
\begin{description}
\item{\rm a)} The linear
space $^\circ{\cal H}_u$ obtained quotienting
$\sum_v(\mu_v,\mu_u)\otimes \tau_v$ by the linear subspace
generated by elements of the form $M\circ\mu(A)\otimes
\psi-M\otimes\tau(A)\circ \psi$ can be naturally completed into a nonzero full
Hilbert
module ${\cal H}_u$ over ${\cal C}$, with a faithful left action of ${\cal C}$ making it into a Hilbert bimodule over ${\cal C}$.
${\cal H}_u$ depends only on $\mu_u$.
 
\item{\rm b)} There is a unique, full,
bimodule representation, $\text{Ind}(\mu_u)$, of $G_\tau$ on ${\cal
H}_u$ with
$$\text{Ind}(\mu_u) M\otimes\psi=M\otimes \hat{v}\psi,$$
for $M\in(\mu_v,\mu_u)$, $\psi\in\tau_v$, $\hat{v}$ being the representation of $G_\tau$ on
$\tau_v$.

  \end{description}

}\medskip

Theorem 6.2. will be  proved in Sections 7 and 8.
 In Sect.\ 8, we shall see that the map $\text{Ind}:\mu_u\to\text{Ind}(\mu_u)$ extends to a $^*$--functor $\text{Ind}:{\cal M}_\mu\to\text{Bimod}_\alpha(G_\tau)$ on the 
 full $C^*$--subcategory  ${\cal M}_{\mu}$
of ${\cal M} $ whose objects are those in the image of $\mu$. It will be called the {\it induction functor}.

In Sect.\ 8, we shall also prove the following  analogue of Swan's theorem in 
our framework. The assumptions are those of Theorem 6.2. $\text{Ind}\mu$ 
denotes the composed functor.
\medskip

\noindent{\bf 6.3. Theorem} {\sl 

\begin{description}

\item{\rm a)} For any object $u$ of ${\cal A}$
 there is a natural isometric  intertwiner of Hilbert module representations 
  $S_u\in( 
 \text{Ind}(\mu_u), \hat{u}\otimes\alpha)$.
In particular if $\mu$ is relaxed tensor, then $S_u$ is unitary.
Hence ${\cal H}_u$ is always  finite projective  as  a right module and free
if $\mu$ is relaxed  tensor.

  \item{\rm b)} For any arrow $A\in(u,v)$,  $\text{Ind}\mu(A)$ corresponds to the restriction of $\tau(A)\otimes I$ to the space of the associated subrepresentation of $\hat{u}\otimes\alpha$ under $S_u$.
  \end{description} }\medskip
  
  In other words, for  any arrow $A\in(u,u')$ in ${\cal A}$,  the following diagram commutes.\medskip 
  
\[
\xymatrix{
\text{Ind}(\mu_u)
   \ar[r]^-{\text{Ind}\mu(A)}
   \ar[d]_-{S_u }
& \text{Ind}(\mu_{u'})
   \ar[d]^-{ S_{u'}}
\\ \hat{u}\otimes{\alpha}
   \ar[r]^-{  \tau(A)\otimes I }
&  \hat{u'}\otimes\alpha}
\]

Note that, if $\mu$ is just quasitensor, ${\cal M}_{\mu}$
may not be a tensor category.
We  ask, however,
whether $\text{Ind}$
 extends {\it tensorially} to the smallest
full tensor subcategory ${\cal M}_\mu^\otimes$ of ${\cal M}$ generated by ${\cal M}_\mu$.
Somewhat surprisingly, the answer is that it does. Sections 9 and 10 will be devoted to discussing
the following result,  a generalization of Theorem 4.1 to a noncommutative framework.
\medskip

\noindent{\bf 6.4. Theorem} {\sl The induction functor $\text{Ind}: \mu_u\in{\cal M}_\mu\to\text{Ind}(\mu_u)\in\text{Bimod}_\alpha(G_\tau)$ 
extends uniquely to a full and faithful strict tensor
functor to a strict tensor category of Hilbert bimodule 
representations 
$$\text{Ind}: 
{\cal M}_\mu^\otimes\twoheadrightarrow \text{Bimod}_\alpha(G_\tau).$$  Furthermore,
$(\text{Ind}, \mu)$  gives rise to an adjoint pair of functors.}\medskip

 Since $\text{Ind}$ is a strict tensor functor,   the composed functor
$\text{Ind}\mu:{\cal A}\to\text{Bimod}_\alpha(G_\tau)$
is quasitensor, relaxed tensor or strict tensor according as $\mu$ is. 
The associated natural transformation $\widetilde{\text{Ind}\mu}$ is computed in Sect.\ 9.
Moreover, since $\text{Ind}$ is full, we may regard it as an isomorphism
between the original functor $(\mu,\tilde{\mu})$ and $(\text{Ind}\mu,
\widetilde{\text{Ind}\mu})$. In the following commutative diagram, the dotted arrows  summarize  
our construction from the given pair $(\tau,\mu)$.\medskip

\[
\xymatrix{
\text{Hilb}&{\cal A} \ar[l]_\tau\ar[r]^\mu\ar@{ .>}[rd]_{\text{Ind}\mu}&
{\cal M}_\mu^{\otimes} \ar@{ .>>}[d]^{\text{Ind}}\\
{}&{}&\text{Bimod}_\alpha(G_\tau)}
\]

\noindent{\it Remark} Combining Theorems 6.3 and 6.4 yields an explicit description
of the quasitensor functor $\mu:{\cal A}\to{\cal M}_\mu^\otimes$ in terms of $\tau$ and the ergodic action $({\cal C}, \alpha)$. This is then used for the embedding results, Theorems 6.7--6.9.

We next  give two applications of our results, that
originally motivated our work.
The first concerns a tensor $C^*$--category with conjugates whose object
set contains a
distinguished generating element. We give two results, corresponding to
the selfconjugate or non-selfconjugate case.
\medskip

\noindent{\bf 6.5. Theorem} {\sl Let ${\cal M}$ be a tensor $C^*$--category with objects
$\iota$, $x$, $x^2, \dots$, where $x$ is a real or pseudoreal object
defined by a solution $R\in(\iota, x^2)$ of the conjugate equations with
$\|R||^2\geq2$.
Let $F\in {\text Mat}_n({\mathbb C})$ be an invertible matrix with
$\text{Tr}(FF^*)=\text{Tr}((FF^*)^{-1})=\|R\|^2$. Then there is a
full and faithful
tensor functor ${\cal M}\to\text{Bimod}_\alpha(A_o(F))$, where
$\alpha$ is the ergodic action of
$A_o(F)$ on $\cC$ associated to $(\tau,\mu)$.}\medskip

\noindent{\bf 6.6. Theorem} {\sl If the set of objects of
${\cal M}$ is generated, as a semigroup, by $x$ and a conjugate $\overline{x}$,
with intrinsic dimension $d(x)\geq2$, then  
conclusions analogous to Theorem 6.5	 hold where the quantum group is now $A_u(F)$.}\medskip

\noindent{\it Examples}
Note that {\it any} spectral functor of an
ergodic action of a compact quantum group arises from some 
pair $(\tau,\mu)$, as we may choose for $\tau:\text{Rep}(G)\to\text{Hilb}$ the embedding functor and for 
$\mu$  the spectral functor of the action, $\mu:=\overline{L}:\text{Rep}(G)\to\text{Hilb}$.
(Recall that $\overline{L}$ is quasitensor by \cite{PR}, cf. Sect.\ 2.) 
On the other hand, many examples of noncommutative ergodic spaces are known to arise from 
pairs $(\tau,\mu)$, with $\tau$ as above, but where $\mu$ is tensorial or relaxed tensorial.
For example,  compact 
quantum quotients (completed in the maximal $C^*$--norm)
$C(K\backslash G)$   arise from the  restriction functor 
$\mu:\text{Rep}(G)\to\text{Rep}(K)$ \cite{PR}. The examples with high multiplicities 
of \cite{BDV} are associated with the composition $\mu$ of a tensorial isomorphism with the embedding functor, 
$\mu:\text{Rep}(G)\simeq\text{Rep}(G')\to\text{Hilb}$.   Examples of categories of the type described in 
Theorems 6.5 and 6.6 arise from  inclusions of $II_1$ factors $N\subset M$ with 
finite 
Jones index. The  ergodic action corresponding to the real object ${}_NM_N$ 
is made explicit in  \cite{PR2}. For any finite index inclusion  of 
 infinite 
factors described by an endomorphism $\rho$ with $d(\rho)\geq 2$,
the tensor $C^*$--category generated by $\rho$ and $\overline{\rho}$
is of the form described in Theorem 6.6.
\medskip

The proofs of Theorems 6.5 and 6.6 will be given at the end of    
Sect.\ 10.
The next application concerns tensor categories ${\cal M}$ extending the representation category of a compact group.
The following results, discussed in Sections 11 and 12, shed light on the problem of recognizing which
tensor
categories  can be embedded into Hilbert spaces. They are obtained  
combining our bimodule construction with the work of Takesaki
\cite{Takesaki}, H\o
egh--Krohn, Landstad and St\o rmer \cite{HLS} and Wassermann \cite{Wassermann3}.

In the following theorem $G$ is a compact Lie group, and we fix
a distinguished faithful representation $u$ such that  every irreducible of
$G$ is a subrepresentation of a tensor power of $u$. We denote by
${\cal S}_G$ be the full subcategory of $\text{Rep}(G)$ with objects $\iota$, $u$, $u^2,\dots$.
\medskip

\noindent{\bf  6.7.\ Theorem} {\sl Let $G$ be a compact Lie group with a distinguished faithful representation $u$, and let
$\mu:{\cal S}_G\to{\cal M}$ be a tensor functor. Let ${\cal
C}$ be the ergodic $C^*$--algebra associated with $\mu$ and the embedding
functor $\tau$ of ${\cal S}_G$ into the category of Hilbert
spaces. Assume that
the von Neumann algebra ${\cal C}''$  generated by ${\cal C}$ in the GNS
representation  of the unique $G$--invariant trace is of type $I$ and let
$K$ be a closed subgroup of $G$ such that $L^\infty(K\backslash G)\simeq Z({\cal C}'')$ as ergodic $W^*$--systems.
 Then
there is a full and faithful tensor functor $\epsilon: {\cal M}_\mu^\otimes\to\text{Rep}(K)$.}\medskip

Notice that in the above theorem ${\cal M}_\mu^\otimes$ is simply the full subcategory
of ${\cal M}$ with objects the tensor powers of $\mu_u$.
\medskip

\noindent{\it Remark} 
 As we shall see in Sect.\ 11, the functor $\epsilon$ is
naturally associated with $\mu$. However, the set of objects in the 
image of $\epsilon$ in general
does not generate all the representations of $K$ under tensor products,
subobjects and direct sums. In fact, in the particular case where each irreducible of $G$ has multiplicity
 in ${\cal C}$ equal to its dimension, then $\epsilon$ maps each object to the
trivial representation of $K$.  Hence ${\cal M}_\mu^\otimes$ admits a tensor functor 
to a full subcategory of the category of Hilbert spaces. Furthermore, any full multiplicity 
ergodic action of $G$ on a type
$I$ von Neumann algebra arises from a relaxed tensor functor $\mu$, its spectral functor.

At the other extreme, if ${\cal C}$ is commutative,  as happens in particular if
${\cal M}$ has permutation symmetry, see \cite{PR}, then we get
the following result, generalizing an important step in \cite{DR} towards proving 
the abstract duality theorem for compact groups. \medskip

\noindent{\bf 6.8.\ Corollary} {\sl  If ${\cal C}$ is commutative, and hence
${\cal C}=C(K\backslash G)$ as ergodic $C^*$--systems, then
$\epsilon(\mu_u)=u\upharpoonright_K$. Hence the completion of the image of 
$\epsilon$ under subobjects contains any irreducible of $K$.}\medskip

It is still an open problem 
whether ergodic actions
of  classical compact Lie groups $G$ on  von Neumann algebras are always of type $I$.
By Theorem 6.7, a positive answer for a specific group $G$ would guarantee the existence of an embedding for all tensor
$C^*$--categories ${\cal M}$ containing the representation category of $G$ and having the same objects.
Wassermann has shown this to be true for $G=SU(2)$ \cite{Wassermann3}. Taking  the known abstract characterization of $\text{Rep}(SU(2))$ into account \cite{DRduals}, see also \cite{Pinzari, PR} we
obtain the following embedding result for   tensor $C^*$--categories 
containing a distinguished pseudoreal object of dimension $2$. No permutation symmetry is assumed.
\medskip

\noindent{\bf 6.9. Corollary} {\sl Let ${\cal M}$ be a tensor
$C^*$--category
whose object semigroup is generated by a pseudoreal object $x$ of
dimension $2$, i.e.
with an intertwiner $R\in(\iota, x^2)$
such that
$$R^*\otimes 1_x\circ 1_x\otimes R=-1_x,$$
$$\|R\|^2=2.$$
Then there is a closed subgroup $K$ of $SU(2)$  and
a full and faithful  tensor  functor ${\cal M}\to\text{Rep}(K)$ .}\medskip

 \noindent{\it Remarks and more results}.
A non-trivial problem is to construct new examples of or even classify the quasitensor functors $\mu:{\cal A}\to{\cal M}$, for a given embeddable tensor $C^*$--category ${\cal A}$.
 Our results connect 
this problem to that of classifying the ergodic $C^*$--actions of  the quantum group $G_\tau$ associated to an embedding $\tau$ of ${\cal A}$.   On one hand,
as recalled in the examples previously discussed, every spectral functor of an  ergodic action of $G_\tau$ on a unital $C^*$--algebra  arises in this
way.   Even for ergodic $C^*$--actions of compact groups, where there are  
  important results, not a lot is 
 known (cf. Sect.\ 2). In the quantum case very little is known, but it is already clear that there are  many new aspects. 

Motivated by
 our applications, we are especially interested in the case  where $\mu$ is tensorial or relaxed tensorial. 
The
reconstruction results, Theorems 6.2 and 6.4, then 
lead to the problem of classifying
those ergodic actions $({\cal C},\alpha)$ where 
 the module representations $u\otimes\alpha$ (or a subrepresentation on a projective module in the quasitensor case) can be made  into full bimodule representations. 

Not all ergodic actions, even of compact groups can arise in this way. 
In Sect.\ 11 we 
classify full bimodule representations arising from ergodic actions of compact groups on type $I$ von Neumann algebras. This provides an obstruction to the existence 
of full bimodule representations (and hence to the existence of relaxed tensor functors
$\mu:{\cal A}\to{\cal M}$) in the case of low but nonzero multiplicities. For example, we derive that neither
the ergodic actions with full spectrum and irreducibles of low multiplicity nor
  the adjoint action 
  by a non-trivial irreducible representation of 
 $SU(2)$ can  arise.

\bigskip

\section{Algebraic bimodules from pairs of functors}

As in the previous section, we start with
  tensor $C^*$--categories ${\cal A}$ and ${\cal M}$ and we assume that ${\cal A}$ has 
conjugates. 
Let
$(\mu,\tilde{\mu}):{\cal A}\to{\cal M}$ be  a quasitensor functor and
 $\tau: {\cal A}\to\text{Hilb}$ a tensor functor into the category of Hilbert spaces.
We have an associated  
unital $C^*$--algebra ${\cal C}$, the completion of a canonical dense 
$^*$--subalgebra ${}^\circ{\cal C}$.
In this section we 
generalize  that construction at the algebraic level to get 
bimodules over  ${}^\circ{\cal C}$. The norm completion and the quantum group action will be considered in the next section.
\medskip

\noindent{\it 7.1. The algebraic bimodules ${}^\circ{\cal H}_u$.}  Pick an object  $u$ of ${\cal A}$. 
Let ${}^\circ{\cal H}_u$ be the linear space
$\sum_v(\mu_v,\mu_u)\otimes\tau_v$, the sum being taken over the objects of ${\cal A}$, quotiented by the linear subspace generated by elements of the form $$M\circ\mu(A)\otimes \psi-M\otimes\tau(A)\circ \psi.$$
Notice that, as the objects involved have conjugates
\cite{PR}, \cite{Ergodic}, and tensor units are irreducible, we are actually taking a sum of finite dimensional
 vector spaces \cite{LR}. 
 It should be noted that the bimodule ${}^\circ\cH_u$ in 
fact depends only on $\mu_u$. 
\medskip

We next introduce a multiplication and adjoint, $$\cdot{}: {}^\circ{\cal H}_u
\times{}^\circ{\cal H}_{u'}\to{}^\circ{\cal H}_{u\otimes u'}$$
$$^*:{}^\circ{\cal H}_u\to{}^\circ{\cal H}_{\overline{u}},$$ 
to get a structure analogous to a $^*$--algebra. These operations will be used in Sect. 8  to simplify 
computations, the  multiplication $\cdot{}$ also plays a role  in Sect. 9.
\medskip

\noindent{\it 7.2. The  multiplication  $\xi\cdot{}\eta$.} For simple tensors $\xi=L\otimes \psi\in{}^\circ{\cal H}_u$, $\eta=M\otimes \phi\in{}^\circ{\cal H}_{u'}$, with $L\in(\mu_w, \mu_u)$, 
$M\in(\mu_v, \mu_{u'})$, 
$\psi\in  \tau_w$, $\phi\in \tau_v$,
set:
$$\xi\cdot\eta:=\tilde{\mu}_{u,u'}\circ(L\otimes 
M)\circ\tilde{\mu}_{w,v}^*\otimes (\psi\otimes \phi).$$
It is easy to check that these maps are well defined and associativity 
follows from that of the functors $\mu$ and $\tau$.
In particular, ${}^\circ{\cal H}_\iota$ is an algebra, 
denoted above by 
$^\circ{\cal C}$.
Note that the multiplication depends in general on $u$ and $u'$ and not only on $\mu_u$, $\mu_{u'}$. 
 However, if $u'$ (or $u$) is the tensor unit, it depends only on $\mu_u$ 
 (or $\mu_{u'}$) as $\tilde{\mu}_{u,\iota}=\tilde{\mu}_{\iota, u}=1_{\mu_u}$. Hence  
$${}^\circ{\cal H}_u \,\, \text{ is a }\,\,
^\circ{\cal C}\text{--bimodule depending just on} \,\, \mu_u.$$\medskip

\noindent{\it 7.3. The functor $\lambda$.}
We define a functor, denoted by $\lambda$ for the moment, from ${\cal M}_\mu$ to the category of 
${}^\circ{\cal C}$--bimodules. After the norm and the quantum group action have been introduced, 
$\lambda$ will be the induction functor $\text{Ind}$.

Given 
$Y\in(\mu_u,\mu_{u'})$, we  define a map
$$\lambda(Y):{}^\circ\cH_{u}\to{}^\circ\cH_{u'},\quad \lambda(Y)(M\otimes \psi):=(Y\circ M)\otimes \psi.$$
It is easily checked that $\lambda(Y)$ is a bimodule map so $\lambda$ is 
a covariant functor,  from  the full subcategory ${\cal M}_\mu$ of ${\cal M}$   whose 
objects are the images of objects of ${\cal A}$ into the category of 
$^\circ{\cal C}$--bimodules.
\medskip

\noindent{\it 7.4. The adjoint $\xi^*$.}
We next   define an 
adjoint on these bimodules. Here matters are slightly more complicated. 

As recalled in subsect.\ 2.2, if $(\mu,\tilde{\mu})$ is a quasitensor functor and if $R$, $\overline{R}$ defines a conjugate for an object $u$, then 
$\hat{R}:=\tilde{\mu}_{\overline{u},u}^*\circ\mu(R)$
and 
$\hat{\overline{R}}:=\tilde{\mu}_{u,\overline{u}}^*\circ\mu(\overline{R})$
is a solution of the conjugate equations for $\mu_u$, the image solution
of $R$, $\overline{R}$
 under $\mu$ \cite{PR}.

\medskip

Fixing an object $u\in{\cal A}$ and a solution $R$, $\overline{R}$ of the conjugate 
equations for  $u$, we associate an antilinear map $^*:{\cal H}_u\to{\cal H}_{\overline{u}}$ in the following way.
Choose solutions of the conjugate equations $v\to R_v,\overline{R}_v$ of  for all  objects of ${\cal A}$.
Set, $$(M\otimes \psi)^*:=M^\bullet\otimes j_v\psi,$$
for $M\in(\mu_v,\mu_u)$, $\psi\in\tau_v$,
 where  $^\bullet:(\mu_v,\mu_u)\to(\mu_{\overline{v}},\mu_{\overline{u}})$ 
 is defined using image solutions under $\mu$  of the chosen solutions 
 $R_v,\overline{R}_v$ for the running objects $v$ appearing in the sum and $R$, $\overline{R}$  for the fixed object $u$  respectively and  $j_v:=j_{\tau_v}$ corresponds to $\tau(R_v)$ and $\tau(\overline{R}_v)$ as in subsect. 2.1.  Notice that
 $^*$ is well defined
  by the compatibility properties with $\circ$, $\mu$ and $\tau$, see subsect. 2.1 and 2.2. Moreover, ${}^*$ is independent of the choice $v\to R_v, \overline{R}_v$ 
in ${\cal A}$ for the running objects $v$, as if $Y\in(\overline{v},\tilde{v})$ is an invertible,  $M^\bullet$ and $j_v\psi$ become  $M^\bullet\circ\mu(Y^{-1})$ and $\tau(Y)j_v\circ\psi$ respectively. However, 
if we change the solution of the conjugate equations for $u$  using an
$X\in(\overline{u},\tilde{u})$, $(M\otimes \psi)^*$ becomes 
$(\mu(X)M^\bullet)\otimes 
(j_v\psi)=\lambda(\mu(X))(M\otimes 
\psi)^*$, hence the associated $^*$ changes.  This unpleasant feature will play no role in the construction of the bimodule representation.

We note that for  $u=\iota$ the $^*$--operation is independent and makes 
${}^\circ{\cal C}$ into a unital $^*$--algebra. \medskip

\noindent{\it 7.5. Compatibility of the various operations.}
\medskip

\noindent{\bf 7.1.\ Proposition} {\sl Let $u$, $u'$, 
be objects of ${\cal A}$. 
If $\xi\in {}^\circ{\cal H}_{u}$ and $\xi'\in {}^\circ{\cal 
H}_{u'}$,
$(\xi\cdot\xi')^*={\xi'}^*\cdot\xi^*$ and 
$\xi^{**}=\xi$,
where we have used
tensor product solutions of the conjugate 
equations for 
$u\otimes u'$ 
and conjugate solutions for $\overline{u}$.

}\medskip

\noindent{\it Proof} Write $\xi:=L\otimes\psi$ and $\xi':=M\otimes\phi$ with $L\in(\mu_w,\mu_u)$, $\psi\in\tau_w$, $M\in(\mu_v,\mu_{u'})$, $\phi\in\tau_{v}$. We may compute
$(\xi\cdot{}\xi')^*$ using image of a tensor product solution of the conjugate equations for $w$ and $v$.
 The first result follows from 
$$(\mu_{u,u'}\circ L\otimes M\circ\tilde\mu_{w,v}^*)^\bullet=\mu_{\overline{u'},\overline{u}}\circ M^\bullet\otimes L^\bullet\circ\tilde\mu_{\overline v,\overline w}^*,$$
see subsect.\ 2.1 for the compatibility properties of $\bullet$ with $\otimes$ and $\circ$ and Cor.\ 13.3 for the explicit computations  of ${\mu_{w,v}^*}^\bullet$ and $\mu_{u,u'}^\bullet$.
 The second result follows 
from $M^{\bullet\bullet}=M$ when a solution of the conjugate equation 
and successively the conjugate solution is used.
\medskip

\noindent{\it Remark} Some care is required in using this proposition. 
For example, if $u$ is pseudoreal and irreducible and we use $R\in(\iota, 
u^2)$ to define 
$\xi^*$ in ${\cal H}_u$, we must use $-R$ to define $\xi^{**}$ in ${\cal 
H}_u$.
\medskip


\noindent{\bf 7.2.\ Proposition} {\sl
For 
$Y\in(\mu_u, \mu_{u'})$, $\xi\in 
{}^\circ{\cal 
H}_u$, $\xi'\in{}^\circ{\cal H}_{u'}$, $A\in(u,z)$, $A'\in(u', z')$,
$$(\lambda(Y)\xi)^*=\lambda(Y^\bullet)\xi^*,$$
$$\lambda(\mu(A\otimes 
A'))\xi\cdot\xi'=\lambda(\mu(A))\xi\cdot\lambda(\mu(A'))\xi'.$$}\medskip
\medskip

The last item to be introduced in this section is a sesquilinear form on the bimodules
${}^\circ{\cal H}_u$.  It will be shown to be  positive in the next section, allowing us
to pass from the algebraic to the analytic level.
\medskip

\noindent{\it 7.6. The sesquilinear form on the bimodules ${}^\circ{\cal H}_u$.} We retain the notation of subsect.\ 7.4 and define a
$^\circ{\cal C}$--valued form on 
${}^\circ{\cal H}_u$ by setting
$$<\xi,\xi'>:=\lambda(\mu(R)^*)(\xi^*\cdot\xi').\eqno(7.1)$$
The explicit formula, for $\xi=M\otimes \psi$, $\xi'=M'\otimes \psi'$,
with $M\in(\mu_v,\mu_u)$, $\psi\in\tau_v$, $M'\in(\mu_{v'},\mu_u)$,
$\psi'\in\tau_{v'}$ is
$$<\xi,\xi'>:=(\hat{R}^*\circ M^\bullet\otimes 
M'\circ\tilde{\mu}^*_{\overline{v}, 
v'})\otimes
(j_v \psi\otimes 
\psi').$$

\noindent{\it Remark} This form does not depend on the chosen solution of the conjugate equations for $u$. Indeed,
if we change solution   using an invertible $X$ then $M^\bullet$ becomes 
$\mu(X)\circ M^\bullet$, and this cancels the simultaneous change of $\hat{R}^*$, which becomes $\hat{R}^*\circ\mu(X^{-1})\otimes 1_{\mu_u}$.

We conclude this section with an explicit computation
of the right hand side needed later.
\medskip

\noindent{\bf 7.3. Lemma} {\sl For $\xi=M\otimes\psi$, $\xi'=M'\otimes\psi'$, $M\in(\mu_v,\mu_u)$, 
$M'\in(\mu_{v'},\mu_u)$, $\psi\in \tau_v$, $\psi'\in\tau_{v'}$
$$<\xi,\xi'>=(\hat{R_v}^*\circ
1_{\mu_{\overline{v}}}\otimes (M^*\circ M')\circ\tilde{\mu}_{\overline{v}, 
v'}^*)\otimes 
(j_v\psi\otimes \psi').$$
In particular, the form depends only on $\mu_u$.}\medskip

\noindent{\it Proof}
$$\hat{R}^*\circ M^\bullet\otimes M'=\hat{R}^*\circ[(\hat{R_v}^*\otimes 
1_{\mu_{\overline{u}}}\circ 1_{\mu_{\overline{v}}}\otimes M^*\otimes 
1_{\mu_{\overline{u}}}\circ1_{\mu_{\overline{v}}}\otimes
\hat{\overline{R}})\otimes 
M']=$$
$$\hat{R_v}^*\otimes \hat{R}^*\circ 1_{\mu_{\overline{v}}}\otimes 
M^*\otimes 1_{\mu_{\overline{u}}}\otimes 1_{\mu_u}\circ
1_{\mu_{\overline{v}}}\otimes \hat{\overline{R}}\otimes 1_{\mu_u}\circ 
1_{\mu_{\overline{v}}}\otimes M'=$$
$$\hat{R_v}^*\circ 1_{\mu_{\overline{v}}}\otimes (M^*\circ M').$$
We have already noted that  the form does not depend on the solution of the conjugate equations for $v$   and $u$ and see now
that it does not change if we replace $u$ by another object $u'$ such that $\mu_u=\mu_{u'}$.

\bigskip

\section{The induced Hilbert bimodule representations}

In this section we consider both the analytic aspect of the bimodules 
${}^\circ{\cal H}_u$ and the quantum group action, leading to 
 the proof of Theorems 6.2  and 6.3.  
 
Recall that $G_\tau$ was defined in Sect.\ 6 as the compact quantum group defined by the functor $\tau$ via Woronowicz duality. Let 
$\alpha$ be  the action  of $G_\tau$ on   ${}^\circ{\cal C}$
defined by $\alpha(M\otimes \psi)=M\otimes 
\hat{v}(\psi)$ for $M\in(\mu_v,\iota)$, $\psi\in\tau_v$, where  
$\hat{v}$ denotes the representation of $G_\tau$ on $\tau_v$.
${}^\circ{\cal C}$ is known to have a maximal $C^*$--norm and $\alpha$ to extend uniquely to  an ergodic action of $G_\tau$ on the completed $C^*$--algebra ${\cal C}$ \cite{PR}.



\medskip

\noindent{\it 8.1. Positivity of the sesquilinear form.}  Given objects $u$, $v\in{\cal A}$, $\phi\in\tau_v$, 
we let
$$L_{u}(\phi):{}^\circ{\cal H}_u\to{}^\circ{\cal H}_{v\otimes u}$$ be the operator of left multiplication by   $1_{\mu_v}\otimes\phi\in{}^\circ{\cal H}_v$ on $^\circ\cH_u$, hence obviously a right module morphism. 
Now set
$$L_{u}(\phi)^*:{}^\circ{\cal H}_{v\otimes u}\to{}^\circ{\cal H}_u,$$
$$L_{u}(\phi)^*:=\lambda(\mu(R^*_{v}\otimes 1_u))L_{v\otimes u}(j_v\phi).$$ 
 If we change solutions of the conjugate equations using an invertible $X\in(\overline{v},\tilde{v})$,  it is not difficult to verify that, by Lemma 8.2 c) below, $L_{u}(\phi)^*$ does not change.
\medskip

\noindent
{\bf 8.1.\  Lemma }  {\sl $<\eta,L_{u}(\phi)^*\xi>=<L_{u}(\phi)\eta,\xi>$.}\medskip  

\noindent
{\it Proof} $$<\eta,L_{u}(\phi)^*\xi>=\lambda(\mu(R_u^*))(\eta^*\cdot(\lambda(\mu(R_v^*\otimes 1_u))(1_{\mu_v}\otimes \phi)^*\cdot\xi))=$$ 
$$\lambda(\mu(R_u^*\circ 1_{\overline u}\otimes R_v^*\otimes 1_u))\eta^*\cdot(
1_{\mu_v}\otimes\phi)^*\cdot\xi=<L_{u}(\phi)\eta,\xi>.$$
The second equality follows from Prop. 7.2 while in the last   we have  chosen product solutions of the conjugate equations  for $v\otimes u$. 
 \medskip

As we shall soon see, the $\cC$--valued form $<\cdot,\cdot>$ is 
positive, so that $L_{u}(\phi)^*$ is the adjoint of $L_{u}(\phi)$ as the 
notation suggests.
These maps satisfy the following properties.\medskip

\noindent{\bf 8.2.\  Lemma} {\sl 
\begin{description}

\item {\rm a)}
$\lambda(\tilde\mu_{r,w}\circ 1_{\mu_r}\otimes Y\circ\tilde{\mu}_{r,u}^*)L_{u}(\phi)=L_{w}(\phi)\lambda(Y), \quad Y\in(\mu_u,\mu_w),$
 \item  {\rm b)} $L_{u}(\phi)^*L_{u}(\psi)=<\phi,\psi>$
\item{\rm c)} $\lambda\mu(A\otimes 1_u)L_u(\phi)=L_u(\tau(A)\phi),$ 
\item{\rm d)}
$\sum_i 
L_u(\phi_i)L_u(\phi_i)^*=\lambda(\tilde{\mu}_{z,u}\circ\tilde{\mu}_{z,u}^*),$
where $\phi_i$ is  an orthonormal basis of $\tau_z$,
\end{description}

}\medskip

For $u=\iota$  the corresponding  maps  
${}^\circ{\cal C}\to{}^\circ{\cal H}_z$,
${}^\circ{\cal H}_z\to{}^\circ{\cal C}$
will simply be  denoted $L(\phi)$ and $L(\phi)^*$.
By d), we have:
$$\sum L(\phi_i)L(\phi_i)^*=1$$
 for any orthonormal basis $(\phi_i)$ of $\tau_z$.
We shall use this relation   to define a faithful right module map 
$S_z:{}^\circ{\cal H}_z\to\tau_z\otimes{}^\circ{\cal C}$
by
$$S_z\xi:=\sum_i\phi_i\otimes L(\phi_i)^*\xi,$$
clearly independent of the choice of the orthonormal basis.

We are now ready to show positivity of the sesquilinear form of ${}^\circ{\cal H}_u$.
\medskip

\noindent{\bf 8.3.\ Proposition} {\sl If $\tau_z\otimes ^\circ{\cal C}$ is considered as a  right prehilbertian $^\circ{\cal C}$--module, the map $S_z$ satisfies:
$$<S_z\xi, S_z\xi'>=<\xi,\xi'>,\quad \xi,\xi'\in{}^\circ{\cal H}_z.$$
Hence $^\circ{\cal H}_z$ is a finite projective,  right prehilbertian 
module over ${}^\circ{\cal C}$ with the sesquilinear form defined in $(7.1)$ and $S_z$ is an isometric right 
$^\circ{\cal C}$--module 
map. Its  adjoint $S_z^*:
\tau_z\otimes 
{}^\circ{\cal C}
\to
{}^\circ{\cal H}_z 
$ is given by
$$S_z^*(\psi\otimes I)=1_{\mu_z}\otimes \psi$$
for $\psi\in\tau_z$.}\medskip

\noindent{\it Proof} 
$$<S_z\xi,S_z\xi'>=\sum_i<L(\phi_i)^*\xi,L(\phi_i)^*\xi'>=\sum_i<\xi,L(\phi_i)L(\phi_i)^*\xi'>=<\xi,\xi'>.$$ 
Hence $<\cdot,\cdot>$ is a faithful, positive, $^\circ{\cal C}$--valued 
inner product on $^\circ{\cal H}_z$ and $S_z$ an isometry.
We next compute the adjoint of $S_z$. If $\xi\in{}^\circ{\cal H}_z$,
$\psi\in \tau_z$ then 
$$<\xi, S_z^*\psi\otimes I>=<S_z\xi, \psi\otimes I>=\sum_i<\phi_i\otimes 
L(\phi_i)^*\xi, \psi\otimes I>=$$ 
$$\sum_i(\phi_i,\psi)<\xi,1_z\otimes\phi_i>=<\xi,1_z\otimes\psi>,$$ 
as required.

\medskip

We next compute the range projection $P_z=S_zS_z^*$ 
to see when $S_z$ is unitary.
If $\psi\in \tau_z$,
$$P_z(\psi\otimes I)=S_z(1_{\mu_z}\otimes\psi)=\sum _i \phi_i\otimes 
L(\phi_i)^*(1_{\mu_z}\otimes\psi)=$$
$$\sum_i\phi_i\otimes((\mu(R_z^*)\circ
\tilde{\mu}_{\overline{z}, 
z}\circ\tilde{\mu}^*_{\overline{z}, z}
)\otimes(j_z\phi_i\otimes\psi)).$$

\noindent{\bf 8.4.\ Corollary} {\sl If 
$\tilde{\mu}_{\overline{z}, 
z}\circ\tilde{\mu}^*_{\overline{z}, z}\circ \mu(R_z)=\mu(R_z)$ (e.g.\ 
when $\mu$ is relaxed tensor),  $S_z$ is 
unitary, and hence $^\circ{\cal H}_z$ is 
a free right $^\circ{\cal C}$--module.}\medskip

\noindent{\it Proof}
We now have
$$\mu(R_z^*)\circ
\tilde{\mu}_{\overline{z}, 
z}\circ\tilde{\mu}^*_{\overline{z}, z}
\otimes(j_z\phi_i\otimes\psi)=\mu(R_z^*)\otimes(j_z\phi_i\otimes\psi)=
1_\iota\otimes\tau(R_z^*)j_z\phi_i\otimes\psi=<\phi_i,\psi>.$$
Hence $P_z$ is the identity projection.\medskip

We conclude the subsection noting that
property c) of Lemma 8.2 implies that $S$ is a natural transformation from $\lambda\mu$ to $\tau\otimes I$, i.e., when $A\in(u,u')$, the following diagram commutes:
\[
\xymatrix{
{}^\circ{\cal H}_u
   \ar[r]^-{\lambda\mu(A)}
   \ar[d]_-{S_u }
&  {}^\circ{\cal H}_{u'}
   \ar[d]^-{ S_{u'}}
\\ \tau_u\otimes{}^\circ{\cal C}
   \ar[r]^-{  \tau(A)\otimes I }
&   \tau_{u'}\otimes{}^\circ{\cal C}}
\]
This will be used at the end of the section when proving Theorem 6.3,to the generalized
form of Swan's Theorem.
\medskip

\noindent{\it  8.2. $^\circ{\cal H}_u$ is algebraically 
full.}\medskip

\noindent{\bf 8.5.\ Proposition} {\sl Let $u$ be an object of ${\cal A}$ with 
$1_{\mu_u}\neq0$, then
the coefficients $<\xi,\xi'>$, $\xi$, $\xi'\in ^\circ{\cal H}_u$, span
$^\circ{\cal C}$.}\medskip

\noindent{\it Proof} Choose  
$v=u$, $v'=u\otimes v''$, $M= 1_{\mu_u}$, $M'=M''\circ\tilde{\mu}^*_{u, 
v''}$ with $M''\in(\mu_u\otimes\mu_{v''},\mu_u)$, $\psi=j_u^{-1}\phi_i$,
$\psi'=j_u^{-1}\phi_i\otimes\psi''$ in $(7.1)$ where $(\phi_i)$ is an orthonormal basis 
of 
$\tau_{\overline{u}}$ and $\psi''\in\tau_{v''}$.
Summing over $i$ gives 
$$(\hat{R}^*_u\circ 1_{\mu_{\overline{u}}}\otimes 
(M''\circ\tilde{\mu}^*_{u, v''})\circ \tilde{\mu}^*_{\overline{u}, 
u\otimes v''})\otimes (\tau(R_u)\otimes\psi'')=$$
$$(\hat{R}^*_u\circ 1_{\mu_{\overline{u}}}\otimes 
(M''\circ\tilde{\mu}^*_{u, v''})\circ \tilde{\mu}^*_{\overline{u}, 
u\otimes v''}\circ\mu(R_u\otimes 1_{v''}))\otimes \psi''=$$
$$(\hat{R}^*_u\circ 1_{\mu_{\overline{u}}}\otimes 
(M''\circ\tilde{\mu}^*_{u, v''})\circ \tilde{\mu}^*_{\overline{u}, 
u\otimes v''}\circ\tilde{\mu}_{u\otimes\overline{u}, v''} 
\circ\mu(R_u)\otimes 1_{\mu_{v''}})\otimes \psi''=$$
$$(\hat{R}^*_u\circ 
1_{\mu_{\overline{u}}}\otimes 
M'' 
\circ\hat{R}_u\otimes 1_{\mu_{v''}}
)\otimes \psi''.$$
Now recall, see e.g.\ \cite{LR} that if $\rho$, $\sigma$, $\tau$ are objects 
of a tensor $C^*$--category with conjugates,
the map $$T\in(\rho\otimes\sigma,\tau)\to 1_{\overline{\rho}}\otimes 
T\circ 
R_\rho\otimes 1_\sigma\in(\sigma, \overline{\rho}\otimes\tau)$$
is a linear isomorphism.
Hence 
$X:=1_{\mu_{\overline{u}}}\otimes 
M'' \circ\hat{R}_u\otimes 1_{\mu_{v''}}$ is a generic element of 
$(\mu_{v''}, \mu_{\overline{u}}\otimes\mu_u)$ and can, in particular,
be any element of the form $X=\hat{R}_u\circ Y$ where $Y\in(\mu_{v''},\iota)$. 
Hence the linear span of
the coefficients of the inner product on $^\circ{\cal H}_u$ is $^\circ{\cal C}$
as it contains any element of the form $Y\otimes\psi''$.\medskip

\noindent{\it 8.3. A useful formula for  the left $^\circ{\cal C}$--action on 
$^\circ{\cal H}_u$.} 
As 
$x_i:=S_u^*(\psi_i\otimes I)=1_{\mu_u}\otimes\psi_i$, where 
$(\psi_i)$ is an orthonormal basis of $\tau_u$,
is a Hilbert module basis, we need only specify 
$<x_i, c\cdot x_j>$ for $c\in{}^\circ{\cal C}$.
\medskip

\noindent{\bf 8.6.\  Proposition} {\sl
If $c=T\otimes \phi\in{}^\circ{\cal C}$, with 
$T\in(\mu_v,\iota)$, $\phi\in\tau_v$,
$$<x_i, c\cdot x_j>=(\hat{R}^*_u\circ 1_{\mu_{\overline{u}}}\otimes 
T\otimes 1_{\mu_u}\circ\tilde{\mu}^*_{\overline{u}, v,u})\otimes 
(j_u\psi_i\otimes\phi\otimes \psi_j).$$ 
}\medskip

\noindent{\it Example} Let $K\subset G$ be an inclusion of a compact 
quantum groups. Then we have a tensor functor
$\mu:\text{Rep}(G)\to\text{Rep}(K)$ taking a representation 
$u$ of 
$G$ to its restriction $u\upharpoonright_K$ to the subgroup, see, e.g., \cite{Pinzari3}.
Hence 
 $^\circ{\cal H}_u$ is a free module.
 $^\circ{\cal C}$ is the 
canonical dense $^*$--subalgebra of $C(K\backslash G)$ . 
The above formula then gives:
$$<x_i, c\cdot x_j>=\sum_r u_{ri}^*cu_{rj},\quad c\in{}^\circ C(K\backslash 
G),$$ 
see \cite{PR}.
This example was discussed at the end of Sect.\ 5.
\medskip

\noindent{\it Example } Let  $G$ be a compact quantum group acting 
ergodically on a unital $C^*$--algebra ${\cal C}$, and let 
$\overline{L}:\text{Rep}(G)\to\text{Hilb}$ 
be the spectral functor of the action as in \cite{PR}, and shown 
there to be a quasitensor functor. Then $^\circ{\cal C}$ 
is the $^*$--algebra spanned by the elements of ${\cal C}$ transforming 
under the action like 
unitary irreducible $G$--representations.  
$^\circ{\cal 
H}_u\neq0$ precisely when $\overline{L}_u\neq 0$ and, if $u$ is irreducible, 
this is equivalent to requiring $u$ to lie in the spectrum of the action.
We thus get a finite projective 
$^\circ{\cal C}$--bimodule $^\circ{\cal H}_u$. Computations similar to those in
the above  example, show that the left $^\circ{\cal 
C}$--action is given by 
$$<x_i, c\cdot x_j>=\sum_r{c^u}_{r,i}^* c c^u_{r,j}$$
where $c^u_r:=(c^u_{rj})_j\in\overline{L_u}$ is an  orthonormal basis of $\overline L_u$.
\medskip

\noindent{\it Remark}
The restriction functor and  the spectral functor of a quantum quotient 
define the same 
algebra $^{\circ}{\cal C}$. However,   the associated 
bimodules  are different in general as, in the first case, they are free and never zero for 
a nonzero object, whilst in the second,  non-spectral representations give zero bimodules. 
\medskip

\noindent{\it 8.4. The completed Hilbert bimodules ${\cal H}_u$.}
Consider $^\circ{\cal C}$ with 
its maximal $C^*$--norm, which is finite by \cite{PR}. 
Completing 
$^\circ{\cal H}_u$ in the norm $\|\xi\|:=
\|<\xi,\xi>\|^{1/2}$, gives a right Hilbert
module ${\cal H}_u$ over the completion ${\cal C}$ of $^\circ{\cal C}$.
There is an isometry ${\cal H}_u\to\tau_u\otimes{\cal C}$ 
extending the algebraic isometry $S_z$. Hence 
${\cal H}_u$ is a finite projective right Hilbert ${\cal C}$--module.
Consequently , every right module map on $^\circ{\cal H}_u$ extends to an adjointable bounded map on ${\cal H}_u$. Hence the left $^\circ{\cal C}$--action extends to a unital
$^*$--homomorphism ${\cal C}\to{\cal L}_{\cal C}({\cal H}_u)$
thus making ${\cal H}_u$ into a Hilbert ${\cal C}$--bimodule.

To show that the left action is faithful 
we need norm continuity of the multiplication of bimodules.\medskip

\noindent{\bf 8.7.\ Proposition} {\sl The multiplication map  
$$ \xi\otimes \xi'\in{}^\circ{\cal H}_u\otimes_{^\circ{\cal C}}{^\circ{\cal H}_{u'}}\to\xi\cdot\xi' \in{}^\circ{\cal H}_{u\otimes u'}$$
is an isometry of prehilbertian $^\circ{\cal C}$--bimodules. 
Hence $\|\xi\cdot\xi'\|\leq\|\xi\|\|\xi'\|$.}\medskip

\noindent{\it Proof} Using a product solution of the conjugate equations,
$$<\xi\cdot\xi', \eta\cdot\eta'>=\lambda(\mu(R_{u\otimes 
u'}^*))\xi'^*\cdot\xi^*\cdot\eta\cdot\eta' 
=\lambda(\mu(R_{u'}^*))\circ\lambda(\mu(1_{\overline{u'}}\otimes 
R_u^*\otimes 1_{u'}))\xi'^*\cdot\xi^*\cdot\eta\cdot\eta'=$$
$$(\lambda(\mu(R_{u'}^*))
\xi'^*)\cdot
(\lambda(\mu(R_u^*))\xi^*\cdot\eta)\cdot\eta'=<\xi', 
<\xi,\eta>\cdot\eta'>=<\xi\otimes\xi', \eta\otimes\eta'>,$$
as required.
\medskip

Consequently, $\cdot$ extends to an associative multiplication 
$\xi\cdot\xi'$
on the completed bimodules ${\cal H}_u$ and ${\cal H}_{u'}$.\medskip

\noindent{\bf 8.8. Proposition} {\sl The extended left action of ${\cal C}$ 
on 
${\cal H}_u$ is faithful whenever $1_{\mu_u}\neq0$.}\medskip

\noindent{\it Proof} If $c\cdot\xi=0$ for all $\xi\in{\cal H}_u$, 
then
$$c\cdot(\lambda(\mu(\overline{R}^*_u))\xi\cdot\eta)=
\lambda(\mu(\overline{R}^*_u))c\cdot\xi\cdot\eta=0,$$
for all $\eta\in{\cal H}_{\overline{u}}$, On the other hand, 
$\lambda(\mu(\overline{R}^*_u))\xi\cdot\eta=<\xi^*,\eta>$, and these 
coefficients span $^\circ{\cal C}$ if 
$1_{\mu_{\overline{u}}}\neq0$
i.e. if $1_{\mu_{u}}\neq0$.\medskip

\noindent{\it 8.5. Quantum group representations on ${\cal H}_u$.}
We next construct quantum group representations on the bimodules ${\cal H}_u$.
Let $G_\tau$ denote, as before, the Woronowicz dual of $\tau:{\cal A}\to\text{Hilb}$.
\medskip

\noindent{\bf 8.9. Proposition} {\sl Given an object $u$ of ${\cal A}$, there is a 
unique bimodule representation 
$\text{Ind}(\mu_u)$ of $G_\tau$
on ${\cal H}_u$ such that $$\text{Ind}(\mu_u)(M\otimes \psi)=M\otimes 
\hat{v}(\psi),$$
$M\in(\mu_v,\mu_u)$, $\psi\in(\iota,\tau_v)$. $\text{Ind}(\mu_u)$ is 
a full bimodule representation.
}\medskip

\noindent{\it Proof}
The relation between the invertible antilinear maps $j_v: \tau_v\to 
\tau_{\overline{v}}$ and the coefficients of the corresponding 
representations of $G_\tau$ is given by
$\hat{\overline{v}}j_v\psi=\sum\phi_i\otimes 
\hat{v}^*_{j_v^*\phi_i,\psi}$, 
where 
$\phi_i$ an orthonormal basis of $\tau_{\overline{v}}$. This relation 
together with $(7.1)$ 
allows us to verify $(3.1)$. $(3.2)$, $(3.3)$ and $(3.5)$ 
follow from straightforward computations, whilst $(3.4)$ is a consequence 
of the corresponding relation for the Hilbert space representation 
$\hat{v}$.
We show that $\text{Ind}(\mu_u)$ is a full representation.  
A $G_\tau$--fixed vector $\xi$ in ${\cal H}_u$ for the underlying module 
representation is a simple tensor of the
form $\xi=T\otimes 1_\iota$, $T\in(\iota, \mu_u)$.
For any $c\in {\cal C}$  
of the form $c=Y\otimes\psi$,
$Y\in(\mu_v,\iota)$, $\psi\in\tau_v$, we have
$$\xi\cdot c=(T\otimes Y)\otimes\psi=(T\otimes 1_\iota\circ 1_\iota\otimes 
Y)\otimes\psi=$$
$$(T\circ Y)\otimes\psi=(1_\iota\otimes T\circ 
Y\otimes1_\iota)\otimes\psi=(Y\otimes T)\otimes\psi=c\cdot\xi.$$
\medskip

  \noindent{\bf 8.10.\ Proposition} {\sl
For any arrow $X\in(\mu_u,\mu_{u'})$,  the norm continuous extension of $\lambda(X)$ to
the completed Hilbert modules lies in the arrow space $(\text{Ind}(\mu_u), \text{Ind}(\mu_{u'}))$
of $\text{Bimod}_\alpha(G_\tau)$.}\medskip

\noindent{\it Proof}
Property $(7.1)$ shows that $\lambda$ is a $^*$--functor from the  $C^*$--category ${\cal M}_\mu$ to the category of prehilbertian $^\circ{\cal C}$--bimodules. Thus $\lambda(X)$ is bounded and
hence extends uniquely to a bimodule map between the completed Hilbert bimodules. On the other hand, the obvious commutation relations between $\lambda(X)$ and the action of $G_\tau$  
 imply that $\lambda(X)$ is an intertwining operator between 
 the corresponding bimodule representations of $G_\tau$. 
\medskip

  \noindent{\it 8.6. The induction functor $\text{Ind}:{\cal M}_\mu\to\text{Bimod}_\alpha(G_\tau)$.}
  \medskip
  
We may thus define a $^*$--functor of $C^*$--categories, 
 $$\text{Ind}: {\cal M}_\mu\to\text{Bimod}_\alpha(G_\tau)$$ 
 taking an object $\mu_u$ to $\text{Ind}(\mu_u)$ and an arrow $X\in(\mu_u,\mu_{u'})$
 to the extension of $\lambda(X)$. This is  the induction functor.
\medskip

\noindent{\it 8.7. The natural transformation $S$ and the generalized Swan's Theorem.} 
 The maps  $S_u$ defined in subsect.  8.1 extend uniquely to isometries 
$S_u:{\cal H}_u\to \tau_u\otimes{\cal C}$ making the following diagrams commute for $A\in(u,u')$,
\[
\xymatrix{
{\cal H}_u
   \ar[r]^-{\lambda\mu(A)}
   \ar[d]_-{S_u }
&  {\cal H}_{u'}
   \ar[d]^-{ S_{u'}}
\\ \tau_u\otimes{\cal C}
   \ar[r]^-{  \tau(A)\otimes I }
&   \tau_{u'}\otimes{\cal C}}
\]
\medskip

\noindent{\bf 8.11. Proposition} {\sl $S_u\in(\text{Ind}(\mu_u), 
\hat{u}\otimes\alpha)$ in the category 
$\text{Mod}_\alpha(G_\tau)$.}\medskip

\noindent{\it Proof} For $\xi= M\otimes \psi_i$, with $M\in(\mu_v,\mu_u)$, 
$(\psi_j)$ an orthonormal basis of $\tau_v$, and  orthonormal bases 
$(\phi_r)$ and $(\eta_p)$ of $\tau_u$ and $\tau_{\overline{u}}$ 
respectively,

$$\hat{u}\otimes{\alpha}\circ S_u(M\otimes\psi_i)=
\sum_r \hat{u}\otimes\alpha(\phi_r\otimes L(\phi_r)^*(M\otimes\psi_i))=$$
$$\sum_{r}\hat{u}\otimes\alpha(\phi_r\otimes(\mu(R_u^*)\circ
\tilde{\mu}_{\overline{u},
u}\circ1_{\mu_{\overline u}}\otimes
M\circ\tilde{\mu}^*_{\overline{u},v})\otimes(j_u\phi_r\otimes\psi_i))=$$
$$\sum_{r,s,p,j}\phi_s\otimes(\mu(R_u^*)\circ
\tilde{\mu}_{\overline{u},
u}\circ1_{\mu_{\overline u}}\otimes
M\circ\tilde{\mu}^*_{\overline{u},v})
\otimes(\eta_p\otimes\psi_j))\otimes\hat{u}_{sr}
({\hat{u}}_{{j_u}^*\eta_p,\phi_r})^*\hat{v}_{ji}=$$
$$\sum_{r,s,p,j,h}\phi_s\otimes((\mu(R_u^*)\circ
\tilde{\mu}_{\overline{u},
u}\circ1_{\mu_{\overline u}}\otimes
M\circ\tilde{\mu}^*_{\overline{u},v})
\otimes(\eta_p\otimes\psi_j))\otimes\hat{u}_{sr}
{\hat{u}}_{h,r}^*\hat{v}_{ji}<\eta_p, j_u\phi_h>=$$
$$\sum_{j,h}\phi_h\otimes(\mu(R_u^*)\circ
\tilde{\mu}_{\overline{u},
u}\circ1_{\mu_{\overline u}}\otimes
M\circ\tilde{\mu}^*_{\overline{u},v})
\otimes(j_u\phi_h\otimes\psi_j))\otimes\hat{v}_{ji}=$$
$$\sum_jS_u(M\otimes\psi_j)\otimes\hat{v}_{ji}=
S_u\otimes1_{{\cal Q}}\circ 
\text{Ind}(\mu_u)(M\otimes\psi_i).$$

\noindent
{\it Remark} $u\to S_u$ is a natural transformation from $\text{Ind}\mu$ to $\tau\otimes 1$   taking values in $\text{Mod}_\alpha(G_\tau)$.

\section{Extending $\text{Ind}$ to a full tensor functor}

As in the previous sections, we consider a pair of $^*$--functors between tensor $C^*$--categories,
$\tau:{\cal A}\to\text{Hilb}$ and $\mu:{\cal A}\to {\cal M}$, $\tau$ is tensorial and $\mu$ quasitensorial, and ${\cal A}$ has conjugates.
 
Now ${\cal M}_\mu$ is a $C^*$--category, but  not a tensor $C^*$--category in general.
This suggests looking for an extension of $\text{Ind}$ to ${\cal M}_\mu^\otimes$.
Here we show that the full tensor subcategory of $\text{Bimod}_\alpha(G_\tau)$ with objects $\text{Ind}(\mu_u)$ admits a natural realization as a strict tensor $C^*$--category
${\cal T}$
and that
 $\text{Ind}$ extends uniquely to a strict tensor isomorphism  
 between ${\cal M}_\mu^\otimes$ and ${\cal T}$.

In the next subsection we construct 
new  bimodules ${\cal H}_{\uu}$  associated with finite
sequences $\underline u=(u_1,\dots, u_n)$ of objects of ${\cal A}$.
This construction reduces to
that of the bimodules ${\cal H}_u$ 
of Sect. 7 and 8 for sequences of length $1$.

If $\uu=(u_1,\dots, u_n)$  is such a sequence and if
$\overline{u}_i$ is a conjugate of $u_i$, we
 write $\overline{\underline{u}}$ for
$(\overline{u}_n,\dots, \overline{u}_1)$.
If $R_{u_i}$, $\overline{R}_{u_i}$  is a 
solution of the conjugate equations for $u_i$,  we denote by 
$R_{\uu}$ and $\overline{R}_{\uu}$ the solutions of
 the conjugate equations for $u_1\otimes\dots\otimes u_n$
given by the product formula. Similarly, starting with the image 
solutions  $\hat R_{u_i}:=\tilde\mu_{\overline
u,u}^*\circ\mu(R_u)$ for $\mu_{u_i}$ in ${\cal M}$,
we use the product formula to define  the solution $\hat R_{\underline u}$ 
for 
$\mu_{u_1}\otimes\dots\otimes\mu_{u_n}$.

\medskip

\noindent
{\bf 9.1.\  Lemma} {\sl $\hat 
R_{\uu}=\tilde\mu_{\overline{\uu},\uu}^*\circ\mu(R_{\uu})$.}\medskip

\noindent
{\it Proof} We prove the lemma by induction on the length of $\uu$. The 
result holds
by definition if this length is one. Suppose $\uu=(\uv, w)$. Then, by construction,
$$\hat R_{\uu}=1_{\mu_{\overline w}}\otimes\hat 
R_{\uv}\otimes1_{\mu_w}\circ\hat R_w.$$
Hence by the induction hypothesis,
$$\hat R_{\uu}=1_{\mu_{\overline w}}\otimes
\tilde\mu_{\overline{\uv}, \uv}^*\otimes1_{\mu_w}\circ 
 1_{\mu_{\overline w}}\otimes\mu(R_{\uv})\otimes1_{\mu_w}\circ
\tilde\mu_{\overline w, w}^*\circ\mu(R_w).$$
But, by naturality,
$$1_{\mu_{\overline 
w}}\otimes\mu(R_{\uv})\otimes1_{\mu_w}\circ\tilde\mu_{\overline w, 
w}^*=
\tilde\mu_{\overline w, z , w}^*\circ
\mu ({1_{\overline w}\otimes R_{\uv}\otimes 1_w}),$$
where $z=\overline{v}_n\otimes\dots\otimes\overline{v}_1\otimes v_1\otimes\dots\otimes  v_n$
if $\uv=(v_1,\dots,v_n)$.
On the other hand,   the following relation follows easily from associativity of $\mu$, 
$$1_{\mu_{\overline 
w}}\otimes\tilde\mu_{\overline{\uv}, \uv}^*\otimes1_{\mu_w}\circ 
\tilde\mu_{\overline w, z,  w}^*=
\tilde\mu_{{\overline\uu}, \uu}^*,$$
thus completing the proof.\medskip

\noindent{\it 9.1. New Hilbert bimodules ${\cal H}_{\uu}$.}

Let $\underline{u}=(u_1,\dots,u_n)$ be a finite sequence  of objects of 
${\cal A}$
and set
$$\mu_{\uu}:=\mu_{u_1}\otimes\dots\otimes\mu_{u_n}.$$ 
Let 
${}^\circ{\cal H}_{\underline u}$
be the linear space
$\sum_v(\mu_v,\mu_{\underline u})\otimes\tau_v$,
the sum being taken over the objects of ${\cal A}$, quotiented
by the linear subspace generated by elements of the form
$M\circ\mu(A)\otimes T-M\otimes\tau(A)\circ T.$

We proceed as in the construction of the bimodules ${}^\circ{\cal H}_u$ of Sect.\ 7, 
defining successively multiplication, the functor $\lambda$, adjoint and sesquilinear form. 

Define bilinear maps
${}^\circ{\cal H}_{\underline u}\times
{}^\circ{\cal H}_{\underline u'}\to
{}^\circ{\cal H}_{\underline 
u, \underline u'}.$
For $\xi=L\otimes \psi$, $\eta=M\otimes \phi$, where $L\in(\mu_w, 
\mu_{\underline u})$,
$M\in(\mu_v, \mu_{\underline u'})$,
$\psi\in \tau_w$, $\phi\in \tau_v$,
set:  
$$\xi\eta:=(L\otimes
M)\circ\tilde{\mu}_{w,v}^*\otimes(\psi\otimes \phi).$$
It is easy to check that these maps are well defined and 
associative.

 For a reason that will soon become clear, this new multiplication 
does {\it not} coincide with the multiplication $\xi\cdot{}\eta$ of Sect. 7
  if $\uu$ or $\uu'$ are objects of ${\cal A}$. We have therefore
used a different notation. 
However, the two multiplications coincide if $\uu$ or $\uu'$ are the tensor 
unit
$\iota$, as $\tilde{\mu}_{\iota,\uu}=\tilde{\mu}_{\uu,\iota}=1_{\mu_{\uu}}$.
Hence ${}^\circ{\cal H}_{(\iota)}$ is again 
the  algebra
${}^\circ{\cal C}$ and
${}^\circ{\cal H}_{\underline u}$ is  a ${}^\circ{\cal C}$--bimodule. Furthermore, as a bimodule, we do have ${}^\circ{\cal H}_{\uu}={}^\circ{\cal H}_u$  if $\uu=(u)$.

Given
$Y\in(\mu_{\underline u},\mu_{\underline u'})$
 a bimodule map,
$\lambda(Y): {}^\circ{\cal H}_{\underline u}\to
{}^\circ{\cal H}_{\underline u'}$ is defined in the obvious way. $\lambda$ is 
a  covariant   functor from
${\cal M}_\mu^\otimes$  to the category of $^\circ{\cal 
C}$--bimodules.
In particular, if $\uu=(u_1,\dots,u_n)$,
$$\lambda(\tilde{\mu}_{\underline u}):
 {}^\circ{\cal H}_{\underline u}\to
{}^\circ{\cal H}_{u_1\otimes\dots\otimes u_n}$$
relates the new and old bimodules.
Moreover, the respective multiplications   are related by
$$\lambda(\tilde{\mu}_{\underline u, \underline u'})(\xi\eta)=(\lambda(\tilde{\mu}_{\underline u})\xi)\cdot(\lambda(\tilde{\mu}_{{\underline u}'})\eta),\eqno(9.1)$$
for $\xi\in{}^\circ{\cal H}_{\uu}$, $\eta\in{}^\circ{\cal H}_{\underline{u}'}$.
In particular
$$\lambda(\tilde{\mu}_{u,v})\xi\eta=\xi\cdot{}\eta,\quad \xi\in{}^\circ{\cal H}_u,\quad \eta\in{}^\circ{\cal H}_v.$$
The adjoint
$^*:{}^\circ{\cal H}_{\underline 
u}\to{}^\circ{\cal
H}_{\overline{\underline u}}$ is defined as before by
$$(M\otimes \psi)^*:=M^\bullet\otimes j_v\psi,$$
but where ${}^\bullet$ now refers to $\hat R_{\uu}$.
The adjoint is 
well defined and independent of the choice of solutions of the 
conjugate
equations for  $v$. However,
if
we change the solution of the conjugate equations for $\underline u$  using a
sequence of invertibles
$(X_1,\dots,X_n)$ wth $X_i\in(\overline{u_i},\tilde{ u}_i)$, 
$(M\otimes \psi)^*$ becomes
$(\mu(X_1)\otimes\dots\otimes\mu(X_n)\circ M^\bullet)\otimes
j_v\psi$.

Nevertheless  $^*$
is an antilinear
map    
satisfying the properties of Prop.\ 7.1  if we use tensor product and conjugate solutions.

The $^\circ{\cal C}$--valued form on
${}^\circ{\cal H}_{\underline u}$ is defined by
$$<\xi,\xi'>:=\lambda(\hat{R}_{\underline u}^*
)(\xi^*\xi').$$
One can  easily check that a formula similar to $(7.1)$ holds,
$$<\xi,\xi'>=(\hat R_v^*\circ
1_{\mu_{\overline v}}\otimes (M^*\circ M')\circ\tilde{\mu}_{\overline{v},
v'}^*)\otimes
 j_v\psi\otimes
\psi',\eqno(9.2)$$
hence the form  reduces to that of  ${}^\circ{\cal H}_u$ if $\uu=(u)$.
As before, this form is independent of the choice of the conjugate of 
${\underline u}$ in view of how $\xi^*$ changes and of Lemma 9.1.
The above expression shows that  $\lambda$ is 
a $^*$-functor.

Since $\tilde{\mu}_{\underline u}$ is an isometry, 
$(9.2)$  shows that 
$\lambda(\tilde{\mu}_{\underline u})$ preserves the corresponding
forms:
$$<\lambda(\tilde{\mu}_{\underline u})\xi, \lambda(\tilde{\mu}_{\underline u})\xi'>=
<\xi,\xi'>,\quad \xi,\xi'\in{}^\circ{\cal H}_{\underline u},$$
generalizing Prop\ 8.7.
On the other hand,   ${}^\circ{\cal H}_{u_1\otimes\dots\otimes u_n}$ is
 a finite projective prehilbertian bimodule, so the same is true of
${}^\circ{\cal H}_{\underline u}$.

Completing ${}^\circ{\cal H}_{\underline u}$
in the norm derived from the maximal 
$C^*$--norm of ${}^\circ{\cal C}$ yields a Hilbert ${\cal C}$--bimodule 
${\cal H}_{\underline u}$. $\lambda$ extends to a $^*$--functor from
 ${\cal M}_{\mu}^\otimes$ to the $C^*$--category of Hilbert ${\cal C}$--bimodules. In this category
${\cal H}_{\underline u}$ is a subobject of ${\cal H}_{u_1\otimes\dots\otimes u_n}$.

\medskip

We next regard the associative multiplication $\xi,\eta\to\xi\eta$  as a bimodule map defined on 
the bimodule  tensor product 
$
{}^\circ{\cal H}_{\underline u}
\otimes_{^\circ{\cal C}}{}^\circ{\cal H}_{{\underline u}'}\to{}^\circ{\cal H}_{\underline u, {\underline u}'}.$
\medskip

\noindent{\bf 9.2. Theorem} {\sl  The multiplication $\xi\otimes\eta\to\xi\eta$
extends uniquely to a unitary map
between Hilbert  bimodules
$$\xi\otimes\eta\in {\cal H}_{\underline u}\otimes_{{\cal C}}{\cal H}_{{\underline u}'}\to\xi\eta\in{\cal H}_{\underline u, {\underline u}'}.$$ 
}\medskip

\noindent{\it Proof}
Using successively that  $\lambda(\tilde{\mu}_{\underline u})$ is isometric,  Prop.\ 8.7, and
relation $(9.2)$, 
  we conclude that
multiplication is a densely defined isometry, extending to an isometry of the completions.
To complete the proof it suffices to show that  the set of all $\xi\eta$, with 
$\xi\in{}^\circ{\cal H}_{\uu}$, $\eta\in{}^\circ{\cal H}_{\uu'}$ span ${}^\circ{\cal H}_{\uu,\uu'}$.

 Since multiplication is associative, it suffices to choose $\underline u$ to be a sequence $(u)$ consisting of a single element.
Consider an element of ${}^\circ{\cal H}_{u,{\underline u}'}$ of the form
$M\otimes\psi$, where $M\in(\mu_v, \mu_u\otimes\mu_{{\underline u}'})$ and $\psi\in\tau_v$.  Using  the explicit linear isomorphism $(\mu_v, \mu_u\otimes\mu_{{\underline u}'})\simeq (\mu_{\overline{u}}\otimes\mu_v, \mu_{{\underline u}'})$,
we may write  $M=1_{\mu_u}\otimes T\circ\hat{\overline{R}}_u\otimes 1_{\mu_v}$ where $T\in(\mu_{\overline{u}}\otimes\mu_v, \mu_{{\underline u}'})$.
We may also write $T=T'\circ\tilde{\mu}_{\overline{u},v}$, with 
$T'\in(\mu_{\overline{u}\otimes v}, \mu_{{\underline u}'})$. Hence 
$$M=1_{\mu_u}\otimes (T'\circ\tilde{\mu}_{\overline{u},v})\circ\hat{\overline{R}}_u\otimes 1_{\mu_v}=1_{\mu_u}\otimes T'\circ 1_{\mu_u}\otimes \tilde{\mu}_{\overline{u},v}\circ\tilde{\mu}^*_{u,\overline{u}}\otimes 1_{\mu_v}\circ \mu({\overline R}_u)\otimes 1_{\mu_v}=$$
$$1_{\mu_u}\otimes T'\circ
\tilde{\mu}^*_{u, \overline{u}\otimes v}\circ\tilde{\mu}_{u\otimes\overline{u}, v}
  \circ \mu({\overline R}_u)\otimes 1_{\mu_v}= 1_{\mu_u}\otimes T'\circ \tilde{\mu}^*_{u, \overline{u}\otimes v}\circ\mu(\overline{R}_u\otimes 1_v).$$
Substituting this into $M\otimes\psi$ gives
$$M\otimes\psi=(1_{\mu_u}\otimes T'\circ \tilde{\mu}^*_{u, \overline{u}\otimes v})\otimes
(\tau(\overline{R}_u)\otimes \psi).$$
Writing $\tau(\overline{R}_u)=\sum_j \phi_j\otimes j_u\phi_j$, for an orthonormal basis $(\phi_j)$ of $\tau_u$, gives 
$$M\otimes\psi=\sum_j \xi_j\eta_j,\quad \xi_j=1_{\mu_u}\otimes \phi_j\in{}^\circ{\cal H}_u\quad \eta_j=T'\otimes(j_u\phi_j\otimes\psi)\in{}^\circ{\cal H}_{{\underline u}'}.$$ 
\medskip

On one hand, as the multiplication maps are isometric
${\cal H}_{\uu}{\cal H}_{\uu'}$ realizes  the tensor product of Hilbert 
bimodules. It has the virtue 
of being {\it strictly associative}, as so are the multiplication 
maps.
We replace the original tensor product of Hilbert bimodules by this
strictly associative tensor product. 
On the other hand, since the multiplication maps are unitary,
we have tensor product decompositions, 
$${\cal H}_{\underline u}={\cal H}_{u_1} \cdots{\cal H}_{u_n},$$ 
for $\uu=(u_1,\dots,u_n).$
(Note that the right hand side is already norm closed, by finite projectivity.)
In particular, if $1_{\mu_{u_i}}\neq0$ for all $i$, ${\cal H}_{\underline u}$
is a full right Hilbert module with a faithful right ${\cal C}$--action. We thus have the following result.
\medskip

\noindent{\bf 9.3.\  Theorem} {\sl  $\lambda$ is a strict tensor 
$^*$--functor from ${\cal M}_\mu^\otimes$ to the tensor
$C^*$--category of Hilbert ${\cal C}$--bimodules (with a strictly associative tensor product).}\medskip

\noindent{\it 9.2. $G_\tau$--representations on the bimodules ${\cal H}_{\underline u}$.} \medskip

\noindent{\bf 9.4.\ Proposition} {\sl Given a finite sequence $\underline 
u=(u_1,\dots,u_n)$ of objects of ${\cal A}$, there is a unique bimodule 
representation $\text{Ind}(\mu_{\underline u})$ 
of $G_\tau$ on ${\cal H}_{\underline u}$ such that, for $M\in(\mu_v,\mu_{\uu})$,
$\psi\in\tau_v$,
$$\text{Ind}(\mu_{\underline u})(M\otimes\psi)=M\otimes\hat{v}\psi.$$
$\text{Ind}(\mu_{\uu})$ is a full bimodule representation.
Under the tensor product decomposition
${\cal H}_{\underline u}={\cal H}_{u_1}\dots{\cal H}_{u_n}$,
$\text{Ind}({\mu_{\uu}})=\text{Ind}(\mu_{u_1})\otimes\dots\otimes\text{Ind}(\mu_{u_n}).$  }\medskip

\noindent{\it Proof} The only non--trivial statement is 
that 
$\text{Ind}(\mu_{\underline u})$ is full. A   
$G_\tau$--fixed vector ${\cal H}_{\underline u}$ is of the form
$M\otimes 1_\iota$, with $M\in(\iota,\mu_{\underline u})$. As for the old 
modules, Prop.\ 8.9,
one shows this element to be central.\medskip

As before, $\lambda(X)$ 
 intertwines
the representations $\text{Ind}(\mu_{\underline u})$ and $\text{Ind}(\mu_{{\underline u}'}))$  for $X\in(\mu_{\underline u}, \mu_{{\underline u}'})$. We therefore have a strict tensor functor 
$$\text{Ind}:{\cal M}_{\mu}^\otimes\to\text{Bimod}_\alpha(G_\tau),$$
the unique tensor extension of the functor $\text{Ind}$ defined on ${\cal M}_\mu$ of the previous section.

 We are now ready to state 
 a central result of this paper, a version of the Frobenius reciprocity theorem 4.1 for quasitensor functors.
 \medskip
 
 \noindent{\bf 9.5.\  Theorem} {\sl $\text{Ind}$ is a full and faithful strict tensor functor from ${\cal M}_\mu^\otimes$ to the category of bimodule representations of $G_\tau$,
 if  the latter is endowed with a strictly associative tensor product.}\medskip
 
 \noindent{\it Proof} ${\cal M}$ is a tensor $C^*$--category with conjugates
 and $\text{Ind}$ a relaxed tensor functor, hence automatically faithful 
 \cite{Ergodic}. It remains to show that $\text{Ind}$ is full.
 This follows from the linear isomorphisms 
  $\gamma: (\mu_{\underline u},\mu_{{\underline u}'})\to(\iota,
\mu_{{\underline u}',\overline{\underline u}})$,
 $T\to X=T\otimes 1_{\mu_{\overline{\underline u}}}\circ\hat{\overline{R}}_{\underline u},$
and 
$\delta: (\text{Ind}(\mu_{\underline u}), \text{Ind}(\mu_{{\underline 
u}'}))
\to(\iota,
\text{Ind}(\mu_{{\underline u}',\overline{\underline u}}))$, defined similarly, where 
$\mu$ is replaced by the quasitensor functor $\text{Ind}\mu$. 
Hence any intertwiner in 
$(\text{Ind}(\mu_{\underline u}), \text{Ind}(\mu_{{\underline u}'}))$ is 
determined by a fixed vector in ${\cal H}_{{\underline u}',\overline{\underline u}}$, which we already know 
to arise from an intertwiner in $(\iota, \mu_{{\underline u}',\overline{\underline u}})$, hence  lying in the image of $\text{Ind}$.\medskip

\noindent{\it Remark} The last proof  uses only the  functor of tensoring
on  the ${\it right}$ by an identity arrow. This also makes sense for ${\it module}$ intertwiners and hence shows the following result.
\medskip

\noindent{\bf 9.6. Theorem} {\sl Any  module intertwiner from $\text{Ind}(\mu_{{\underline u}})$ to
 $\text{Ind}(\mu_{{\underline u}'})$, namely an intertwiner in the 
$C^*$--category $\text{Mod}_\alpha(G_\tau)$, is automatically a bimodule intertwiner.}

\medskip

\noindent{\it 9.3. The functor $\text{Ind}\mu:{\cal A}\to\text{\rm Bimod}_\alpha(G_\tau)$.}
We finally define the composed  functor 
$$\text{Ind}\mu:{\cal A}\to\text{ Bimod}_\alpha(G_\tau).$$
\medskip

\noindent{\bf 9.7. Theorem} {\sl If  $(\mu,\tilde{\mu})$ is a quasitensor (relaxed tensor, tensor) functor,
$\text{Ind}\mu$ is a quasitensor (relaxed tensor, tensor) functor too,
with natural transformation $\widetilde{\text{Ind}(\mu})$ given by
the $\cdot{}$--multiplication maps,
$$\text{Ind}(\tilde{\mu}_{u,u'}):\xi\xi'\in{\cal H }_u{\cal H}_{u'}\to\xi\cdot\xi'\in{\cal H}_{u\otimes u'}.$$
Furhermore, $\text{Ind}$ is tensor isomorphism from $(\mu,\tilde\mu)$ to $(\text{Ind}\mu, \widetilde{\text{Ind}\mu})$.
}\medskip

\noindent{\it Proof} 
Since $\text{Ind}$ is a strict tensor functor  and $\mu$ is quasitensor (relaxed tensor, tensor, resp.), 
their composition $\text{Ind}\mu$,  with the composed natural transformation, is quasitensor
(relaxed tensor, tensor, resp.), see subsect.\ 2.2. This natural transformation is precisely the map $\xi\xi'\to\xi\cdot{}\xi'$.
The last statement is clear.
\medskip


We next prove Theorems 6.5 and 6.6.
\medskip

{\it Proof of Theorems  6.5 and 6.6.}
 We briefly recall from  
\cite{Ergodic} 
how to get a pair of functors $\mu$ and $\tau$.
Assume that ${\cal M}$ has conjugates and an irreducible tensor unit $\iota$, 
and fix an object $x$ in ${\cal M}$ with intrinsic dimension $>1$ and a 
standard solution $R$, 
$\overline{R}$ of the conjugate equations for $x$. By Jones's result
\cite{Jones1} the intrinsic dimension  of $x$ can only take the values
$d=2\cos\frac{\pi}{\ell}$, for $\ell=3,4,\dots$ or $d\geq2$.
Consider the universal tensor $^*$--category ${\cal T}_{d}$ 
with objects the finite words in $u$ and $\overline{u}$ and whose arrows
are generated by two arrows $S\in(\iota, \overline{u}\otimes u)$ and 
$\overline{S}\in(\iota, u\otimes\overline{u})$ subject to the relations
expressing $(S, \overline{S})$ as a normalized solution of the conjugate 
equations 
for $u$. $\iota$ is the empty word 
and acts as a 
tensor unit. ${\cal T}_d$ is a tensor 
$C^*$--category for the allowed values of $d$. Furthermore there is 
a tensor functor $\tau$ from ${\cal T}_d$ to the category of Hilbert spaces
if and only if $d\geq2$, and all such functors can be easily classified.
Picking an embedding $\tau$, we get an associated compact quantum group 
$G_\tau=A_u(F)$, where $F$ is an invertible matrix such that 
$\text{Tr}(FF^*)=\text{Tr}((FF^*)^{-1})=R^*R$.  Furthermore we 
have a canonical tensor functor $\mu:{\cal T}_d\to{\cal M}$
such that $\mu(u)=x$, $\mu(\overline{u})=\overline{x}$, $\mu(S)=R$,
$\mu(\overline{S})=\overline{R}$, and we may now apply our main result.

Similarly, given a real or pseudoreal solution of the conjugate equations
in ${\cal M}$, namely $R\in(\iota, x^2)$ with $R^*\otimes 1_x\circ 
1_x\otimes R=\pm 1_x$,  we consider the associated 
universal Temperley--Lieb
categories 
${\cal T}_{rd}$ and
${\cal T}_{pd}$ with generating arrow $S\in(\iota, u^2)$. If $R^*\circ 
R\geq2$, a choice of an embedding of ${\cal T}_{pd}$ or ${\cal T}_{rd}$ 
into the Hilbert 
spaces provides a quantum group $A_o(F)$ with $F$ an invertible matrix
satisfying $F\overline{F}=\pm I$, 
$\text{Tr}(FF^*)=\text{Tr}((FF^*)^{-1})=R^*\circ R$.
\medskip

\section{An adjoint pair of functors}

Recall that a pair of functors $F:{\Phi}\to{\Phi}'$ and 
$F':{\Phi}'\to{\Phi}$ between  categories is an adjoint pair 
if, for any 
pair of objects $\phi\in{\Phi}$, $\phi'\in{\Phi}'$, there is a  
isomorphism
$\beta_{\phi',\phi}: (\phi', F_\phi)\to(F'_{\phi'}, \phi)$ natural in 
$\phi$ and $\phi'$.

In this section we show that, essentially by construction, the pair $(\text{Ind}, \mu)$ gives rise to an adjoint pair.

To this end, we assume as before  that $\tau:{\cal A}\to\text{Hilb}$ is a tensor 
functor 
into the category of Hilbert spaces, 
so that ${\cal A}$ is a
category of representations of a compact quantum group $G_\tau$, and that $\mu:{\cal A}\to{\cal M}$ is a 
quasitensor functor of  strict tensor $C^*$--categories with 
irreducible tensor units and construct the corresponding Hilbert 
$C^*$--bimodules.

Following Mackey's construction of the induced representation for locally 
compact groups,
we  consider the scalar--valued inner product on
${\cal H}_{\underline{u}}$ given by composing the ${\cal C}$--valued 
inner product with the unique $G_\tau$--invariant faithful state. We thus get
a Hilbert space, $H_{\underline{u}}$ and the bimodule representation 
$\text{Ind}(\mu_{\underline{u}})$ of $G_\tau$ defines a densely defined 
representation of $G_\tau$ on $H_{\underline{u}}$. 
This representation is isometric as the state is invariant and  
extends uniquely to a unitary 
representation of $G_\tau$ 
again denoted by $\text{Ind}(\mu_{\underline{u}})$. However, although 
we start with a finite dimensional representation, the Hilbert space of the induced 
representation in general fails to be finite dimensional, hence we need to
work with the category of not necessarily finite dimensional unitary
representations of $G_\tau$, denoted by $\widetilde{\text{Rep}}(G_\tau)$. We thus have a  functor, 
$\text{Ind}:{\cal 
M}_\mu^\otimes \to
\widetilde{\text{Rep}}(G)$. We let $\widetilde{{\cal M}_\mu^\otimes}$ be the tensor $W^*$--category 
completion of ${\cal M}_\mu^\otimes$ under infinite direct sums, cf. \cite{LR}. 
Then $\mu$ and $\text{Ind}$ extend uniquely to  $^*$--functors on
$\widetilde{\text{Rep}}(G)$ and $\widetilde{{\cal M}_\mu^\otimes}$ 
respectively.

\medskip

\noindent{\bf 10.1.\ Theorem} {\sl The pair of functors
$\text{Ind}:\widetilde{{\cal 
M}_\mu^\otimes}\to\widetilde{\text{Rep}}(G_\tau)$ and
$\mu:\widetilde{\text{Rep}}(G_\tau)\to\widetilde{{\cal M}_\mu^\otimes}$ 
is an adjoint 
pair.}\medskip

\noindent{\it Proof}
Note that the linear span 
of the images of elements of the form
$T\otimes\psi$, where $T\in(\mu_v,\mu_{\underline{u}})$ and 
$\psi\in\tau_v$, where $v$ runs over the irreducible representations of $G_\tau$,  
 is dense in $H_{\underline{u}}$.  
If we fix an irreducible $v$, the space of intertwiners
$(v,\text{Ind}(\mu_{\underline{u}}))$ is given precisely by 
the set of maps $\hat{T}: \psi\in\tau_v\to T\otimes\psi\in 
H_{\underline{u}}$,
with $T\in
(\mu_{v},\mu_{\underline{u}})$. Hence there is
a linear isomorphism $(\mu_{v},\mu_{\underline{u}})
\to
(v,\text{Ind}(\mu_{\underline{u}}))$, natural in 
$\mu_{\underline{u}}$.
This isomorphism extends uniquely 
to a linear isomorphism natural in $v$.
\medskip

\section {Full bimodule representations from group 
actions}

Which ergodic actions $({\cal C},\alpha, G)$ can arise from a pair of  functors $\text{Hilb}\stackrel\tau\leftarrow{\cal A}\stackrel\mu\rightarrow{\cal M}$ with $\mu$ relaxed tensor? 
In this section we attack this  problem   
when $G_\tau$ is a   compact {\it group}. 

By Theorems 6.1 and 6.2, a necessary condition is that $u\otimes\alpha$
can be made into a full bimodule representation of $G$ for every $u\in\text{Rep}(G)$.
Using the results of \cite{Takesaki, HLS}, given an ergodic action $({\cal C},\alpha)$ of $G$ we
shall construct, a canonical
full  bimodule representation on every $G$--submodule $X_u\subset H_u\otimes{\cal C}$, where $u$ is an object of $\text{Rep}(G)$.
The submodule is the full module $H_u\otimes{\cal C}$ for all $u$ precisely when the multiplicity of $u$ is maximal.

It turns out that $u\to X_u$ is   a  quasitensor functor  in general , 
related to the spectral functor of the ergodic action, and is relaxed tensor when the ergodic action is of full multiplicity.

We
shall classify the full bimodule structures on the 
intermediate projective $G$--submodules    
$X_u\subset Y\subset u\otimes\alpha$ assuming that the 
weak completion of ${\cal C}$ in the GNS representation 
of the invariant trace is 
  a finite type $I$ von 
Neumann 
algebra, and hence induced by an ergodic action $\beta$ of a closed 
subgroup $K$ on a matrix algebra ${\cal F}$, cf.\ subsect.\ 2.3 (Theorem 
11.10).  

This provides an obstruction to the existence of full bimodule representation 
structures when multiplicities of the primitive action of $K$ are low but nonzero.
More precisely, we shall show that   
certain low multiplicity actions as well as certain ergodic actions of $SU(2)$ are excluded
(Cor.\ 11.12, Example 11.11).
\medskip

\noindent{\it 11.1. The minimal full bimodule representations.} 
Let $({\cal F},\beta, K)$ be an ergodic action of a compact group $K$ on 
a unital $C^*$--algebra ${\cal F}$.
 Given 
 an object $v$ of $\text{Rep}(K)$, 
recall from Sect.\ 2 that   the  spectral
space $\overline{L}_v$ is the complex conjugate
of the set $L_v$ of all linear maps $T:H_v\to{\cal F}$ intertwining $v$ 
with 
$\beta$ and that, by the  multiplicity bound theorem, subsect.\ 2.3, 
$\text{dim}(\overline{L}_v)\leq\text{dim}(v)$.
If $v$ is irreducible, $\text{dim}(\overline{L}_v)$ is   the
multiplicity of $v$ in $\beta$, denoted by $\text{mult}(v)$.
 Recall also that $\overline{L}_v$ is a 
 Hilbert space with inner product 
$<\overline{S},\overline{T}>:=\sum_i S(\psi_i)T(\psi_i)^*$, where 
$(\psi_i)$ is an orthonormal basis of $H_v$.

If $\overline{L}_v\neq0$, we  construct
 a natural nonzero $K$--module subrepresentation   of
 $v\otimes\beta$ 
with a full bimodule structure.

Consider the linear map $Z_v: \overline{T}\in\overline{L}_v\to\sum\psi_i\otimes T(\psi_i)^*\in H_v\otimes{\cal F}$, clearly independent of the choice of the orthonormal basis. The range of $Z_v$ is the space of $K$--fixed vectors in $H_v\otimes{\cal F}$ for the action 
$v\otimes\beta$. 
We may identify $Z_v$ with a rectangular matrix
$(T_k(\psi_i)^*)$, where $\overline{T}_k$ is an orthonormal basis  of $\overline{L}_v$.
Since  $\iota\otimes\beta_k(Z_v)=v(k)^*\otimes IZ_v$,   $Z_v^*$
is referred to as 
  the  {\it eigenmatrix}
of $\beta$ in \cite{Wassermann1}.

Now, $Z_v$   takes the  inner product of $\overline{L}_v$ to the inner product of the range inherited 
from $H_v\otimes{\cal F}$ as a right  Hilbert ${\cal F}$--module, 
which on 
that subspace indeed  takes values in ${\mathbb C}$ 
since $\beta$ is ergodic.

Extend  $Z_v$ uniquely to an adjointable bounded map between Hilbert 
modules $\overline{L}_v\otimes{\cal F}
\to  H_v\otimes{\cal F}$, still denoted $Z_v$. 
It is easy to verify that $Z_v$ is an isometry, $Z_v^*Z_v=I$, intertwining the module representations $\iota_{\overline{L}_v}\otimes\beta$ and $v\otimes\beta$, where $\iota_{\overline{L}_v}$ is the trivial representation of $K$ on the Hilbert space $\overline{L}_v$.
We may clearly identify  
${\cal L}_{\cal F}({\overline{L}_v}\otimes{\cal F}, H_v\otimes{\cal F})\simeq{\cal L}({\overline{L}_v}, H_v)\otimes{\cal F}$. 
Note that $X_v:=Z_v(\overline{L}_v\otimes{\cal F})$
is a projective  ${\cal F}$--submodule of $H_v\otimes{\cal F}.$  
\medskip

\noindent{\it Remark} An easy computation shows that, if $\tau$ is the unique invariant normalized tracial state
of ${\cal F}$ and $\text{Tr}$ is the non-normalized trace of ${\cal L}(H_v)$,
$\text{Tr}\otimes\tau(E_v)=\text{dim}(\overline{L}_v)$, where $E_v:=Z_vZ_v^*$. Hence $Z_v$ is a unitary
if and only if $\text{dim}(L_v)=\text{dim}(H_v)$, i.e.\ $\overline{L}_v$ must have maximal dimension.
\medskip

By the intertwining property of $Z_v$, 
$X_v$ is
$K$--invariant.
There is  a faithful unital $^*$--homomorphism $$\zeta:{\cal F}\to 
E_v{\cal 
L}(H_v)\otimes{\cal F}E_v,$$  
$$\zeta(f)=Z_vI\otimes fZ_v^*,$$  making
$X_v$ into  a ${\cal F}$--bimodule.

\medskip

\noindent{\bf 11.1.\ Proposition} {\sl For any  representation $v$ of    $K$ 
with $\overline{L}_v\neq0$, 
 $\zeta$ makes $X_v$ into a nonzero full bimodule 
$K$--representation, 
  isomorphic to 
$\overline{L}_v\otimes{\cal F}$ via $Z_v$ with trivial left and right ${\cal 
F}$--actions, where
 $K$ acts 
as $\iota_{\overline{L}_v}\otimes\beta$.}\medskip

\noindent{\it Proof} 
By construction, $K$-fixed vectors are 
${\cal F}$--central in $X_v$, as they correspond via $Z_v$ to the  fixed vectors for $\iota_{\overline{L}_v}\otimes\beta$, namely to $\overline{L}_v\otimes{\mathbb C}$, clearly central for the trivial bimodule action.  
Property $(3.5)$ follows just as easily.
\medskip

\noindent{\it 11.2. The intermediate full bimdule representations.}
Given $({\cal F}, K,\beta)$ as before, 
we look for extensions of $X_v$ to full bimodule structures on 
intermediate projective $K$--module subrepresentations 
$$X_v\subset Y\subset H_v\otimes{\cal F}.$$
Clearly, such submodules are the ranges of projections
$E\in{\cal L}(H_v)\otimes{\cal F}$ satisfying
$$E\geq E_v,\eqno(11.1)$$
and the $K$--invariance condition
$$\text{Ad}v(k)\otimes\beta_k(E)=E,\quad k\in K.\eqno(11.2)$$
In what follows, we  set $Z_v=0$ and $X_{v}=\{0\}$ if 
$\overline{L}_v=\{0\}$.

\medskip

\noindent{\bf 11.2.\ Proposition} {\sl
Given $v\in\text{Rep}(K)$ and a  projection
$E\in{\cal L}(H_v)\otimes{\cal F}$ satisfying $(11.1)$ and $(11.2)$,   a 
unital 
$^*$--homomorphism $\eta:{\cal F}\to E{\cal L}(H_v)\otimes{\cal F}E$ 
defines
a full $K$--bimodule representation   on $Y=E(H_v\otimes{\cal F})$
   if and only if
$$\eta(\beta_k(f))=\text{Ad}v(k)\otimes\beta_k(\eta(f)),\quad k\in 
K,\eqno(11.3)$$
$$\eta(f)Z_v=Z_vI\otimes f,\quad f\in{\cal F}.\eqno(11.4)$$
}\medskip

\noindent{\it Proof} The proof is straightforward. We just note that 
$(11.3)$ corresponds to left $K$--equivariance in the sense of $(3.5)$,
whilst the property of being a full representation
is expressed by  $(11.4)$, as for $X_v$.
\medskip

\noindent{\bf 11.3. Corollary} {\sl If $\overline{L}_v$ has maximal dimension, 
  $v\otimes\beta$ 
becomes a full 
bimodule $K$--representation in a unique way.}\medskip

\noindent{\it Proof}
 By the previous 
 remark, 
 $Z_v$ is a unitary in   
${\cal L}(H_v)\otimes {\cal F}$, and $X_v=H_v\otimes{\cal F}$.
$\eta$ is uniquely determined by $(11.4)$. This formula defines
a $^*$--homomorphism clearly satisfying $(11.1)$--$(11.3)$ for $E=I$.
\medskip

\noindent{\it 11.3. Intermediate full bimodules for induced $C^*$--actions}
Now assume that $K$ is a closed subgroup of a compact group $G$ 
acting on ${\cal F}$,
which may be either a $C^*$--algebra or a von Neumann algebra. This action, $\beta$,
is supposed to be continuous in the appropriate topology. 
Consider the induced  algebra $\text{Ind}({\cal F})$ defined as in subsect.\ 2.3,
the action $\rho$ of $G$ being given by right translation.

Let  $v$ 
be a f.d.\ unitary representation of $G$. In the next,  known, proposition we  
determine the spectral spaces for 
the action $\rho$ in terms of those of the original 
action $\beta$.\medskip

\noindent{\bf 11.4. Proposition} {\sl The map  $T\in 
L^\beta_{v\upharpoonright_K}\to T' \in L^\rho_v$, with 
$T':H_v\to\text{Ind}({\cal F})$ 
defined by 
$T'(\psi)(g):=T(v(g)\psi)$, is unitary. As a consequence, 
$$Z^\rho_v(g)=v(g)^*\otimes I Z^\beta_{v\upharpoonright_K},$$
 hence
$$E^\rho_v(g)=v(g)^*\otimes I E^\beta_{v\upharpoonright_K} v(g)\otimes I.$$}\medskip

\noindent{\it Proof} 
Let us extend $\beta$ and $\rho$ 
to unitary representations of $K$ and $G$, respectively, on the 
$L^2$-completions of ${\cal F}$ and $\text{Ind}({\cal F})$
for the unique invariant traces. The extension of $\rho$ is clearly the
 representation induced from the extension of $\beta$ in the sense of  Mackey.
Extending in this way does not increment the spectra. Hence 
$L^\rho_v$
may be determined
by the classical Frobenius reciprocity theorem, showing that $T\to T'$ is a linear isomorphism. It is easily checked to be an isometry.
Therefore
$$Z^\rho_v(\overline{T'})(g)=\sum_i\psi_i\otimes T(v(g)\psi_i)^*,$$
showing that, if  $(\psi_i)$ is an orthonormal basis of $H_v$
and $(\overline{T}_j)$ an orthonormal basis of $\overline{L}_v$, 
then $(\xi_j')$, where
$\xi_j'(g):=\sum \psi_k\otimes T_j(v(g)\psi_k)^*$, is an orthonormal basis of the Hilbert space 
of $G$--fixed vectors in $H_v\otimes\text{Ind}({\cal F})$, hence the 
$jr$-- entry of $Z^\rho_v$ is the function
$$T_r(v(g)\psi_j)^* =(v(g)^*\otimes IZ^\beta_{v\upharpoonright_K})_{jr}.$$
\medskip

If $z$ and $v$ are representations of $K$, we identify the space of bounded adjointable
${\cal F}$--module maps ${\cal L}_{\cal F}(H_z\otimes{\cal F}, 
H_v\otimes{\cal F})$ 
with ${\cal L}(H_z, H_v)\otimes{\cal F}$. Hence $$(z\otimes\beta, v\otimes\beta)=\{T\in{\cal 
L}(H_z, H_v)\otimes{\cal F}: 
\iota\otimes\beta_k(T)=v(k)^*\otimes IT
z(k)\otimes I\}.$$
As a $K$-space, this space is linearly isomorphic to 
$$(H_{v\otimes\overline{z}}\otimes{\cal 
F})^{v\otimes\overline{z}\otimes\beta}=\overline{L}^\beta_{v\otimes\overline{z}}$$
and therefore finite dimensional.
This remark, combined with the previous proposition, 
shows the following result, needed later. A
 module map $T\in{\cal L}_{\text{Ind}({\cal 
F})}(H\otimes\text{Ind}({\cal F}), H'\otimes\text{Ind}({\cal F}))$
will be regarded as a function $T:G\to{\cal L}(H, H')\otimes{\cal F}$. 

\medskip

\noindent{\bf 11.5. Corollary} {\sl
There is a full and faithful $^*$--functor from the full subcategory of 
$\text{Mod}_{\rho}(G)$
with objects $v\otimes\rho$, $v\in\text{Rep}(G)$,
to the category   $\text{Mod}_\beta(K)$, given by
$$v\otimes\rho\to v\upharpoonright_K\otimes\beta, \quad
T\in(v\otimes\rho, v'\otimes\rho)\to T(1)\in(v\upharpoonright_K\otimes\beta, v'\upharpoonright_K\otimes\beta).$$ The inverse map on arrows is given by
$A\to A'$ with $A'(g):=v'(g)^*\otimes IAv(g)\otimes I.$
}\medskip

The functor $T\to T(1)$ defined in the above corollary will be referred to as the 
{\it evaluation functor}.

Given a   projection $E\in{\cal L}(H_v)\otimes{\cal F}$ and a 
unital $^*$--homomorphism $\eta:{\cal F}\to E{\cal L}(H_v)\otimes{\cal 
F}E$ 
defining a full bimodule structure on the intermediate $K$--module
$Y=EH_v\otimes{\cal F}$, i.e.\ satisfying conditions $(11.1)$--$(11.4)$, we 
may consider  the projection
 $\tilde{E}\in C(G,{\cal L}(H_v)\otimes{\cal F})\simeq {\cal 
L}(H_v)\otimes C(G, {\cal F})$,
$$\tilde{E}(g):=v(g)^*\otimes I E v(g)\otimes I,$$
which 
clearly satisfies 
$$\iota\otimes\beta_k(\tilde{E}(g))=$$
$$v(g)^*\otimes 
I\iota\otimes\beta_k(E) 
v(g)\otimes I=v(kg)^*\otimes I E v(kg)\otimes I=\tilde{E}(kg),$$
hence
 $\tilde{E}\in {\cal L}(H_v)\otimes\text{Ind}({\cal F})$.
 We may also consider the map taking a continuous function 
$f$ on $G$ with values in ${\cal F}$ 
to the function
$$\tilde{\eta}(f)(g):=v(g)^*\otimes I\eta(f(g))v(g)\otimes I.$$
Similar computations and $(11.3)$ show that if 
$f\in\text{Ind}({\cal F})$ then $\tilde{\eta}(f)\in{\cal 
L}(H_v)\otimes\text{Ind}({\cal F})$ and 
$\tilde{E}\tilde{\eta}(f)=\tilde{\eta}(f)=\tilde{\eta}(f)\tilde{E}$,
hence $\tilde{\eta}$ is  in fact a unital $^*$--homomorphism between
$$\tilde{\eta}: \text{Ind}({\cal 
F})\to \tilde{E}{\cal L}(H_v)\otimes\text{Ind}({\cal F})\tilde{E},$$
and $(\tilde{E},\tilde{\eta})$ defines a bimodule
over the induced algebra $\text{Ind}({\cal F})$.
We shall refer to it as the {\it induced 
bimodule}.
\medskip

\noindent{\bf 11.6.\ Theorem} {\sl The induced bimodule  
$(\tilde{E},\tilde{\eta})$ satisfies $(11.1)$--$(11.4)$ if $(E,\eta)$  does. 
Furthermore, if ${\cal F}$ is the completion of the dense $^*$--subalgebra of $K$--finite elements in the maximal $C^*$--norm, any intermediate projective $G$--module 
$X^\rho_v\subset Y\subset H_v\otimes\text{Ind}({\cal F})$ with a full 
bimodule 
structure is defined by such a pair $(E,\eta)$.}\medskip

\noindent{\it Proof} The validity of $(11.1)$--$(11.4)$ for  a  bimodule induced from 
one with analogous properties
follows easily from the previous proposition. Conversely, let $(E', \eta')$ satisfy 
$(11.1)$--$(11.4)$ with respect to the automorphism group $\rho$ of the induced algebra. By $(11.2)$, $v(g)\iota\otimes \rho_g(E')v(g)^*\otimes I=E'$. Evaluating in $g'$ gives
$v(g)\otimes IE'(g'g)v(g)^*\otimes I=E'(g')$, hence
$E'(g)=v(g)^*\otimes I E v(g)\otimes I$, where $E:=E'(1).$ It is now clear that $E$ satisfies $(11.1)$. Moreover, for $k\in K$,
$$v(k)^*\otimes IEv(k)\otimes I=E'(k)=\iota\otimes\beta_k(E'(1)),$$ hence
$E$ satisfies $(11.2)$.

On the other hand, $E'{\cal L}(H_v)\otimes\text{Ind}({\cal F})E'$, with $G$--action 
$\text{Ad} v\otimes\rho$, is  
isomorphic to the $C^*$--system induced by $E{\cal L}(H_v)\otimes{\cal 
F}E$ with $K$--action $\text{Ad} v\upharpoonright_K\otimes \beta$. An 
explicit $G$--equivariant isomorphism takes  $f\in E'{\cal L}(H_v)\otimes\text{Ind}({\cal F})E'$ to the element
of $C(G, E{\cal L}(H_v)\otimes{\cal F}E)$ defined by
$g\in G\to \text{Ad} v(g)\otimes I f(g).$
Therefore condition $(11.3)$ can be regarded as an intertwining relation 
between induced group representations. Hence,  by Frobenius reciprocity, there
is a map, a priori just linear, and densely defined on the $^*$--subalgebra of $K$--finite elements,
$\eta:{\cal F}\to E{\cal L}(H_v)\otimes{\cal F}E$ satisfying the 
intertwining relation
$$\eta(\beta_k(f))=\text{Ad}v_k\otimes\beta_k(\zeta(f)),$$
and hence $(11.3)$,
for $f\in{\cal F}$, $k\in K$, inducing $\eta'$ via
$$\eta'(f)(g)=\text{Ad}v(g)^*\otimes I\eta(f(g)).$$
We show that $\eta$ is a unital $^*$--homomorphism. It is well known that
for any $K$--finite element $f_1\in{\cal F}$, there is an element $f\in\text{Ind}({\cal F})$
with $f_1=f(1)$. Thus, $(11.4)$ follows. On the other hand, since $\eta'$ is a unital $^*$--homomorphism,  the above formula, 
 evaluated in $1$, shows that $\eta$ is a unital $^*$--homomorphism on the 
dense $^*$--subalgebra of $K$--finite elements. 
Since ${\cal F}$ is the completion in the maximal $C^*$--norm, we may 
conclude 
that the unique extension of $\eta$ to ${\cal F}$ has the required properties.
\medskip

\noindent{\it 11.4 Classification of intermediate full bimodule representations for type $I$ ergodic actions.}
Since a type $I$ ergodic action of a compact group $G$ is induced by an ergodic action $\beta$ of a closed subgroup $K$ on a matrix algebra  and since all intermediate full submodule
representations
for the ergodic action of $G$ are induced by similar submodules for 
the action of $K$ (Theorem 11.6), it suffices to classify the the intermediate full submodule representations for the action of the subgroup.

 Let the compact group $K$ act on
${\cal F}$. We first give a simple method of constructing 
extensions of $X^\beta_v$ to 
full 
bimodule representations on projective submodules of $H_v\otimes{\cal F}$. 
 \medskip

\noindent{\bf 11.7.\ Proposition} {\sl Pick a representation $v$
of $K$.

\begin{description}
\item{\rm a)}
If there is a unitary representation $z$ of $K$ with 
$$\text{dim}(z)\leq\text{dim}(v)-\text{dim}(\overline{L}_v)\eqno(11.5)$$ and an 
isometry
$$W\in(z\otimes\beta, v\otimes\beta)\quad 
\text{such that }\quad W^*Z_v=0,\eqno(11.6)$$  
then the intermediate $K$--module subrepresentation $X_v\subset Y\subset 
H_v\otimes{\cal F}$ 
defined by
the projection $E:=Z_vZ_v^*+WW^*\in{\cal L}(H_v)\otimes{\cal F}$ 
becomes a full bimodule 
$K$--representation  
with left action $\eta(f):=Z_vI\otimes f Z_v^*+W
I\otimes fW^*$, 

\item{\rm b)}
if we can choose $z$ with 
$\text{dim}(z)=\text{dim}(v)-\text{dim}(\overline{L}_v)$ then
we get a full $K$--bimodule representation  for
$v\otimes\beta$. 
\end{description}
}\medskip

\noindent{\it  Proof} 
$E$ and $\eta$ defined as in 
a) certainly satisfy the assumptions of Proposition 11.2, hence 
we get a full 
bimodule representation $Y$, and b) clearly follows.\medskip

We next provide 
a complete 
list if ${\cal F}$ is a matrix algebra.
\medskip

\noindent{\bf 11.9.\ Theorem} {\sl Let $\beta$ be an ergodic action of a  compact group $K$  
on a factor ${\cal F}$ and $v$ a representation of $K$.

\begin{description}
\item{\rm a)}
Two pairs $(z, W)$, $(z', W')$ satisfying $(11.5)$ and $(11.6)$
define the same intermediate $K$--bimodule representation $Y$ if and only if there is a unitary
intertwiner $U\in(z,z')$ such that $W=W'U\otimes I$,

\item{\rm b)}
if ${\cal F}$ is a matrix algebra, then 
 any full intermediate 
bimodule representation $X_v\subset Y\subset H_v\otimes {\cal F}$ arises from a pair 
$(z, W)$. In particular, full $K$--bimodule representations 
 on $H_v\otimes{\cal F}$
correspond to pairs $(z,W)$ satisfying $(11.5)$ and $(11.6)$ where the inequality of $(11.5)$ is strengthened to an equality.
\end{description}

}\medskip

\noindent{\it Proof} a) Obviousy  two equivalent pairs $(z,W)$, $(z', W')$, as in a), give rise to the same intermediate $K$--module $Y$ with the same left action $\eta$. Conversely, suppose $(z,W)$ and $(z',W')$ define the 
same $K$--bimodule representation $Y$. Then clearly $WW^*=W'{W'}^*$. Since the two left actions 
coincide, $WI\otimes fW^*=W'I\otimes f {W'}^*$. Hence the unitary 
${W'}^*W\in(z\otimes\beta, z'\otimes\beta)$ is of 
the form $U\otimes I$, with $U:H_z\to H_{z'}$, as ${\cal F}$ is a factor.
Thus $W=W'U\otimes I$. Making the intertwining property of $W$ and $W'$ explicit shows that $U\in(z, z')$.

b) Assume that ${\cal F}=\text{Mat}_r({\mathbb C})$. If $Y$ is 
defined by $E$ and $\eta$,
then $E$ needs to be 
of rank $qr$ with $q$ integer, as $\eta$ is unital, and 
$q\geq\text{dim}(\overline{L}_v)$ as $E\geq E_v$.
Set $\eta_1(f):=\eta(f)(E-E_v)$. We can write $\eta_1$ in the form
$\eta_1(f)=W I\otimes f W^*$ with $W$ a partial isometry such that
$WW^*=E-E_v=:E_1$ and $W^*W\in{\cal L}(H_v)\otimes{\mathbb C}$. The relation 
$W^*Z_v=0$ implies 
$\text{dim}(W^*WH_v)+\text{dim}(\overline{L}_v)\leq \text{dim}(v)$.
The covariance condition $(3.5)$ 
for $Y$ becomes 
$$\text{Ad}v(k)\otimes\beta_k\eta(f)=\eta(\beta_k(f)),\quad 
f\in\text{Mat}_r({\mathbb C}),$$
and is equivalent to requiring
an analogous relation for $\eta_1$:
$$v(k)\otimes I\iota\otimes\beta_k(W)I\otimes 
\beta_k(f)\iota\otimes\beta_k(W^*)v(k)^*\otimes 
I=WI\otimes\beta_k(f)W^*,$$
or
$$W^*v(k)\otimes I\iota\otimes\beta_k(W)\in
{\cal 
L}(H_v)\otimes{\mathbb C}.$$
On the other hand,
the  map
$k\to z(k)$ with $z(k)$ defined by 
$$z(k)\otimes I:=W^*v(k)\otimes I\iota\otimes\beta_k(W)$$ is
 a unitary representation of $K$ on the subspace $W^*WH_v$, completing
the proof of c).
\medskip

The following result summarizes  the classification 
of full bimodule representations for type $I$ ergodic 
actions  achieved here.
\medskip

\noindent{\bf 11.10.\ Theorem} {\sl Let ${\cal F}$ be a matrix 
algebra, and let $\beta$ be an ergodic action
of a closed subgroup $K$ of a compact group $G$ on ${\cal F}$. Pick a unitary f.d.\ 
representation $v$ of $G$. Then 
\begin{description}
\item{\rm a)}
the full bimodule $G$--representations over intermediate projective $G$--module subrepresentations
$X_v\subset Y \subset H_v\otimes\text{Ind}({\cal F})$ are classified by equivalence classes of pairs $(z,W)$, where  $z$ is a unitary f.d.\ representation of $K$ and $W\in(z\otimes\beta, v\upharpoonright_K\otimes\beta)$
an isometry satisfying 
$$W^*Z^\beta_{v\upharpoonright_K}=0.$$
$(W, z)$ and $(W', z')$ are equivalent if there is a unitary intertwiner  $U\in(z,z')$ 
with $W=W'U\otimes I$. 
\item{\rm b)}
In particular, the full $G$--bimodule representations 
on $H_v\otimes\text{Ind}({\cal F})$ correspond to 
pairs $(W,z)$ where $W$ satisfies
$$WW^*+Z^\beta_{v\upharpoonright_K}{Z^\beta_{v\upharpoonright_K}}^*=I.$$
The corresponding left module structure $\tilde{\eta}:\text{Ind}({\cal F})\to{\cal L}(H_v)\otimes\text{Ind}({\cal F})$  is given by 
$$\tilde{\eta}(f)(g)=v(g)^*\otimes I(Z^\beta_{v\upharpoonright_K} I\otimes 
f(g){Z^\beta_{v\upharpoonright_K}}^*+WI\otimes f(g)W^*)v(g)\otimes I.$$
\end{description}}\medskip

\noindent{\it Remark}
If
the module $G$--representation $v\otimes\rho$ over 
$\text{Ind}({\cal F})$ can be made into  a full bimodule 
$G$--representation,  
and if it is induced by the pair $(z,W)$,  
we may form the $K$--representation $z':=\iota_{\overline{L}_v}\oplus z$ 
of the same dimension as $v$. Then $Z^\beta_{v\upharpoonright_K}\oplus W$ 
is a 
unitary equivalence  from the original full bimodule structure 
for
$v\upharpoonright_K\otimes\beta$ inducing the given full bimodule 
structure for $v\otimes\rho$, in the sense of Theorem 11.6, 
to $z'\otimes\beta$
with the {\it trivial} left module structure. This remark will play a role 
in the proof of Theorem 6.7. 
\medskip

As a consequence of b) of Theorem 11.9, the module representation
$v\otimes\beta$, with $v$ in the  spectrum,  in some cases,
does not admit any full bimodule $K$--representation unless $v$ has full
multiplicity. We discuss a class of examples.
\medskip

\noindent{\bf 11.11. Example} Consider the adjoint action $\beta_r$ of the unique
$r+1$--dimensional irreducible representation $v_r$ of the group $K=SU(2)$ 
acting on the matrix algebra
$\text{Mat}_{r+1}({\mathbb C})$. We show that if $r\geq1$, $v\otimes\beta_r$ becomes a 
full bimodule representation  only for certain $v$. Hence none of the 
actions $\beta_r$ arise from a relaxed tensor functor 
$\text{Rep}(SU(2))\to{\cal M}$ to a tensor $C^*$--category, as 
this functor would make
all $v\otimes\beta_r$ 
into full bimodule representations
by Theorem 6.2. 

The spectrum of $\beta_r$ may be determined by the Clebsch--Gordan rule
$$v_r\otimes v_s\simeq v_{r-s}\oplus v_{r-s+2}\oplus\dots\oplus v_{r+s},\quad r\geq s,$$
after regarding $\beta_r$ as a    Hilbert space representation with respect to the inner product defined by the ($K$--invariant) trace of $\text{Mat}_{r+1}({\mathbb C})$.  $v_r$  being  selfconjugate, we have
$\beta_r\simeq v_r\otimes v_r\simeq v_0\oplus v_2\oplus\dots\oplus v_{2r}$.
Hence any spectral representation has multiplicity $1$. 

In particular, $v_1$  is never in the spectrum of $\beta_r$, and the full 
bimodule structures on $H_{v_1}\otimes\text{Mat}_{r+1}({\mathbb C})$ are described by pairs $(z,W)$ with $\text{dim}(z)=\text{dim}(v_1)=2$.
Since $z$ can never contain the trivial representation, we necessarly have $z=v_1$.
Hence we need to specify a unitary
$$W\in(v_1\otimes\beta_r, 
v_1\otimes\beta_r)\simeq$$
$$(v_1\otimes v_r, v_1\otimes v_r)\simeq(v_{r-1}\oplus v_{r+1}, v_{r-1}\oplus v_{r+1})\simeq{\mathbb C}\oplus{\mathbb C}.$$
Hence $v_1\otimes\beta_r$ admits full bimodule 
structures, and they are classified by ${\mathbb T}$. 

On the other hand, low multiplicity of a representation in the spectrum in general rules out full bimodule structures on $v\otimes\beta_r$ as the following simple argument shows.
If there were a structure of a full $K$--bimodule representation  on 
$v_2\otimes\beta_r$
defined by $(z, W)$, then we must have 
$\text{dim}(z)=\text{dim}(v_2)-\text{mult}(v_2)=2$. Since $z$ cannot 
contain the trivial representation, $z=v_1$. On the other hand 
the space of module intertwiners 
$(v_1\otimes\beta_r, 
v_2\otimes\beta_r)$   is isomorphic to
$(v_1\otimes v_r, v_2\otimes v_r)$ which is trivial, again by the Clebsch--Gordan rule.
Hence $v_2\otimes\beta_r$ 
admits no 
full bimodule structure, and actually $X_{v_2}$ admits no proper extension to a full bimodule representation.

\medskip

\noindent{\bf 11.12.\ Corollary} {\sl Let $K$ act ergodically on a 
matrix 
algebra.
\begin{description}
\item{\rm a)}
 Let $v$ be a  representation with $\overline{L}_v\neq\{0\}$ and assume 
that
 any 
irreducible of smaller dimension has full multiplicity. If 
$\text{dim}(\overline{L}_v)<\text{dim}(v)$ 
then
$X_v$ does not admit any proper 
extension to a full bimodule 
$K$--representation.

\item{\rm b)} If $\beta$ has full spectrum (hence $K$ is finite) and each 
$v\otimes\beta$ can be made into a full bimodule 
representation 
then each irreducible is of full multiplicity in $\beta$.
\end{description}

}\medskip

\noindent{\it Proof} Let $z$ and 
$W$ be as required in a) of Prop.\ 11.7. Since
$\text{dim}(z)<\text{dim}(v)$, any irreducible 
subrepresentation of $z$ has full multiplicity in $\beta$. Hence there is 
a 
unitary $U\in{\cal L}(H_z)\otimes\text{Mat}_r({\mathbb C})$
 with $\iota\otimes\beta_k(U)=z(k)^*\otimes IU$. Hence every column 
of 
$WU$ gives an element of $\overline{L}_v$ orthogonal to 
$\overline{L}_v$ itself, as $W^*Z_v=0$. So $WU=0$ and $W=0$. This completes the 
proof of a) and b) follows easily.

\medskip

\section{Tensorial properties of the evaluation functor}

In   this section we use the classification of full Hilbert bimodule 
structures on type $I$ von Neumann algebras obtained in the previous section to prove Theorem  6.7, and Corollaries 6.8 and 6.9.

We need a few simple lemmas that clarify the tensorial properties of the evaluation  functor defined in the previous section.
We thus assume that we are given an action $\beta$ of a closed subgroup 
$K$ of a compact group $G$ on a $C^*$--algebra  ${\cal F}$ and that for 
each $v\in\text{Rep}(G)$ we have a full bimodule structure for
$v\upharpoonright_K\otimes\beta$ defined by the $^*$--homomorphism 
$\eta_v:{\cal F}\to {\cal L}(H_v)\otimes{\cal F}$. We consider the full bimodule structure $\tilde{\eta}_v$ for $v\otimes\rho$ induced by $\eta_v$ as in subsect.\ 11.3.
\medskip

\noindent{\bf 12.1. Lemma} {\sl If $T\in(v\otimes\rho, v'\otimes\rho)$ 
is a bimodule map 
 then
$T(1)\in(v\upharpoonright_K\otimes\beta, 
v'\upharpoonright_K\otimes\beta)$ is a bimodule map as well.}\medskip

\noindent{\it Proof} The proof is straightforward. By Cor.\ 11.5, we may 
write $T$ in the form $T(g)= v'(g)^*\otimes I T(1)v(g)\otimes I$, with 
$T(1)\in(v\upharpoonright_K\otimes\beta, v'\upharpoonright_K\otimes\beta)$. The intertwining relation for $T$ evaluated at $1$ gives the  intertwining relation for $T(1)$. \medskip

Let us now consider a unital $C^*$--algebra ${\cal C}$ and two f.d.\ Hilbert spaces $H$ and $L$. Consider the right $C^*$--modules $H\otimes{\cal C}$ and $L\otimes{\cal C}$.
If $L\otimes{\cal C}$ also has a left ${\cal C}$--module structure defined by 
a unital $^*$--homomorphism $\eta:{\cal C}\to{\cal L}_{\cal C}(L)\otimes{\cal C}$
then we may form the tensor product right Hilbert $C^*$--module $(H\otimes{\cal C})\otimes_{\cal C} (L\otimes{\cal C})$, to be identified with $(H\otimes L)\otimes{\cal C}$. We may thus form tensor products $T\otimes S$ of a module intertwiner $T\in{\cal L}(H, H')\otimes{\cal C}$ with a bimodule intertwiner $S\in{}_{\cal C}{\cal L}_{\cal C}(L\otimes{\cal C},
 L'\otimes{\cal C})$ giving an element of ${\cal L}(H\otimes L, H'\otimes L')\otimes{\cal C}$.
 
 \medskip

\noindent{\bf 12.2.\  Lemma} {\sl Let us consider
 $H_v\otimes\text{Ind}({\cal F})$ and $H_{v'}\otimes\text{Ind}({\cal F})$
 as right $\text{Ind}({\cal F})$--modules. Let 
$\tilde{\eta}_u$, $\tilde{\eta}_{u'}$ 
make 
$H_{u}\otimes\text{Ind}({\cal F})$ and 
$H_{u'}\otimes\text{Ind}({\cal F})$ into  $\text{Ind}({\cal 
F})$--bimodules. 
For a module intertwiner $T\in(v\otimes\rho, v'\otimes\rho)$
and a bimodule intertwiner $S\in (u\otimes\rho, u'\otimes\rho)$,
we have 
$$(T\otimes S)(1)=T(1)\otimes S(1).$$}\medskip

\noindent{\it Proof} Notice that $S(1)$ is a bimodule intertwiner by the previous lemma, hence the right hand side makes sense.  Let $H$, $H'$, $L$, $L'$ be  f.d.\ Hilbert spaces and $\eta$, $\eta'$   left ${\cal C}$--module structures on  $L\otimes{\cal C}$ and $L'\otimes{\cal C}$ respectively.
Given a module intertwiner $T\in{\cal L}(H, H')\otimes{\cal C}$, and a bimodule intertwiner $S\in{}_{\cal C}{\cal L}_{\cal C}(L\otimes{\cal C},
 L'\otimes{\cal C})$, a simple computation shows that if $T$ is represented by the ${\cal C}$--valued  matrix $(t_{rs})$, in the sense 
 that $T=\sum_{rs} e_{rs}\otimes t_{rs}$, where $e_{rs}$ are matrix units, and if $S$ is represented  by $(s_{pq})$ then the module intertwiner
 $T\otimes S$ regarded as an element of ${\cal L}(H\otimes L, H'\otimes L')\otimes{\cal C}$ is represented by the matrix whose $(rp)(sq)$-entry is 
 $\sum_h\eta'(t_{rs})_{ph}s_{hq}$. 
 We apply this to $H_v$, $H_{v'}$, $H_u$, $H_{u'}$ and $\text{Ind}({\cal F})$ respectively.
 By Cor.\ 11.5, we may write $(t_{rs})(g)=v'(g)^*T(1)v(g)$, $(s_{pq})(g)=u'(g)^*S(1)u(g)$, where $T(1)$ and $S(1)$ are now represented by ${\cal F}$--valued matrices.
Recalling how $\tilde{\eta}_{u'}$ was defined before Theorem 11.6, the
 $(rp)(sq)$-entry of $T\otimes S$ is the 
 function 
 $$\sum_h\tilde{\eta}_{u'}(t_{rs})_{ph}(g)s_{hq}(g)=\sum_{h,l,m}\overline{u'(g)}_{lp}\eta_{u'}(t_{rs}(g))_{lm}u'(g)_{mh}s_{hq}(g)=$$
 $$\sum_{h,l,m,i,j,k,t}\overline{u'(g)}_{lp}\overline{v'(g)}_{ir}\eta_{u'}(T(1)_{ij})_{lm}v(g)_{js}u'(g)_{mh}\overline{u'(g)}_{kh}S(1)_{kt}u(g)_{tq}=$$
 $$\sum_{l,i,j,k,t}\overline{u'(g)}_{lp}\overline{v'(g)}_{ir}\eta_{u'}(T(1)_{ij})_{lk}v(g)_{js}S(1)_{kt}u(g)_{tq}=$$
 $$\sum_{l,i,j,k,t} v'\otimes u'(g)^*_{(rp)(il)}\eta_{u'}(T(1)_{ij})_{lk}S(1)_{kt}v\otimes u(g)_{(jt)(sq)}=$$
 $$\sum_{l,i,j,k,t} v'\otimes u'(g)^*_{(rp)(il)}(T(1)\otimes S(1))_{(il)(jt)}
 v\otimes u(g)_{(jt)(sq)}.$$
 Hence $(T\otimes S)(1)=T(1)\otimes S(1).$
 \medskip

 Note that if $H\otimes{\cal C}$ and $L\otimes{\cal C}$ 
have left bimodule structures defined by 
$\eta:{\cal C}\to{\cal L}(H)\otimes{\cal C}$ and 
$\zeta:{\cal C}\to{\cal L}(L)\otimes{\cal C}$ then 
under the  unitary module map 
$(H\otimes{\cal C})\otimes_{\cal C}(L\otimes{\cal C})\simeq 
(H\otimes L)\otimes{\cal C}$ the left module structure 
${\cal C}\to{\cal L}(H\otimes L)\otimes{\cal C}$ corresponding to 
the tensor product bimodule is given by $\iota_{{\cal 
L}(H)}\otimes\zeta\circ\eta$, $\iota_{{\cal L}(H)}$ being the identity map on ${\cal L}(H)$. 
This tensor product left action will be denoted 
by $\eta\otimes\zeta$.
\medskip

\noindent{\bf 12.3.\ Lemma} {\sl If the induced set of left actions $\{\tilde{\eta}_u, u\in\text{Rep}(G)\}$ on the $C^*$--modules $H_u\otimes\text{Ind}({\cal F})$ is tensorial, i.e. $\tilde{\eta}_{u\otimes v}=\tilde{\eta}_u\otimes\tilde{\eta}_v$ for $u,v\in\text{Rep}(G)$ then  the original set  $\{\eta_u, u\in\text{Rep}(G)\}$ is tensorial too.}\medskip

\noindent{\it Proof} It suffices to evaluate  the tensorial relation for the $\tilde{\eta}_u$'s at $1$.
\medskip

We summarize the above lemmas as follows.
\medskip

\noindent{\bf 12.4.\ Theorem} {\sl 
Let $\beta$ be an action of a closed subgroup 
$K$ of a compact group $G$ on a $C^*$--algebra  ${\cal F}$. Assume that   
for each $v\in\text{Rep}(G)$ we have a full bimodule structure for 
$v\upharpoonright_K\otimes\beta$ defined by the $^*$--homomorphism $\eta_v:{\cal F}\to {\cal L}(H_v)\otimes{\cal F}$. If the set of induced bimodule structures $\tilde{\eta}_v$ for $v\otimes\rho$  is tensorial then the evaluation functor $T\to T(1)$  restricts to a faithful tensor functor from the full tensor  $C^*$--subcategory of $\text{Bimod}_\rho(G)$ with objects $v\otimes\rho$ to $\text{Bimod}_\beta(K)$. }\medskip

\medskip

\noindent{\it Proof of Theorem  6.7. and Cor.\ 6.8.}
 Theorem 6.2, applied 
to the given tensor functor $\mu:{\cal S}_G\to{\cal M}$ and to
the embedding functor $\tau:{\cal S}_G\to\text{Hilb}$, allows us to 
identify ${\cal 
M}_\mu^\otimes$ with the full subcategory of $\text{Bimod}_\alpha(G)$ with objects $u^r\otimes\alpha$, $r=0,1,2,\dots$, where $u$ is the distinguished representation of $G$ and, as
before, $\alpha$ is the ergodic action of $G$ on the associated $C^*$--algebra 
${\cal C}$.
That theorem provides us with a full $G$--bimodule representation  for each 
$u^r\otimes\alpha$ and the collection of these left module structures is tensorial.
Since the von Neumann completion of ${\cal 
C}$ in the GNS representation of the $G$--invariant trace state is of type
$I$, we may identify the completed ergodic system with a  von 
Neumann ergodic system  $(\text{Ind}({\cal F}),\rho)$ induced from a closed 
subgroup $K$, unique up to conjugation, where ${\cal F}$ is a 
matrix algebra with an ergodic  action $\beta$ of $K$.
The left ${\cal C}$--action on $H_{u^r}\otimes{\cal C}$ is defined by a 
unital $^*$--homomorphism $\eta:{\cal C}\to{\cal L}(H_u)\otimes{\cal C}$
intertwining $\alpha$ with $\text{Ad}(u)\otimes\alpha$ by Proposition 11.2.
Hence, if $\text{tr}$ and $\tau$ are  the normalized $G$--invariant
traces on ${\cal L}(H_u)$ and ${\cal C}$ respectively, 
$(\text{tr}\otimes\tau)\circ\eta$
 is a $G$--invariant trace on ${\cal C}$. Such a trace is unique so 
$(\text{tr}\otimes\tau)\circ\eta=\tau$. Thus $\eta$ induces a normal
$^*$--homomorphism from $\text{Ind}({\cal F})$ to ${\cal L}(H_u)\otimes\text{Ind}({\cal F})$.
Correspondingly, we get a set of tensorial full bimodule structures 
for $u^r\otimes\rho$. Thus by Theorem 12.4 there is  
a faithful tensor functor from the full 
subcategory of $\text{Bimod}_\rho(G)$ with objects $u^r\otimes\rho$
to the full subcategory ${\cal T}$ of $\text{Bimod}_\beta(K)$ with objects 
$u^r\upharpoonright_K\otimes\beta$. 
We next apply Theorem 11.10 to the full 
bimodule $K$--representation $u\upharpoonright_K\otimes\beta$ 
fixing a pair $(z, W)$.
We set $z':=\iota_{\overline{L}^\beta_{u\upharpoonright_K}}\oplus {z}$ and
$U:=Z_{u\upharpoonright_K}\oplus W$, a $K$--bimodule unitary in $(z'\otimes\beta, u\upharpoonright_K\otimes\beta)$ if $z'$ has the trivial left ${\cal C}$--action.
We define a $^*$--functor ${\cal T}\to\text{Rep}(K)$ taking $u^r\upharpoonright_K\otimes\beta$ to ${z'}^r$ and a bimodule  intertwiner $T\in(u^r\upharpoonright_K\otimes\beta, u^s\upharpoonright_K\otimes\beta)$ to ${U^*}^{\otimes s}TU^{\otimes r}$, which is 
tensorial to the category of Hilbert bimodule representations.
We need to show that any arrow is in fact an arrow in the category $\text{Rep}(K)$ regarded as embedded into the category of bimodule representations as a tensor $C^*$--category. In other words, we need to show that ${U^*}^{\otimes s}TU^{\otimes r}$ lies in the subspace ${\cal L}(H_{{z'}^r}, H_{{z'}^s})\otimes{\mathbb C}$ of ${\cal L}(H_{{z'}^r}, H_{{z'}^s})\otimes{\cal F}$.
To this end, recall that Theorem 6.2 ensures that any module
 $G$--intertwiner is in fact a bimodule intertwiner, see   
Theorem 9.6. 
The same property holds 
for the bimodule structures of the $u^r\upharpoonright_K\otimes\beta$'s and hence for the bimodule structures of the ${z'}^r\otimes\beta$, unitarily related to them, since the evaluation functor is full and faithful, see Cor.\ 11.5.
But now each  ${z'}^r\otimes\beta$ has the trivial left module structure over ${\cal F}$,
hence a bimodule intertwiner lies in $ {\cal L}(H_{{z'}^r}, H_{{z'}^s})\otimes Z({\cal F})={\cal L}(H_{{z'}^r}, H_{{z'}^s})\otimes{\mathbb C}$ since ${\cal F}$ is a factor, see the discussion following Prop.\ 5.1.
This argument completes the proof of Theorem  6.7. If in particular ${\cal C}$ is commutative then ${\cal F}={\mathbb C}$, and $z'=u\upharpoonright_K$,
completing the proof of Cor.\ 6.8.
\medskip

\noindent{\it Proof of Cor.\ 6.9} The condition on $R$ allows us to define 
a tensor functor from ${\cal S}_{SU(2)}$ to ${\cal M}$ taking the defining representation $u$ to $x$ and the determinant element to $R$, see \cite{DRduals}. We may now apply
Theorem  6.7.
\medskip

\section{Appendix}

In this appendix we collect some computations with quasitensor functors that we have used throughout the paper.\medskip

\noindent
{\bf 13.1.\ Proposition} {\sl If we take $1_{\mu_{\overline v}}\otimes \hat R_u\otimes 1_{\mu_v}\circ \hat R_v$ as a solution of the conjugate equations for $\mu_u\otimes\mu_v$ and $\hat R_{u\otimes v}$ as the solution for $\mu_{u\otimes v}$, where 
$\hat R_{u\otimes v}$ is the image solution of the tensor product solution for $u\otimes v$, 
then 
$$\tilde\mu_{u,v}^\bullet =\tilde\mu_{\overline v,\overline u},\quad 
\tilde\mu_{u,v}^{*\bullet}=\tilde\mu_{\overline v,\overline 
u}^*.$$}\medskip

\noindent
{\it Proof} 
$$\tilde\mu_{u,v}^\bullet=(\hat R_v^*\circ 1_{\mu_{\overline v}}\otimes\hat R_u^*\otimes 1_{\mu_v})\otimes 1_{\mu_{\overline v\otimes\overline u}}\circ 1_{\mu_{\overline v}}\otimes 1_{\mu_{\overline u}}\otimes\tilde\mu_{u,v}^*\otimes 1_{\mu_{\overline v\otimes\overline u}}\circ 1_{\mu_{\overline v}}\otimes 1_{\mu_{\overline u}}\otimes\hat{\overline R}_{u\otimes v}$$
$$=(\hat R_v^*\circ 1_{\mu_{\overline v}}\otimes\hat R_u^*\otimes 1_{\mu_v})
\otimes 1_{\mu_{\overline v\otimes\overline u}}\circ 1_{\mu_{\overline v}}\otimes 1_{\mu_{\overline u}}\otimes(\tilde\mu_{u,v}^*\otimes 1_{\mu_{\overline v\otimes\overline u}}\circ\tilde\mu_{u\otimes v,\overline v\otimes\overline u}^*\circ\mu(\overline R_{u\otimes v})=$$
$$(\hat R_v^*\circ 1_{\mu_{\overline v}}\otimes\hat R_u^*\otimes 1_{\mu_v})
\otimes 1_{\mu_{\overline v\otimes\overline u}}\circ 1_{\mu_{\overline v}}\otimes 1_{\mu_{\overline u}}\otimes(1_{\mu_u}\otimes\tilde\mu_{v,\overline v\otimes\overline u}^*\circ\tilde\mu_{u,v\otimes\overline v\otimes\overline u}^*\circ\mu(\overline R_{u\otimes v}))=$$
$$(\hat R_v^*\circ 1_{\mu_{\overline v}}\otimes\hat R_u^*\otimes 1_{\mu_v})
\otimes 1_{\mu_{\overline v\otimes\overline u}}\circ 
1_{\mu_{\overline v}}\otimes 1_{\mu_{\overline u}}\otimes
(1_{\mu_u}\otimes\tilde\mu_{v,\overline v\otimes\overline u}^*\circ 
1_{\mu_u}\otimes\mu(\overline R_v\otimes 1_{\overline u})
\circ\tilde\mu_{u,\overline u}^*\circ\mu(\overline R_u))=$$ 
$$(\hat R_v^*\circ 1_{\mu_{\overline v}}\otimes\mu(R_u^*)\otimes 
1_{\mu_v})\otimes 1_{\mu_{\overline v\otimes\overline u}}\circ 1_{\mu_{\overline v}}\otimes 1_{\mu_{\overline u\otimes u}}\otimes(\tilde\mu_{v,\overline v\otimes\overline u}^*\circ\mu(\overline R_v\otimes 1_{\overline u}))$$
$$\circ 1_{\mu_{\overline v}}\otimes\tilde\mu_{\overline u,u}\otimes 1_{\mu_{\overline u}}\circ1_{\mu_{\overline v}}\otimes 1_{\mu_{\overline u}}\otimes\hat R_u=$$
$$\hat R_v^*\otimes 1_{\mu_{\overline v\otimes\overline u}}\circ 
1_{\mu_{\overline v}}\otimes\tilde
\mu_{v,\overline v\otimes\overline u}^*\circ 1_{\mu_{\overline 
v}}\otimes\mu(\overline R_v\otimes 1_{\overline u})\circ 1_{\mu_{\overline v}}\otimes\mu(R_u^*)\otimes 1_{\mu_{\overline u}}\circ 1_{\mu_{\overline v}}\otimes\tilde\mu_{\overline u,u}^*\otimes 1_{\mu_{\overline u}}\circ 1_{\mu_{\overline v}}\otimes 1_{\mu_{\overline u}}\otimes\hat{\overline R}_u.$$ 
Now 
$$\mu(R_u^*)\otimes 1_{\mu_{\overline u}}\circ\tilde\mu_{\overline u,u}\otimes 1_{\mu_{\overline u}}\circ 1_{\mu_{\overline u}}\otimes\hat{\overline R}_u=$$
$$\mu(R_u^*)\otimes 1_{\mu_{\overline u}}\circ\tilde\mu_{\overline u,u}\otimes 1_
{\mu_{\overline u}}\circ 1_{\mu_{\overline u}}\otimes\tilde\mu_{u,\overline u}^*\circ 1_{\mu_{\overline u}}\otimes\mu(\overline R_u)=$$
$$\mu(R_u^*)\otimes 1_{\mu_{\overline u}}\circ\tilde\mu_{\overline u\otimes u,\overline u}^*\circ\tilde\mu_{\overline u,u\otimes u\otimes\overline u}\circ 1_{\mu_{\overline u}}\otimes\mu(\overline R_u)=$$ 
$$\mu(R_u^*\otimes 1_{\overline u})\circ\mu(1_{\overline u}\otimes\overline R_u)= 1_{\mu_u}.$$ 
Substituting this into our calculation gives 
$$\tilde\mu_{u,v}^\bullet=\hat R_v^*\otimes 1_{\mu_{\overline v\otimes\overline u}}\circ 1_{\mu_{\overline v}}\otimes\tilde\mu_{v,\overline v\otimes\overline u}^*\circ 1_{\mu_{\overline v}}\otimes\mu(\overline R_v\otimes 1_{\overline u})=$$
$$\mu(R_v^*)\otimes 1_{\mu_{\overline v\otimes\overline u}}\circ\tilde\mu_{\overline v,v}\otimes 1_{\mu_{\overline v\otimes\overline u}}\circ 1_{\mu_{\overline v}}\otimes\tilde\mu_{v,\overline v\otimes\overline u}^*\circ 1_{\mu_{\overline v}}\otimes\mu(\overline R_v\otimes 1_{\overline u})=$$ 
$$\mu(R_v^*)\otimes 1_{\mu_{\overline v\otimes\overline u}}\circ\tilde\mu_{\overline v\otimes v,\overline v\otimes u}^*\circ\tilde\mu_{\overline v,v\otimes \overline v\otimes\overline u}\circ 1_{\mu_{\overline v}}\otimes\mu(\overline R_v\otimes 1_{\overline u})=$$
$$\mu(R_v^*\otimes 1_{\overline v\otimes u})\circ\mu(1_{\overline v}\otimes\overline R_v\otimes 1_u)\circ\tilde\mu_{\overline v,\overline u}=\tilde\mu_{\overline v,\overline u}.$$ 
Dualizing with respect to $\otimes$ yields 
$\tilde\mu_{u,v}^{*\bullet*}=\tilde\mu_{\overline v,\overline u}$ and 
taking adjoints completes the proof.\medskip

\noindent{\bf 13.2.\ Corollary} {\sl For $M\in(\mu_u,\mu_{u'})$, 
$N\in(\mu_v,\mu_{v'}),$ 
$$(\tilde{\mu}_{u',v'}\circ M\otimes 
N\circ\tilde{\mu}_{u,v}^*)^\bullet=
\tilde{\mu}_{\overline{v'},\overline{u'}}\circ N^\bullet\otimes 
M^\bullet\circ\tilde{\mu}^*_{\overline{v},\overline{u}},$$
with respect to the image of a tensor product 
solution of the conjugate equations.}\medskip

\noindent{\it Proof} By the previous proposition,
$$(\tilde{\mu}_{u',v'}\circ M\otimes 
N\circ\tilde{\mu}_{u,v}^*)^\bullet=
\tilde{\mu}_{u',v'}^\bullet\circ (M\otimes N)^\bullet
\circ{\tilde\mu}_{u,v}^{*\bullet}=\tilde{\mu}_{\overline{v'},\overline{u'}}\circ N^\bullet\otimes 
M^\bullet\circ \tilde{\mu}_{\overline{v},\overline{u}}.$$
\medskip

\noindent{\bf 13.3.\ Proposition} {\sl If we take the conjugate solution ${R}_{\overline u}={\overline R}_u$ as a solution of the conjugate equations for $\overline{u}$ and
the tensor product solution
$R_{\overline{u}\otimes u}=1_{\overline{u}}\otimes R_{\overline u}\otimes 1_u\circ R_u$
for $\overline{u}\otimes u$ then $R_u^\bullet=R_u$ and $R_u^{*\bullet}=R_u^*$.}\medskip

\noindent{\it Proof}
$$R_u^{\bullet}= R_u^*\otimes 1_{\overline{u}\otimes u}\circ\overline{R}_{\overline{u}\otimes u}=$$
$$R_u^*\otimes 1_{\overline{u}\otimes u}\circ 1_{\overline{u}}\otimes\overline{R}_u\otimes 1_u\circ R_u=R_u.$$
Dualizing again with respect to $\otimes$ gives $R_u^{*\bullet}=R_u^*.$

\medskip

\noindent{\it Acknowledgements} We are grateful to S.\ Doplicher for numerous
discussions on an early stage of this line of research.  He pursued the idea of embedding 
tensor $C^*$--categories into categories of Hilbert bimodules.

Part of the results of this paper have been announced at a conference held 
in Leuven in the fall 2008. C.P.\ would like to thank S.\ Vaes for the invitation and for
discussions. She would also like to thank P.\ Hajac and R.\ Tomatsu for 
discussions.

Finally, we are grateful to the referee. His comments led to a simplified 
and improved presentation 
of  the paper.

\end{document}